\documentclass[11pt]{amsart}
\usepackage{amssymb,amsmath,a4wide,graphicx,stmaryrd,color}
\usepackage[all]{xy}
\usepackage[pagebackref=true]{hyperref}

\title{Combinatorics of labelling in higher dimensional automata} 
\author[P. Gaucher]{Philippe Gaucher}
\address{Laboratoire PPS  (CNRS UMR 7126)\\ Universit{\'e} Paris 7--Denis Diderot\\
  Site Chevaleret\\ Case 7014\\ 75205 PARIS Cedex 13 \\ France}
\urladdr{http://www.pps.jussieu.fr/{\~{}}gaucher/} 
\thanks{This work has been supported by the ANR
  ``Invariants alg\'ebriques des syst\`emes informatiques''
  ANR-05-BLAN-0267.}
\thanks{I thank very much Ronnie Brown, Philippe
Malbos and Tim Porter for helpful conversations.}

\swapnumbers

\newcommand{\C}{\mathcal{C}}
\newcommand{\D}{\mathcal{D}}

\newcommand{\N}{\mathbb{N}}

\newcommand{\de}{\partial}
\newcommand{\p}\times
\renewcommand{\vec}{\overrightarrow}
\renewcommand{\P}{\mathbb{P}}

\newcounter{numerothm}[section]

\numberwithin{section}{part}

\newtheorem{thm}[numerothm]{Theorem}
\newtheorem{prop}[numerothm]{Proposition}
\newtheorem{conj}[numerothm]{Conjecture}
\newtheorem{lem}[numerothm]{Lemma}
\newtheorem{cor}[numerothm]{Corollary}

\newtheorem{defn}[numerothm]{Definition}
\newtheorem{nota}[numerothm]{Notation}

\newcommand{\bd}{\begin{defn}}
\newcommand{\ed}{\end{defn}}
\newcommand{\bp}{\begin{prop}}
\newcommand{\ep}{\end{prop}}
\newcommand{\bth}{\begin{thm}}
\renewcommand{\eth}{\end{thm}}
\newcommand{\bpf}{\begin{proof}}
\newcommand{\epf}{\end{proof}}

\newcommand{\fl}[1]{\ar@{->}[ll]_-{#1}}
\newcommand{\fr}[1]{\ar@{->}[rr]^-{#1}}
\newcommand{\fd}[1]{\ar@{->}[dd]_-{#1}}
\newcommand{\fu}[1]{\ar@{->}[uu]^-{#1}}
\newcommand{\f}[2]{\ar@{->}[#1]|{#2}}
\newcommand{\ff}[2]{\ar@2{->}[#1]|{#2}}
\newcommand{\frr}[1]{\ar@{->}[rrrr]^-{#1}}

\renewcommand{\top}{{\mathbf{Top}}}

\newcommand{\iso}{\cong}
\newcommand{\lp}{\left(}
\newcommand{\rp}{\right)}
\newcommand{\ot}{\otimes}

\renewcommand{\leq}{\leqslant}
\renewcommand{\geq}{\geqslant}

\def\cartesien{%
  \ar@{-}[]+R+<6pt,-2pt>;[]+RD+<6pt,-6pt>%
  \ar@{-}[]+D+<2pt,-6pt>;[]+RD+<6pt,-6pt>%
}
\def\cocartesien{%
  \ar@{-}[]+L+<-6pt,+2pt>;[]+LU+<-6pt,+6pt>%
  \ar@{-}[]+U+<-2pt,+6pt>;[]+LU+<-6pt,+6pt>%
}
\def\hocartesien{%
  \ar@{-}[]+R+<6pt,-2pt>;[]+RD+<6pt,-6pt>_{h}%
  \ar@{-}[]+D+<2pt,-6pt>;[]+RD+<6pt,-6pt>%
}
\def\hococartesien{%
  \ar@{-}[]+L+<-6pt,+2pt>;[]+LU+<-6pt,+6pt>_{h}%
  \ar@{-}[]+U+<-2pt,+6pt>;[]+LU+<-6pt,+6pt>%
}

\newcommand{\brm}[1]{\rm{\mathbf{#1}}}

\newcommand{\dtop}{{\brm{Flow}}}
\newcommand{\set}{{\brm{Set}}}
\newcommand{\poset}{{\brm{PoSet}}}

\newcommand{\proc}{{\brm{Proc}}}

\DeclareMathOperator{\rec}{rec}
\DeclareMathOperator{\id}{Id}
\DeclareMathOperator{\sh}{Sh}

\DeclareMathOperator{\card}{card}

\newcommand{\liminj}{\varinjlim}

\newcommand{\cat}{{\mathbf{Cat}}}

\makeatletter
\def\varholim@#1#2{%
  \vtop{\m@th\ialign{##\cr
    \hfil$#1\operator@font holim$\hfil\cr
    \noalign{\nointerlineskip\kern1.5\ex@}#2\cr
    \noalign{\nointerlineskip\kern-\ex@}\cr}}%
}
\def\holimproj{%
  \mathop{\mathpalette\varholim@{\leftarrowfill@\textstyle}}\nmlimits@
}
\def\holiminj{%
  \mathop{\mathpalette\varholim@{\rightarrowfill@\textstyle}}\nmlimits@
}
\makeatother

\DeclareMathOperator{\cosk}{cosk}
\DeclareMathOperator{\COSK}{\vec{\cosk}}

\newcommand{\ddownarrow}{{\downarrow}}

\begin{document}

\begin{abstract}
  The main idea for interpreting concurrent processes as labelled
  precubical sets is that a given set of $n$ actions running
  concurrently must be assembled to a labelled $n$-cube, in exactly
  one way.  The main ingredient is the non-functorial construction
  called labelled directed coskeleton. It is defined as a subobject of
  the labelled coskeleton, the latter coinciding in the unlabelled
  case with the right adjoint to the truncation functor. This
  non-functorial construction is necessary since the labelled
  coskeleton functor of the category of labelled precubical sets does
  not fulfil the above requirement.  We prove in this paper that it is
  possible to force the labelled coskeleton functor to be well-behaved
  by working with labelled transverse symmetric precubical
  sets. Moreover, we prove that this solution is the only one. A
  transverse symmetric precubical set is a precubical set equipped
  with symmetry maps and with a new kind of degeneracy map called
  transverse degeneracy. Finally, we also prove that the two settings
  are equivalent from a directed algebraic topological viewpoint. To
  illustrate, a new semantics of CCS, equivalent to the old one, is
  given.
\end{abstract}

\maketitle

\tableofcontents

\part{Introduction}

\section{Presentation of the results}

Directed algebraic topology is a field of research aiming at modelling
time flows of concurrent processes and their properties by various
algebraic topological models \cite{MR1683333} \cite{model3}
\cite{diCW} \cite{mg} \cite{SK} \cite{FR} \cite{survol}
(cf. \cite{rvg} for other references). In this work, we are interested
in concurrent processes arising from process algebras \cite{MR1365754}
\cite{0683.68008} \cite{0628.68025}, and more precisely in the
labelling process of these objects, which is related to combinatorics
in a non-trivial way. By borrowing several ideas coming from
\cite{exHDA} and \cite{labelled} (see also \cite{Pratt}
\cite{MR1461821} \cite{rvg} \cite{fahrenberg05-hda-long}
\cite{fahrenberg05-hda}), with several slight modifications, the paper
\cite{ccsprecub} presented a semantics of process algebras in terms of
\textit{labelled precubical sets}. We consider in this paper only the
case of Milner's calculus of communicating systems (CCS). The
adaptation to other synchronization algebras and therefore to other
process algebras is straightforward and is left to the reader.

The principle of this semantics is that the concurrent execution of
$n$ actions is abstracted by a full labelled $n$-cube. Each coordinate
corresponds to one of the $n$ actions, and therefore two opposite
faces are labelled by the same action (e.g., Figure~\ref{concab}
represents the concurrent execution of two actions $a$ and $b$). The
core of the construction of \cite{ccsprecub} is the non-functorial
notion of \textit{labelled directed coskeleton}. It is applied to the
fibered product of the $1$-dimensional parts of two full labelled
cubes representing two higher dimensional transitions. This
construction is the key ingredient to defining the parallel
composition with synchronization of CCS in \cite{ccsprecub}. It is
defined as a subobject of the labelled ($1$-dimensional) coskeleton.
The latter coincides with the usual coskeleton, i.e. the right adjoint
to the truncation functor, when the set of labels is a singleton.  The
labelled directed coskeleton construction $\COSK^\Sigma$ takes a
particular kind of $1$-dimensional labelled precubical set $K$ (the
set of vertices $K_0$ must be a cube) to a higher dimensional labelled
precubical set $\COSK^\Sigma(K) $ such that each set of $n$ actions
running concurrently is assembled to an $n$-cube, in exactly one
way. This role cannot be played by the full labelled $1$-dimensional
coskeleton functor $\cosk_1^{\square,\Sigma}$ (see
Proposition~\ref{i1_relatif}) of the category of labelled precubical
sets since the latter may add several different $n$-cubes for the same
set of $n$ actions running concurrently.

The purpose of this paper is to introduce the notion of
\textit{labelled transverse symmetric precubical set}. A transverse
symmetric precubical set is a precubical set equipped with symmetry
maps as in \cite{MR1988396} and with a new kind of degeneracy map
called transverse degeneracy. To the best of our knowledge, the latter
maps seem to be new. In this new category of precubical sets, the
labelled coskeleton functor is well behaved, as explained in
Theorem~\ref{but_math} and Theorem~\ref{cavavraiment}. Indeed, the
labelled transverse symmetric precubical set
$\mathcal{L}(\COSK^\Sigma(K))$ freely generated by the labelled
directed coskeleton $\COSK^\Sigma(K)$ of $K$ is isomorphic to the
labelled coskeleton functor $\cosk_1^{\widehat{\square},\Sigma}(K)$ of
the category of labelled transverse symmetric precubical sets applied
to $K$ if $K$ is the $1$-dimensional part of an $n$-cube or a fibered
product over a synchronization algebra. Since the labelled transverse
symmetric precubical set $\mathcal{L}(\COSK(K))$ and the labelled
precubical set $\COSK(K)$ generate the same topological space of
execution paths by Proposition~\ref{factor_rea} and
Figure~\ref{conclusion}, this result gives a functorial interpretation
of the labelled directed coskeleton construction which is equivalent
to the non-functorial construction from a directed algebraic
topological point of view.

The labelled coskeleton functor in the category of labelled transverse
symmetric precubical sets is therefore a categorical machinery
allowing the understanding of the combinatorics of the labelling
process in the parallel composition with synchronization of CCS. The
advantage of this labelled coskeleton functor is twofold: 1) it is a
functorial construction; 2) it is defined for \textit{any} labelled
$1$-dimensional [transverse symmetric]~\footnote{The words ``transverse
  symmetric'' can be omitted here by Proposition~\ref{1} and
  Proposition~\ref{2}.} precubical set, allowing future
generalizations.

This enables us to give a semantics of CCS in terms of labelled
transverse symmetric precubical sets which is equivalent to the one of
\cite{ccsprecub} in terms of labelled precubical sets from a directed
algebraic topological point of view: see Theorem~\ref{LLL} and
Figure~\ref{conclusion}.

\section{Outline of the paper and reading guide}

The paper is divided in three parts:
\begin{enumerate}
\item Section~\ref{generalisation}, Section~\ref{geo} and
  Section~\ref{cosk} generalize notions previously introduced in
  \cite{ccsprecub} to any category of cubes.  
\item Section~\ref{forma} contains the mathematical treatment.  The
  reader only interested in the computer-scientific applications will
  only have to read the statement of Theorem~\ref{but_math}.
\item Section~\ref{OK}, Section~\ref{compare} and
  Section~\ref{compare2} are the computer-scientific part of the
  paper.
\end{enumerate}


The core of the paper is the categorical interpretation of the
non-functorial labelled directed coskeleton construction using a
generalization of the notion of labelled precubical set.  The notion
of category of cubes, and the generalized notion of labelled
precubical sets are presented in Section~\ref{generalisation}. The
main difficulty is the definition of the generalized precubical set of
labels.  Section~\ref{geo} proves that all the notions of labelled
precubical sets are equivalent from a directed algebraic topological
point of view, in particular that they generate the same path space of
execution paths. This section is the only topological one of the
paper.  Proposition~\ref{factor_rea} is only used in Theorem~\ref{LLL}
to conclude that the two semantics of CCS generate the same spaces of
execution paths. There is also a small application
(Proposition~\ref{maximal}) which is used inside the proof of
Theorem~\ref{cavavraiment}. The topological material of
Section~\ref{geo} is not necessary for the proof of
Proposition~\ref{maximal} but a pure combinatorial proof would be far
more complicated.  Section~\ref{cosk} generalizes to all categories of
precubical sets the labelled coskeleton functor. It is defined as a
right adjoint of a truncation functor, as in the setting of labelled
precubical sets.

Section~\ref{forma} is the mathematical core of the paper. It proves
that all labelled coskeleton functors but one are defective. Indeed,
the labelled coskeleton of the $1$-dimensional part of the $n$-cube is
never contractible in a directed algebraic topological sense, except
for the unique \textit{shell-complete} category of cubes, the maximal
one containing all adjacency-preserving maps. This is the key property
to obtaining a well-behaved labelled coskeleton functor (see
Theorem~\ref{but_math}).  A presheaf over the unique shell-complete
category of cubes is called a \textit{transverse symmetric precubical
  set}.

Section~\ref{OK} is the first section of the computer-scientific part
of the paper. It explains how one can use the preceding constructions
to represent the parallel composition in CCS of a $m$-transition with
a $n$-transition, modelled by a full labelled $m$-cube and a full
labelled $n$-cube respectively. In other terms, it studies parallel
composition in the local case. It is shown that the definition of the
fibered product in CCS must be slightly modified to allow the use of
the labelled coskeleton functor of the category of labelled transverse
symmetric precubical sets.  Section~\ref{compare} then studies
parallel composition in CCS in the global case. It compares the two
notions of synchronized tensor products in the category of labelled
precubical sets and in that of labelled transverse symmetric ones. It
is then proved in Section~\ref{compare2} that the two semantics of CCS
in terms of labelled precubical sets and labelled transverse symmetric
ones are equivalent from a directed algebraic topological point of
view. 

Finally, Section~\ref{labsym} is an additional section treating the
particular case of labelled symmetric precubical sets. This formalism
will enable us to establish a link between concurrent processes viewed
as precubical sets and Cattani-Sassone higher dimensional transition
systems in \cite{hdts}.

\section{Prerequisites}

The paper \cite{ccsprecub} contains an introduction to CCS for
mathematicians which is enough to understand Section~\ref{compare} and
\ref{compare2} of this paper. Computer scientists might prefer
\cite{0683.68008} and \cite{MR1365754}. For the rest of the paper,
only general knowledge in category theory \cite{MR1712872}
\cite{MR1300636} is required, in particular in presheaf theory and in
the theory of locally presentable categories \cite{MR95j:18001}. A few
model category techniques are used in Section~\ref{geo}. In fact,
except for Section~\ref{geo}, the rest of the paper is purely
combinatorial.  Possible references for model categories are
\cite{MR1361887} \cite{MR99h:55031} and \cite{ref_model2}.

\part{About labelled precubical sets over categories of cubes}
\label{tronc_commun}

\section{Labelled precubical set over a category of cubes}
\label{generalisation}

We want to generalize the notion of labelled precubical set introduced
in \cite{ccsprecub} by working on a category of cubes $\mathcal{A}$
(see Definition~\ref{def_cube}) instead on the reduced box category
$\square$ (see Definition~\ref{reduced}) as in \cite{ccsprecub}. The
particular case $\mathcal{A} = \square$ will give back the notion of
labelled precubical set.

\subsection*{Category of cubes (definition and examples)} The category
of partially ordered sets or posets together with the strictly
increasing maps ($x<y$ implies $f(x)<f(y)$) is denoted by $\poset$. It
is worth noting that it is not the usual category of partially ordered
sets since we restrict to strictly increasing maps.  Let $[0] =
\{()\}$ and $[n] = \{0,1\}^n$ for $n \geq 1$.  By convention, one has
$\{0,1\}^0 = [0] = \{()\}$. The set $[n]$ is equipped with the product
ordering $\{0<1\}^n$: $(\epsilon_1, \dots, \epsilon_n) \leq
(\epsilon'_1, \dots, \epsilon'_n)$ if and only if for every $1 \leq i
\leq n$ , one has $\epsilon_i \leq \epsilon'_i$.  The poset $[n]$ is
also called the \textit{$n$-cube}.

\bd Let $\delta_i^\alpha : [n-1] \rightarrow [n]$ be the set map
defined for $1\leq i\leq n$ and $\alpha \in \{0,1\}$ by
$\delta_i^\alpha(\epsilon_1, \dots, \epsilon_{n-1}) = (\epsilon_1,
\dots, \epsilon_{i-1}, \alpha, \epsilon_i, \dots, \epsilon_{n-1})$.
These maps are called the {\rm face maps}. \ed

They satisfy the cocubical relations $\delta_j^\beta \delta_i^\alpha =
\delta_i^\alpha \delta_{j-1}^\beta $ for $i<j$ and for all
$(\alpha,\beta)\in \{0,1\}^2$.

\bd \label{reduced} The {\rm reduced box category}, denoted by
$\square$, is the subcategory of $\poset$ with the set of objects
$\{[n],n\geq 0\}$ and generated by the morphisms
$\delta_i^\alpha$. \ed

It is well-known that the face maps together with the cocubical
relations give a presentation by generators and relations of the small
category $\square$ \cite{MR1988396}.

\bp \label{distance} Let $n\geq 1$. Let
$(\epsilon_1,\dots,\epsilon_n)$ and $(\epsilon'_1,\dots,\epsilon'_n)$
be two elements of the poset $[n]$ with $(\epsilon_1,\dots,\epsilon_n)
\leq (\epsilon'_1,\dots,\epsilon'_n)$. Then there exist
$i_1>\dots>i_{n-r}$ and $\alpha_1,\dots,\alpha_{n-r} \in \{0,1\}$ such
that $(\epsilon_1,\dots,\epsilon_n) = \delta_{i_1}^{\alpha_1} \dots
\delta_{i_{n-r}}^{\alpha_{n-r}}(0\dots 0)$ and
$(\epsilon'_1,\dots,\epsilon'_n) = \delta_{i_1}^{\alpha_1} \dots
\delta_{i_{n-r}}^{\alpha_{n-r}}(1\dots 1)$ where $r\geq 0$ is the
number of $0$ (resp. $1$) in the arguments $0\dots 0$ (resp. $1\dots
1$).  In other terms, $(\epsilon_1,\dots,\epsilon_n)$ is the bottom
element and $(\epsilon'_1,\dots,\epsilon'_n)$ the top element of a
$r$-dimensional subcube of $[n]$.  \ep

\bpf The set $\{1,\dots,n\}$ is equal to the disjoint union
\[\{i\in \{1,\dots,n\}, \epsilon_i = \epsilon'_i\} \sqcup \{i\in
\{1,\dots,n\}, \epsilon_i < \epsilon'_i\}.\] In the latter case, one
necessarily has $\epsilon_i = 0$ and $\epsilon'_i = 1$.  \epf

\bd Let $n\geq 1$. Let $(\epsilon_1,\dots,\epsilon_n)$ and
$(\epsilon'_1,\dots,\epsilon'_n)$ be two elements of the poset
$[n]$. The integer $r$ of Proposition~\ref{distance} is called the
{\rm distance} between $(\epsilon_1,\dots,\epsilon_n)$ and
$(\epsilon'_1,\dots,\epsilon'_n)$.  Let us denote this situation by $r
= d((\epsilon_1,\dots,\epsilon_n),(\epsilon'_1,\dots,\epsilon'_n))$.
By definition, one has \[r = \sum_{i=1}^{i=n}
|\epsilon_i-\epsilon'_i|.\] \ed

\bd A set map $f:[m] \rightarrow [n]$ is {\rm adjacency-preserving} if
it is strictly increasing and if
$d((\epsilon_1,\dots,\epsilon_m),(\epsilon'_1,\dots,\epsilon'_m)) = 1$
implies
$d(f(\epsilon_1,\dots,\epsilon_m),f(\epsilon'_1,\dots,\epsilon'_m)) =
1$.  \ed

An adjacency-preserving map does not necessarily preserve
distance. For example, the map $\gamma_1:[2] \rightarrow [2]$ defined
by $\gamma_1(\epsilon_1,\epsilon_2) = (\max(\epsilon_1,\epsilon_2),
\min(\epsilon_1,\epsilon_2))$ is adjacency-preserving and not
distance-preserving because $\gamma_1(1,0) = \gamma_1(0,1) =
(1,0)$. We shall later see that $\gamma_1$ is an example of transverse
degeneracy map (cf. Definition~\ref{def_trans}).

\bp \label{example_preservation_distance} For any $n\geq 1$, the set
map $\delta_i^\alpha : [n-1] \rightarrow [n]$ is adjacency-preserving. Any
strictly increasing map from $\{0 < 1\}^n$ to itself is
adjacency-preserving as well.  \ep

\bpf That the set map $\delta_i^\alpha : [n-1] \rightarrow [n]$ is
adjacency-preserving is clear. Let $f$ be a strictly increasing map
from $\{0 < 1\}^n$ to itself. Let $x$ and $y$ be two elements of $\{0
< 1\}^n$ with $d(x,y) = 1$ and, for example, $x < y$. Then there
exists a strictly increasing chain $(0,\dots,0) = x_0 < x_1 < \dots <
x_n = (1,\dots,1)$ of $\{0<1\}^n$ with $x = x_{i-1}$ and $y = x_{i}$
for some $i \geq 1$. Then $f(x_0) < f(x_1) < \dots < f(x_n)$ is a
strictly increasing chain of $\{0<1\}^n$. Therefore one has $f(x_0) =
x_0$ and $f(x_n) = x_n$. It is easy to see that $n = d(f(x_0),f(x_n))
= \sum_{i=1}^{n} d(f(x_{i-1}),f(x_{i}))$. So for all $i\geq 1$, one
has $d(f(x_{i-1}),f(x_{i})) = 1$.  Thus, $f$ is
adjacency-preserving. \epf

\bd \label{def_cube} A {\rm category of cubes} $\mathcal{A}$ is a
subcategory of $\poset$ such that:
\begin{itemize}
\item the set of objects is $\{[n],n\geq 0\}$
\item there is the inclusion $\square \subset \mathcal{A}$
\item every morphism of $\mathcal{A}$ is adjacency-preserving.
\end{itemize} \ed

The minimal category of cubes for inclusion is the reduced box
category $\square$. 

\begin{nota} Let us denote by $\widehat{\square}$ the subcategory of
  $\poset$ containing all adjacency-preserving maps. \end{nota}

The category $\widehat{\square}$ is the maximal category of cubes for
inclusion. In other terms, a small category $\C$ is a category of
cubes if and only there are the inclusions $\square \subset \C \subset
\widehat{\square}$.

\begin{nota} For the sequel, $\mathcal{A}$ always denotes a category
  of cubes. \end{nota}

\bd \cite{MR1988396} \label{def_sym} Let $\sigma_i:[n] \rightarrow
[n]$ be the set map defined for $1\leq i\leq n-1$ and $n\geq 2$ by
$\sigma_i(\epsilon_1, \dots, \epsilon_{n}) = (\epsilon_1, \dots,
\epsilon_{i-1},\epsilon_{i+1},\epsilon_{i},
\epsilon_{i+2},\dots,\epsilon_{n})$. These maps are called the {\rm
  symmetry maps}. \ed

The symmetry maps are clearly adjacency-preserving. 

\begin{nota} Let us denote by $\square_S$ the smallest category of
  cubes containing the symmetry maps. \end{nota}

We have the inclusions of categories of cubes $\square \subset
\square_S \subset \widehat{\square}$.

\subsection*{Unlabelled $\mathcal{A}$-set}

\bd An {\rm (unlabelled) $\mathcal{A}$-set} is a presheaf over
$\mathcal{A}$. The corresponding category is denoted by
$\mathcal{A}^{op}\set$.\ed

Let $K$ be an object of $\mathcal{A}^{op}\set$. The set $K([n])$ will
be also denoted by $K_n$. A map $f:K \rightarrow L$ of
$\mathcal{A}^{op}\set$ will be also denoted by $(f_n)_{n\geq 0}$ where
$f_n:K_n \rightarrow L_n$ is the corresponding set map. For any map
$k:[m] \rightarrow [n]$ of $\mathcal{A}$ and any $\mathcal{A}$-set
$K$, denote by $k^*:K_n \rightarrow K_m$ the set map induced by $k$.

Let $p\geq 0$. The \textit{$p$-dimensional $\mathcal{A}$-cube or
  $p$-cube} $\mathcal{A}[p]$ is by definition the presheaf
$\mathcal{A}(-,[p])$. In other terms, $\mathcal{A}[p]_k$ is the set of
maps from $[k]$ to $[p]$ in the category of cubes $\mathcal{A}$.  The
\textit{boundary} $\de\mathcal{A}[p]$ of the $p$-dimensional
$\mathcal{A}$-cube is the presheaf defined by
$\de\mathcal{A}[p]_k=\mathcal{A}[p]_k$ if $k<p$ and
$\de\mathcal{A}[p]_p = \varnothing$ otherwise. In particular, the
boundary of the $0$-dimensional $\mathcal{A}$-cube is the empty
presheaf.

Let $\mathcal{A}_n\subset \mathcal{A}$ be the full subcategory of
$\mathcal{A}$ whose set of objects is $\{[k],k\leq n\}$.  The category
of presheaves over $\mathcal{A}_n$ is denoted by
$\mathcal{A}_n^{op}\set$.  Its objects are called the
\textit{$n$-dimensional $\mathcal{A}$-sets}. The category of
$n$-dimensional $\mathcal{A}$-sets can be identified with the full
subcategory of the category of $\mathcal{A}$-sets $K$ such that
$K_p=\varnothing$ for $p>n$.

Let $K$ be an $\mathcal{A}$-set. Let $K_{\leq n}$ be the
$\mathcal{A}$-set obtained from $K$ by keeping the $p$-dimensional
cubes of $K$ only for $p\leq n$. In particular, $K_{\leq 0}=K_0$. Note
that one has $\de\mathcal{A}[n] = \mathcal{A}[n]_{\leq n-1}$ for every
$n\geq 0$ since our precubical sets contain no degeneracy maps in the
usual sense.

\bd A $\square$-set is called a {\rm precubical set}
\cite{Brown_cube}.  A $\square_S$-set is called a {\rm symmetric
  precubical set} \cite{MR1988396}.  A $\widehat{\square}$-set is
called a {\rm transverse symmetric precubical set}.~\footnote{Note that
  the last notion is new.} \ed

The inclusion functor $\square \subset \mathcal{A}$ induces a
forgetful functor $\omega_\mathcal{A} : \mathcal{A}^{op}\set
\rightarrow \square^{op}\set$ which has both a left and a right
adjoint obtained respectively as a left and a right Kan extension
along the inclusion $\square^{op} \subset \mathcal{A}^{op}$.  The
right adjoint is denoted by $\mathcal{R}_\mathcal{A}:\square^{op}\set
\rightarrow \mathcal{A}^{op}\set$.  The left adjoint
$\mathcal{L}_\mathcal{A}:\square^{op}\set \rightarrow
\mathcal{A}^{op}\set$ is of special interest since it formally adds
all additional operators defining an $\mathcal{A}$-set. The two
following propositions state some elementary remarks about
$\mathcal{L}_\mathcal{A}$ which will be reused later.

\bp \label{carre} Let $K$ be a precubical set. Then one has the
isomorphism
\[\mathcal{L}_\mathcal{A}(K) \iso \liminj_{\square[n] \rightarrow K}
\mathcal{A}[n].\] In particular, there is the isomorphism of
$\mathcal{A}$-sets $\mathcal{L}_\mathcal{A}(\square[n])\iso
\mathcal{A}[n]$.  \ep

\bpf For every $\mathcal{A}$-set $K$, one has $K_n=(\omega_\mathcal{A}
K)_n$ for all $n\geq 0$ since the inclusion functor $\square \subset
\mathcal{A}$ is the identity on objects. So one has the bijections of
sets
\[\mathcal{A}^{op}\set(\mathcal{L}_\mathcal{A}(\square[n]),K) \iso
\square^{op}\set(\square[n],\omega_\mathcal{A} K) =
(\omega_\mathcal{A} K)_n = K_n =
\mathcal{A}^{op}\set(\mathcal{A}[n],K).\] By the Yoneda lemma, one
obtains the isomorphism $\mathcal{L}_\mathcal{A}(\square[n]) \iso
\mathcal{A}[n]$ for all $n\geq 0$.  Since $\mathcal{L}_\mathcal{A}$ is
a left adjoint, it preserves colimits. So one obtains for every
precubical set $K$
\[\mathcal{L}_\mathcal{A}(K) = \mathcal{L}_\mathcal{A}\lp \liminj_{\square[n] \rightarrow K}
\square[n] \rp \iso \liminj_{\square[n] \rightarrow K} \mathcal{A}[n].\] 
\epf 

\bp \label{pintuitive0} Let $K$ be a precubical set. The identity map
$\id_{\mathcal{L}_{\mathcal{A}}(K)}$ induces by adjunction an
inclusion of presheaves $i_K : K \subset
\omega_\mathcal{A}\mathcal{L}_\mathcal{A}(K)$.  \ep

\bpf Since the functor $K\mapsto K_p$ from precubical sets to sets is
colimit-preserving for every $p\geq 0$, one has the bijections
\[K_p\iso \liminj_{\square[n]\rightarrow K} \square[n]_p\] 
and 
\[\omega_\mathcal{A}\mathcal{L}_\mathcal{A}(K)_p\iso
\liminj_{\square[n]\rightarrow K} \omega_\mathcal{A}\mathcal{A}[n]_p =
\liminj_{\square[n]\rightarrow K}\mathcal{A}([p],[n]).\] Each set map
$\square[n]_p \rightarrow \omega_\mathcal{A}\mathcal{A}[n]_p$ is
one-to-one because of the inclusions of sets $\square([p],[n]) \subset
\mathcal{A}([p],[n])$ for every $p\geq 0$. For any map $g: [n]
\rightarrow [n']$ of $\square$, one has the commutative diagram of
sets
\[
\xymatrix{
\square[n]_p \fr{\subset} \fd{\subset} && \square[n']_p \fd{\subset}\\
&&\\
\omega_\mathcal{A}\mathcal{A}[n]_p = \mathcal{A}([p],[n])\fr{}  &&
\omega_\mathcal{A}\mathcal{A}[n']_p= \mathcal{A}([p],[n']).}
\] 
The bottom map is one-to-one since it consists of composing by $g$
which is one-to-one as any map of $\square$. So each set map
$\omega_\mathcal{A}\mathcal{A}[n]_p \rightarrow
\omega_\mathcal{A}\mathcal{A}[n']_p$ of the diagram calculating
$\omega_\mathcal{A}\mathcal{L}_\mathcal{A}(K)_p$ is one-to-one as
well. One deduces that the map $K_p \rightarrow
\omega_\mathcal{A}\mathcal{L}_\mathcal{A}(K)_p$ is one-to-one. \epf

\begin{nota} Since $K \subset \omega_\mathcal{A
  }\mathcal{L}_\mathcal{A}(K)$ is an inclusion, $i_K(y)$ will be
  simply denoted by $y$ for any $y\in K$.
\end{nota}

\subsection*{The $1$-dimensional case}

This paragraph proves that the $1$-dimensional case does not depend on
the choice of the category of cubes. The crucial facts are that a
category of cubes contains all face maps and that all morphisms are
adjacency-preserving.

\bp \label{ex} For every $m>n$, one has $\mathcal{A}([m],[n]) =
\varnothing$. For every $n\geq 0$, the inclusion $\square \subset
\mathcal{A}$ implies the bijections $[n] \iso \square([0],[n]) \iso
\mathcal{A}([0],[n])$ and $\square([1],[n]) \iso
\mathcal{A}([1],[n])$.  \ep

Note that this implies that $\mathcal{A}$ cannot have any
degeneracies.

\bpf It is clear that $\mathcal{A}([m],[n]) \neq \varnothing$ implies
$m\leq n$. One has the inclusions \[\square([0],[n]) \subset
\mathcal{A}([0],[n]) \subset \widehat{\square}([0],[n]) =
\square([0],[n]) = \{\delta^{\epsilon_n}_n\ldots \delta^{\epsilon_1}_1,
(\epsilon_1,\dots,\epsilon_n)\in[n]\}\] hence the second
assertion. For every $n\geq 0$, the inclusion $\square([1],[n])\subset
\mathcal{A}([1],[n])$ is a bijection since every map of $\mathcal{A}$
is adjacency-preserving by definition of a category of cubes, hence
the third assertion.  \epf

\bp \label{restrict1_1} Let $K$ be a precubical set. Then the
inclusion of precubical sets $K \subset \omega_\mathcal{A}
\mathcal{L}_\mathcal{A}(K)$ induces the isomorphism of $1$-dimensional
precubical sets $K_{\leq 1} \iso \omega_\mathcal{A}
\mathcal{L}_\mathcal{A}(K)_{\leq 1}$. \ep

\bpf We already know by Proposition~\ref{ex} that for every $n\geq 0$,
the inclusions $\square([0],[n])\linebreak[4]\subset
\mathcal{A}([0],[n])$ and $\square([1],[n])\subset
\mathcal{A}([1],[n])$ are bijective. So the inclusion of presheaves
$\square[n]_{\leq 1} \subset \omega_\mathcal{A} \mathcal{A}[n]_{\leq
  1}$ is an isomorphism for every $n\geq 0$. Since the forgetful
functor $\omega_\mathcal{A} : \mathcal{A}^{op}\set \rightarrow
\square^{op}\set$ is a left adjoint, it is colimit-preserving. Hence
the proof is complete.  \epf

\bp \label{1} The category of $1$-dimensional precubical sets is
equivalent to the category of $1$-dimensional $\mathcal{A}$-sets. \ep

\bpf The adjunction $\mathcal{L}_\mathcal{A} : \square^{op}\set
\leftrightarrows \mathcal{A}^{op}\set:\omega_\mathcal{A}$ induces an
adjunction $(\mathcal{L}_\mathcal{A})_{\leq 1} : \square_1^{op}\set
\leftrightarrows \mathcal{A}_1^{op}\set:(\omega_\mathcal{A})_{\leq 1}$
by Proposition~\ref{carre}. We already know by
Proposition~\ref{restrict1_1} that there is the isomorphism
$(\omega_\mathcal{A})_{\leq 1} (\mathcal{L}_\mathcal{A})_{\leq 1} \iso
\id_{\square_1^{op}\set}$.  One has $(\mathcal{L}_\mathcal{A})_{\leq
  1}(\omega_\mathcal{A})_{\leq 1}(\mathcal{A}[0]) \iso \mathcal{A}[0]$
by Proposition~\ref{carre} and $(\mathcal{L}_\mathcal{A})_{\leq
  1}(\omega_\mathcal{A})_{\leq 1}(\mathcal{A}[1]) \iso \mathcal{A}[1]$
by Proposition~\ref{ex} and Proposition~\ref{carre}. Hence the
isomorphism of functors $(\mathcal{L}_\mathcal{A})_{\leq
  1}(\omega_\mathcal{A})_{\leq 1} \iso \id_{\mathcal{A}_1^{op}\set}$.
\epf

\subsection*{Labelled $\mathcal{A}$-set}

We fix a non-empty set $\Sigma$ of \textit{labels} or of
\textit{actions}. It always contains a distinguished label denoted by
$\tau$. We want to label the cubes of a $\mathcal{A}$-set with the
elements of $\Sigma$. A labelled $\mathcal{A}$-set will be a map of
$\mathcal{A}$-sets $K \rightarrow L$ where $L$ is the
$\mathcal{A}$-set of labels. Let us start by recalling the
construction of the precubical set of labels.

\bp (Variant of Goubault's construction
\cite{labelled}) \label{cubeetiquette} Let
\begin{itemize}
\item $(!\Sigma)_0=\{()\}$ (the empty word)
\item for $n\geq 1$, $(!\Sigma)_n=\Sigma^n$
\item $\de_i^0(a_1,\dots,a_n) = \de_i^1(a_1,\dots,a_n) =
  (a_1,\dots,\widehat{a_i},\dots,a_n)$ where the notation
  $\widehat{a_i}$ means that $a_i$ is removed.
\end{itemize}
Then these data generate a precubical set denoted by $!\Sigma$. \ep

\begin{figure}
\[
\xymatrix{
& () \ar@{->}[rd]^{(b)}&\\
()\ar@{->}[ru]^{(a)}\ar@{->}[rd]_{(b)} & (a,b) & ()\\
&()\ar@{->}[ru]_{(a)}&}
\]
\caption{Concurrent execution of $a$ and $b$}
\label{concab}
\end{figure}
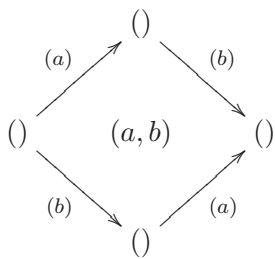

\bd Let $K$ be an $\mathcal{A}$-set. Let $x\in K_p$ with $p\geq
1$. The {\rm boundary of} $x$ is the composite map $\de x:
\de\mathcal{A}[p] \subset \mathcal{A}[p] \stackrel{x} \longrightarrow
K$.  \ed

The main feature of the precubical set $!\Sigma$ is that for every
$p\geq 2$, a $p$-cube of $!\Sigma$, which labels the concurrent
execution of $p$ actions like in Figure~\ref{concab}, is determined by
its boundary. In other terms, a commutative square of precubical sets
of the form
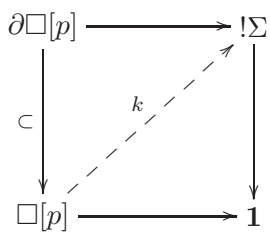
\[
\xymatrix{
\de\square[p]\fr{} \fd{\subset}&& !\Sigma \fd{} \\
&&\\
\square[p] \fr{} \ar@{-->}[rruu]^-{k} && \mathbf{1}}
\] 
with $p\geq 2$ where $\mathbf{1}$ is the terminal precubical set
admits at most one lift $k$. An equivalent mathematical formulation of
the preceding condition is that for every commutative square of
precubical sets of the form
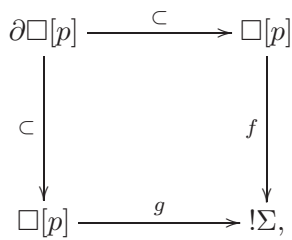
\[
\xymatrix{
\de\square[p] \fr{\subset} \fd{\subset}&& \square[p] \fd{f} \\
&&\\
\square[p] \fr{g}  && !\Sigma,}
\] 
one has $f = g$. So every commutative square of precubical sets of the
form
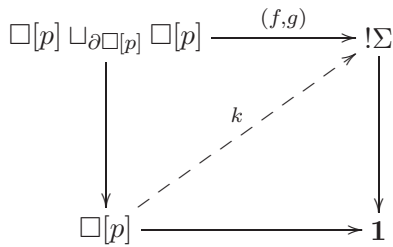
\[
\xymatrix{
\square[p]\sqcup_{\de\square[p]}\square[p]\fr{(f,g)} \fd{}&& !\Sigma \fd{} \\
&&\\
\square[p] \fr{} \ar@{-->}[rruu]^-{k} && \mathbf{1}}
\] 
with $p\geq 2$ admits exactly one lift $k = f =g$. In other terms, the
precubical set $!\Sigma$ turns out to be orthogonal to the set of maps
$\{\square[p]\sqcup_{\de\square[p]}\square[p] \rightarrow
\square[p],p\geq 2\}$ in the sense of
\cite[Definition~1.32]{MR95j:18001}.

Because of the inclusion $K \subset \omega_\mathcal{A}
\mathcal{L}_\mathcal{A}(K)$ for every precubical set $K$, we need more
cubes for the $\mathcal{A}$-set of labels as soon as the inclusion
$\square \subset \mathcal{A}$ is strict. Indeed, we must be able to
label all cubes of $\mathcal{L}_\mathcal{A}(K)$ for every labelled
precubical set $K \rightarrow !\Sigma$.  The first candidate for the
$\mathcal{A}$-set of labels is then the $\mathcal{A}$-set
$\mathcal{L}_\mathcal{A}(!\Sigma)$ freely generated by
$!\Sigma$. However, it is not well-behaved. Consider the two set
involutions $\sigma_1 : [2] \rightarrow [2]$ and $\id_{[2]} : [2]
\rightarrow [2]$. Let us suppose that $\sigma_1 \in \mathcal{A}$. Then
the two $2$-cubes $(\sigma_1)^{*}(\tau,\tau)$ and $(\tau,\tau)$ of
$\mathcal{L}_\mathcal{A}(!\Sigma)$ have the same boundary. This means
that the commutative square of $\mathcal{A}$-sets
\[
\xymatrix{
  \de\mathcal{A}[p]\ar@{->}[rrrr]^-{\de(\sigma_1)^{*}(\tau,\tau) = \de(\tau,\tau)} \fd{\subset}&&&& \mathcal{L}_\mathcal{A}(!\Sigma) \fd{} \\
  &&\\
  \mathcal{A}[p] \ar@{->}[rrrr] \ar@{-->}[rrrruu]^-{k} &&&& \mathbf{1}}
\] 
has two distinct lifts $k = (\sigma_1)^{*}(\tau,\tau)$ and $k =
(\tau,\tau)$. In other terms, the $\mathcal{A}$-set
$\mathcal{L}_\mathcal{A}(!\Sigma)$ is never orthogonal to the set of
morphisms $\{\mathcal{A}[p]\sqcup_{\de\mathcal{A}[p]}\mathcal{A}[p]
\rightarrow \mathcal{A}[p], p\geq 2\}$ as soon as $\sigma_1:[2]
\rightarrow [2]$ belongs to $\mathcal{A}$. In fact, the
$\mathcal{A}$-set $\mathcal{L}_\mathcal{A}(!\{\tau\})$ is even not the
terminal $\mathcal{A}$-set in this case. Yet, the notion of
$\mathcal{A}$-set must coincide with the unlabelled notion if the set
of labels is equal to $\{\tau\}$. The full subcategory
$\{\mathcal{A}[p]\sqcup_{\de\mathcal{A}[p]}\mathcal{A}[p]\rightarrow
\mathcal{A}[p], p\geq 2\}^\bot$ of $\mathcal{A}$-sets orthogonal to
the set of maps
$\{\mathcal{A}[p]\sqcup_{\de\mathcal{A}[p]}\mathcal{A}[p]\rightarrow
\mathcal{A}[p], p\geq 2\}$ is a full reflective subcategory of the
locally presentable category of $\mathcal{A}$-sets by
\cite[Theorem~1.39]{MR95j:18001}. Let
\[\boxed{\sh_\mathcal{A}:\mathcal{A}^{op}\set \rightarrow
\{\mathcal{A}[p]\sqcup_{\de\mathcal{A}[p]}\mathcal{A}[p]\rightarrow
\mathcal{A}[p], p\geq 2\}^\bot}\] be the left adjoint to the inclusion
functor
$\{\mathcal{A}[p]\sqcup_{\de\mathcal{A}[p]}\mathcal{A}[p]\rightarrow
\mathcal{A}[p], p\geq 2\}^\bot \subset \mathcal{A}^{op}\set$.

\bd The {\rm $\mathcal{A}$-set of labels} is the $\mathcal{A}$-set
$\sh_\mathcal{A}\mathcal{L}_\mathcal{A}(!\Sigma)$. \ed

In $\sh_\mathcal{A}\mathcal{L}_\mathcal{A}(!\Sigma)$, the two
$2$-cubes $(\sigma_1)^{*}(\tau,\tau)$ and $(\tau,\tau)$ are forced to
be equal. Note that there is the isomorphism of precubical sets
$\sh_\square\mathcal{L}_\square(!\Sigma) \iso !\Sigma$.

\bd A {\rm labelled $\mathcal{A}$-set (over $\Sigma$)} is an object of
the comma category \[\mathcal{A}^{op}\set \ddownarrow
\sh_\mathcal{A}\mathcal{L}_\mathcal{A}(!\Sigma).\] That is, an object
is a map of $\mathcal{A}$-sets $\ell:K \rightarrow
\sh_\mathcal{A}\mathcal{L}_\mathcal{A}(!\Sigma)$ and a morphism is a
commutative diagram \[ \xymatrix{ K \ar@{->}[rr]\ar@{->}[rd]&& L
  \ar@{->}[ld]\\ & \sh_\mathcal{A}\mathcal{L}_\mathcal{A}(!\Sigma).&}
\]
The $\ell$ map is called the {\rm labelling map}.  The
$\mathcal{A}$-set $K$ is sometimes called the {\rm underlying
  $\mathcal{A}$-set} of the labelled $\mathcal{A}$-set. \ed

The functor $\mathcal{L}_\mathcal{A} : \square^{op}\set \rightarrow
\mathcal{A}^{op}\set$ induces a functor (denoted in the same way)
\[\mathcal{L}_\mathcal{A} : \square^{op}\set\ddownarrow !\Sigma \rightarrow
\mathcal{A}^{op}\set\ddownarrow
\sh_\mathcal{A}\mathcal{L}_\mathcal{A}(!\Sigma)\]
which takes $\ell:K\rightarrow !\Sigma$ to the composite
$\mathcal{L}_\mathcal{A}(K)
\stackrel{\mathcal{L}_\mathcal{A}(\ell)}\longrightarrow
\mathcal{L}_\mathcal{A}(!\Sigma) \rightarrow
\sh_\mathcal{A}\mathcal{L}_\mathcal{A}(!\Sigma)$.

\bp\label{restrict1_2} Let $K$ be an $\mathcal{A}$-set. Then the map
of $\mathcal{A}$-sets $K \rightarrow \sh_\mathcal{A}(K)$ induces the
isomorphism of $1$-dimensional $\mathcal{A}$-sets $K_{\leq 1} \iso
\sh_\mathcal{A}(K)_{\leq 1}$.  \ep

\bpf For every $p\geq 2$ and for every commutative diagram of solid
arrows
\[
\xymatrix{ \mathcal{A}[p]\sqcup_{\de\mathcal{A}[p]}\mathcal{A}[p]
  \fr{}\fd{} && K \fd{} \\
  && \\
  \mathcal{A}[p] \ar@{-->}[rruu]^-{k} \fr{} && \mathbf{1},}
\] 
there exists at most one lift $k$. So an $\mathcal{A}$-set $K$ is
orthogonal to the set of morphisms
$\{\mathcal{A}[p]\sqcup_{\de\mathcal{A}[p]}\mathcal{A}[p] \rightarrow
\mathcal{A}[p], p\geq 2\}$ if and only if the canonical map
$K\rightarrow \mathbf{1}$ satisfies the right lifting property with
respect to the same set of morphisms. So the $\mathcal{A}$-set
$\sh_\mathcal{A}(K)$ can be obtained by a small object argument by
factoring the map $K\rightarrow \mathbf{1}$ as a composite
$K\rightarrow \sh_\mathcal{A}(K) \rightarrow\mathbf{1}$ where $K
\rightarrow \sh_\mathcal{A}(K)$ is a relative
$\{\mathcal{A}[p]\sqcup_{\de\mathcal{A}[p]}\mathcal{A}[p]\rightarrow
\mathcal{A}[p], p\geq 2\}$-cell complex and where the map
$\sh_\mathcal{A}(K) \rightarrow\mathbf{1}$ satisfies the right lifting
property with respect to the same set of morphisms. The small object
argument is possible by \cite[Proposition~1.3]{MR1780498} since the
category of $\mathcal{A}$-sets is locally presentable, as every
presheaf category. Since for every $p\geq 2$, the map of
$\mathcal{A}$-sets
$\mathcal{A}[p]\sqcup_{\de\mathcal{A}[p]}\mathcal{A}[p] \rightarrow
\mathcal{A}[p]$ induces an isomorphism
\[\lp\mathcal{A}[p]\sqcup_{\de\mathcal{A}[p]}\mathcal{A}[p]
\rp_{\leq 1} \iso \mathcal{A}[p]_{\leq 1},\] one deduces that the
canonical map $K\rightarrow \sh_\mathcal{A}(K)$ induces an isomorphism
$K_{\leq 1} \iso \sh_\mathcal{A}(K)_{\leq 1}$.  \epf

\bp \label{1objet} There is the isomorphism
$\sh_\mathcal{A}\mathcal{L}_\mathcal{A}(!\{\tau\}) \iso
\mathbf{1}$. Therefore when $\Sigma=\{\tau\}$, the category of
labelled $\mathcal{A}$-sets is equivalent to the category of
unlabelled $\mathcal{A}$-sets. \ep
 
\bpf Indeed, both the functors $\sh_\mathcal{A}$ and
$\mathcal{L}_\mathcal{A}$ do not modify the set of $0$-cubes and the
set of $1$-cubes by Proposition~\ref{restrict1_1} and
Proposition~\ref{restrict1_2}. Moreover, for any $\mathcal{A}$-set $K$
such that $K_1$ is a singleton, it is clear by induction on $p\geq 1$
that the set $(\sh_\mathcal{A} K)_p$ is a singleton. So
$\sh_\mathcal{A}\mathcal{L}_\mathcal{A}(!\{\tau\}) \iso \mathbf{1}$
(the terminal object of $\mathcal{A}^{op}\set$). \epf

\begin{nota} Let $(a_1,\dots,a_n)\in \Sigma^n$ with $n\geq 1$. The
  labelled precubical set $\square[a_1,\dots,a_n]$ denotes the map
  $\ell:\square[n]\rightarrow !\Sigma$ such that $\ell(\id_{[n]}) =
  (a_1,\dots,a_n)$. \end{nota}

Figure~\ref{concab} gives the example of the labelled $2$-cube
$\square[a,b]$. It represents the concurrent execution of $a$ and
$b$. It is important to notice that two opposite faces of
Figure~\ref{concab} have the same label.

\begin{nota} Let $(a_1,\dots,a_n)\in \Sigma^n$ with $n\geq 1$. The
  labelled $\mathcal{A}$-set $\mathcal{A}[a_1,\dots,a_n]$ denotes the
  labelled $\mathcal{A}$-set
  $\mathcal{L}_\mathcal{A}(\square[a_1,\dots,a_n])$. \end{nota}

\bp \label{2} The category of labelled $1$-dimensional precubical sets
is equivalent to the category of labelled $1$-dimensional
$\mathcal{A}$-sets.  \ep

\bpf This is a consequence of Proposition~\ref{1} and
Proposition~\ref{restrict1_2}.  \epf

\section{Geometric realization of labelled $\mathcal{A}$-set}
\label{geo}

The purpose of this section is to prove that the geometric realization
functor $\square^{op}\set\ddownarrow !\Sigma \rightarrow \dtop
\ddownarrow ?\Sigma$ of \cite{ccsprecub} which takes a labelled
precubical set to the corresponding labelled flow factors as a
composite $\square^{op}\set \ddownarrow !\Sigma \rightarrow
\mathcal{A}^{op}\set \ddownarrow
\sh_\mathcal{A}\mathcal{L}_\mathcal{A}(!\Sigma) \rightarrow \dtop
\ddownarrow ?\Sigma$, where the left-hand functor is induced by
$\mathcal{L}_\mathcal{A}$. This result ensures that all the notions of
labelled $\mathcal{A}$-sets are equivalent from a directed algebraic
topological point of view. The results of this section are used only
in Theorem~\ref{LLL} and in Theorem~\ref{cavavraiment}. 

\subsection*{Unlabelled flow}

The category $\top$ of \textit{compactly generated topological spaces}
(i.e.\ of weak Hausdorff $k$-spaces) is complete, cocomplete and
cartesian closed (more details for these kinds of topological spaces
are in \cite{MR2273730}, \cite{MR2000h:55002}, the appendix of
\cite{Ref_wH} and also in the preliminaries of \cite{model3}). For the
sequel, all topological spaces will be supposed to be compactly
generated. A \textit{compact space} is always Hausdorff.  

\bd \cite{model3} A {\rm (time) flow} $X$ is a small topological
category without identity maps. The set of objects is denoted by
$X^0$.  The topological space of morphisms from $\alpha$ to $\beta$ is
denoted by $\P_{\alpha,\beta}X$. The elements of $X^0$ are also called
the {\rm states} of $X$. The elements of $\P_{\alpha,\beta}X$ are
called the {\rm (non-constant) execution paths from $\alpha$ to
  $\beta$}. A flow $X$ is {\rm loopless} if for every $\alpha\in X^0$,
the space $\P_{\alpha,\alpha}X$ is empty. \ed

\begin{nota} Let $\P X = \bigsqcup_{(\alpha,\beta)\in X^0\p X^0}
  \P_{\alpha,\beta}X$.  The topological space $\P X$ is called the
  \textit{path space} of $X$.  The source map (resp.\ the target map)
  $\P X\rightarrow X^0$ is denoted by $s$ (resp.\ $t$).
\end{nota}

\bd Let $X$ be a flow, and let $\alpha \in X^0$ be a state of $X$. The
state $\alpha$ is \textit{initial} if $\alpha\notin t(\P X)$, and the
state $\alpha$ is \textit{final} if $\alpha\notin s(\P X)$. \ed

\bd A morphism of flows $f: X \rightarrow Y$ consists in a set map
$f^0: X^0 \rightarrow Y^0$ and a continuous map $\P f: \P X
\rightarrow \P Y$ such that $s(\P f(x)) = f^0(s(x))$, $t(\P f(x)) =
f^0(t(x))$ and $\P f(x * y) = \P f(x) * \P f(y)$ for every $x,y\in \P
X$. The corresponding category is denoted by $\dtop$. \ed

The strictly associative composition law
\[
\left\{ \begin{array}{c} \P_{\alpha,\beta}X \p \P_{\beta,\gamma}X \longrightarrow \P_{\alpha,\gamma}X \\
(x,y) \mapsto x*y \end{array} \right.
\]
models the composition of non-constant execution paths. The
composition law $*$ is extended in a usual way to states, that is to
constant execution paths, by $x*t(x) = x$ and $s(x)*x = x$ for every
non-constant execution path $x$.

Here are two fundamental examples of flows:
\begin{enumerate}
\item Let $S$ be a set. The flow associated with $S$, still denoted by
  $S$, has $S$ as a set of states and the empty space as path space.
  This construction induces a functor $\set \rightarrow \dtop$ from
  the category of sets to that of flows. The flow associated with a
  set is loopless.
\item Let $(P,\leq)$ be a poset. The flow associated with $(P,\leq)$,
  and still denoted by $P$ is defined as follows: the set of states of
  $P$ is the underlying set of $P$; the space of morphisms from
  $\alpha$ to $\beta$ is empty if $\alpha\geq \beta$ and equal to
  $\{(\alpha,\beta)\}$ if $\alpha<\beta$ and the composition law is
  defined by $(\alpha,\beta)*(\beta,\gamma) = (\alpha,\gamma)$. This
  construction induces a functor $\poset \rightarrow \dtop$ from the
  category of posets together with the strictly increasing maps to the
  category of flows. The flow associated with a poset is loopless as
  well.~\footnote{and must be loopless ! This is one of the reasons for
    working with small categories without identity maps.}
\end{enumerate}

There is an important model structure on $\dtop$ which is
characterized as follows \cite{model3}:
\begin{itemize}
\item The weak equivalences are the \textit{weak S-homotopy
    equivalences}, i.e.\ the morphisms of flows $f:
  X\longrightarrow Y$ such that $f^0: X^0\longrightarrow Y^0$ is
  a bijection of sets and such that $\P f: \P X\longrightarrow \P
  Y$ is a weak homotopy equivalence.
\item The fibrations are the morphisms of flows $f: X\longrightarrow
  Y$ such that $\P f: \P X\longrightarrow \P Y$ is a Serre
  fibration\footnote{that is, a continuous map having the right
    lifting property with respect to the inclusion $\mathbf{D}^n\p
    0\subset \mathbf{D}^n\p [0,1]$ for any $n\geq 0$ where
    $\mathbf{D}^n$ is the $n$-dimensional disk.}.
\end{itemize}
This model structure is cofibrantly generated. The cofibrant
replacement functor is denoted by $(-)^{\textit{cof}}$.

\subsection*{Labelled flow}

\bd \cite{ccsprecub} The {\rm flow of labels} $?\Sigma$ is defined as
follows: $(?\Sigma)^0 = \{0\}$ and $\P ?\Sigma$ is the discrete free
commutative semigroup generated by the elements of $\Sigma$.  \ed

\bd \cite{ccsprecub} A {\rm labelled flow} is an object of the comma
category $\dtop \ddownarrow ?\Sigma$. That is an object is a map of
flows $\ell: X \rightarrow ?\Sigma$ and a morphism is a commutative
diagram
\[ \xymatrix{ X \ar@{->}[rr]\ar@{->}[rd]&& Y \ar@{->}[ld]\\ &
  ?\Sigma.&}
\]
The $\ell$ map is called the {\rm labelling map}.  The flow $X$ is
sometimes called the {\rm underlying flow} of the labelled flow.  \ed


\subsection*{Geometric realization of a labelled precubical set}

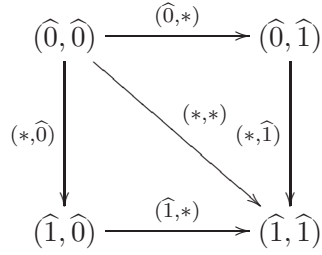
\begin{figure}
\[
\xymatrix{
  (\widehat{0},\widehat{0}) \fr{(\widehat{0},*)}\fd{(*,\widehat{0})}\ar@{->}[ddrr]^-{(*,*)} && (\widehat{0},\widehat{1})\fd{(*,\widehat{1})}\\
  && \\
  (\widehat{1},\widehat{0}) \fr{(\widehat{1},*)}&&
  (\widehat{1},\widehat{1})}
\]
\caption{The flow $\{\widehat{0}<\widehat{1}\}^2$ (Note that
  $(*,*)=(\widehat{0},*)*(*,\widehat{1}) = (*,\widehat{0}) *
  (\widehat{1},*)$)}
\label{cube2}
\end{figure}

A state of the flow associated with the poset
$\{\widehat{0}<\widehat{1}\}^n$ (i.e. the product of $n$ copies of
$\{\widehat{0}<\widehat{1}\}$) is denoted by an $n$-tuple of elements
of $\{\widehat{0},\widehat{1}\}$. By convention,
$\{\widehat{0}<\widehat{1}\}^0 = \{()\}$. The unique
morphism/execution path from $(x_1,\dots,x_n)$ to $(y_1,\dots,y_n)$ is
denoted by an $n$-tuple $(z_1,\dots,z_n)$ of
$\{\widehat{0},\widehat{1},*\}$ with $z_i = x_i$ if $x_i = y_i$ and
$z_i = *$ if $x_i < y_i$. For example in the flow $\{\widehat{0} <
\widehat{1}\}^2$ (cf. Figure~\ref{cube2}), one has the algebraic
relation $(*,*) = (\widehat{0},*)*(*,\widehat{1}) = (*,\widehat{0}) *
(\widehat{1},*)$.

Let $\square \rightarrow \poset \subset \dtop$ be the functor defined
on objects by the mapping $[n]\mapsto \{\widehat{0}<\widehat{1}\}^n$
and on morphisms by the mapping \[\delta_i^\alpha \mapsto \lp
(\epsilon_1, \dots, \epsilon_{n-1}) \mapsto (\epsilon_1, \dots,
\epsilon_{i-1}, \alpha, \epsilon_i, \dots, \epsilon_{n-1})\rp,\] where
the $\epsilon_i$'s are elements of
$\{\widehat{0},\widehat{1},*\}$. The functor $[n] \mapsto
\{\widehat{0}<\widehat{1}\}^n$ from $\square$ to $\dtop$ induces a bad
realization functor from $\square^{op}\set$ to $\dtop$ defined
by \[\boxed{|K|_{bad}:=\liminj_{\square[n]\rightarrow K}
  \{\widehat{0}<\widehat{1}\}^n}.\]

\bth (\cite[Theorem~7.1]{ccsprecub} and
\cite[Proposition~8.1]{ccsprecub}) \label{badsigma} For all $n\geq 3$,
the inclusion $\de\square[n] \subset \square[n]$ induces an
isomorphism of flows $|\de\square[n]|_{bad} \iso
|\square[n]|_{bad}$. One has the isomorphism of flows $|!\Sigma|_{bad}
\iso\ ?\Sigma$.  \eth

\bd \cite{ccsprecub} Let $K$ be a precubical set. By definition, the
{\rm geometric realization} of $K$ is the flow
\[\boxed{|K| := \liminj_{\square[n]\rightarrow
    K} (\{\widehat{0}<\widehat{1}\}^n)^{\textit{cof}}} \]
\ed

The natural trivial fibrations $(\{\widehat{0} <
\widehat{1}\}^n)^{\textit{cof}} \longrightarrow \{\widehat{0} <
\widehat{1}\}^n$ for $n\geq 0$ induce a natural map $|K|
\longrightarrow |K|_{bad}$ for any precubical set $K$. Let
$K\rightarrow !\Sigma$ be a labelled precubical set. Then the
composition $|K| \rightarrow |!\Sigma| \rightarrow |!\Sigma|_{bad}
\iso ?\Sigma$ gives rise to a labelled flow.

\subsection*{Geometric realization of a labelled $\mathcal{A}$-set}

Let $\mathcal{A} \rightarrow \poset \subset \dtop$ be the functor
defined on objects by the mapping $[n]\mapsto
\{\widehat{0}<\widehat{1}\}^n$ and on morphisms as follows. Let $f:[m]
\rightarrow [n]$ be a map of $\mathcal{A}$ with $m,n\geq 0$. Let
$(\epsilon_1,\dots,\epsilon_m)\in \{\widehat{0},\widehat{1},*\}^m$ be
a $r$-cube. Since $f$ is adjacency-preserving, the two elements
$f(s(\epsilon_1,\dots,\epsilon_m))$ and
$f(t(\epsilon_1,\dots,\epsilon_m))$ are respectively the initial and
final states of a unique $r$-dimensional subcube denoted by
$f(\epsilon_1,\dots,\epsilon_m)$ of $[n]$ with
$f(\epsilon_1,\dots,\epsilon_m)\in
\{\widehat{0},\widehat{1},*\}^n$. Note that the composite functor
$\square \subset \mathcal{A} \rightarrow \poset \subset \dtop$ is the
functor defined above. The functor $[n] \mapsto
\{\widehat{0}<\widehat{1}\}^n$ from $\mathcal{A}$ to $\dtop$ induces a
bad realization functor from $\mathcal{A}^{op}\set$ to $\dtop$ defined
by \[\boxed{|K|_{bad}:=\liminj_{\mathcal{A}[n]\rightarrow K}
  \{\widehat{0}<\widehat{1}\}^n}.\]

\bd Let $K$ be an $\mathcal{A}$-set. By definition, the {\rm geometric
  realization} of $K$ is the flow
\[\boxed{|K| := \liminj_{\mathcal{A}[n]\rightarrow
    K} (\{\widehat{0}<\widehat{1}\}^n)^{\textit{cof}}} \]
\ed

Note that the two geometric realizations of $\mathcal{A}$-sets are
colimit-preserving. In fact, it is easy to prove that both are left
adjoints.

\bp \label{factor_rea} Let $K$ be a precubical set. Then there are the
natural isomorphisms of flows $|\mathcal{L}_\mathcal{A}(K)|_{bad} \iso |K|_{bad}$
and $|\mathcal{L}_\mathcal{A}(K)| \iso |K|$. \ep

\bpf Since all functors involved in the statement of the proposition
are left adjoint and therefore colimit-preserving, it suffices to
check the isomorphism for $K = \square[n]$. The proof is complete
after Proposition~\ref{carre}. \epf

\begin{cor} \label{badaussi} For all $n\geq 3$, the inclusion
  $\de\mathcal{A}[n] \subset \mathcal{A}[n]$ induces an
  isomorphism of flows $|\de\mathcal{A}[n]|_{bad} \iso
  |\mathcal{A}[n]|_{bad}$. \end{cor}

\bpf Since $\mathcal{L}_\mathcal{A}$ is colimit-preserving, one has
$\mathcal{L}_\mathcal{A}(\de\square[n]) \iso
\de\mathcal{A}[n]$. So by Proposition~\ref{factor_rea} and
Theorem~\ref{badsigma}, one obtains $|\de\mathcal{A}[n]|_{bad}
\iso |\de\square[n]|_{bad} \iso |\square[n]|_{bad} \iso
|\mathcal{A}[n]|_{bad}$.  \epf

\bp Let $K$ be an $\mathcal{A}$-set. The canonical map $K \rightarrow
\sh_\mathcal{A}(K)$ induces an isomorphism of flows $|K|_{bad} \iso
|\sh_\mathcal{A}(K)|_{bad}$. \ep

\bpf We already know that the map $K \rightarrow \sh_\mathcal{A}(K)$
is obtained by factoring the canonical map $K \rightarrow \mathbf{1}$
as the composite $K \rightarrow \sh_\mathcal{A}(K) \rightarrow
\mathbf{1}$ where $K\rightarrow \sh_\mathcal{A}(K)$ is a relative
$\{\mathcal{A}[p]\sqcup_{\de\mathcal{A}[p]}\mathcal{A}[p]\rightarrow
\mathcal{A}[p], p\geq 2\}$-cell complex and the map
$\sh_\mathcal{A}(K) \rightarrow\mathbf{1}$ satisfies the right lifting
property with respect to the same set of morphisms. So the map
$|K|_{bad} \rightarrow |\sh_\mathcal{A}(K)|_{bad}$ is a relative
$\{|\mathcal{A}[2]|_{bad}\sqcup_{|\de\mathcal{A}[2]|_{bad}}|\mathcal{A}[2]|_{bad}\rightarrow
|\mathcal{A}[2]|_{bad}\}$-cell complex by
Corollary~\ref{badaussi}. Figure~\ref{cube2} explains why the map of
flows
$|\mathcal{A}[2]|_{bad}\sqcup_{|\de\mathcal{A}[2]|_{bad}}|\mathcal{A}[2]|_{bad}\rightarrow
|\mathcal{A}[2]|_{bad}$ is in fact an
isomorphism.~\footnote{Intuitively, adding an algebraic relation is an
  idempotent operation.} Hence the proof is complete. \epf

The commutative diagram of flows of Figure~\ref{conclusion} concludes
the section. It proves that labelled precubical sets and labelled
$\mathcal{A}$-sets are equivalent from a directed algebraic
topological point of view, $K$ being any labelled precubical set.

\begin{figure}
\[
\xymatrix{
 |K| \fr{\iso} \fd{}&& |\mathcal{L}_\mathcal{A}(K)|\fd{} &&\\
&&&&\\
 |!\Sigma| \fr{\iso}\fd{} && |\mathcal{L}_\mathcal{A}(!\Sigma)|\fd{} \fr{} &&
|\sh_\mathcal{A}\mathcal{L}_\mathcal{A}(!\Sigma)| \fd{}\\
&&&&\\
?\Sigma\iso |!\Sigma|_{bad} \fr{\iso} &&
|\mathcal{L}_\mathcal{A}(!\Sigma)|_{bad}\fr{\iso} &&
|\sh_\mathcal{A}\mathcal{L}_\mathcal{A}(!\Sigma)|_{bad}}
\]
\caption{Labelled precubical sets and labelled $\mathcal{A}$-sets
  equivalent from the directed algebraic topological point of view}
\label{conclusion}
\end{figure}

\subsection*{An application}

We give now a small application of the notion of geometric realization
of labelled $\mathcal{A}$-set which will be reused later. The
following proposition could of course be proved without using the
topological material of this section. However, the proof would be more
complicated (see the proof of \cite[Theorem~7.1]{ccsprecub}).

\bp \label{maximal} Let $\ell : \mathcal{A}[p] \rightarrow
\sh_\mathcal{A}\mathcal{L}_\mathcal{A}(!\Sigma)$ be a full labelled
$p$-dimensional $\mathcal{A}$-cube with $p\geq 2$. Then there exists
$(a_1,\dots,a_p)\in \Sigma^p$ such that for every maximal path
$(c_1,\dots,c_p)$ of $\mathcal{A}[p]$, i.e. for any $p$-tuple of
$1$-cubes of $\mathcal{A}[p]$ with $\de_1^0(c_1) = (0,\dots,0)$,
$\de_1^1(c_i) = \de_1^0(c_{i+1})$ for $1\leq i \leq p-1$ and
$\de_1^1(c_{p}) = (1,\dots,1)$, one has $\ell(c_1) * \dots * \ell(c_p)
= a_1 * \dots * a_p$. \ep

\bpf Let $(c_1,\dots,c_p)$ and $(c'_1,\dots,c'_p)$ be two maximal
paths. Since there is a unique morphism from
$(\widehat{0},\dots,\widehat{0})$ to $(\widehat{1},\dots,\widehat{1})$
in $|\mathcal{A}[p]|_{bad}$ (this is the key point !), one has
$\ell(c_1)*\dots*\ell(c_p) = \ell(c'_1)*\dots*\ell(c'_p)$ in the flow
$|\sh_\mathcal{A}\mathcal{L}_\mathcal{A}(!\Sigma)|_{bad}$.  But the
semigroup $\P
(|\sh_\mathcal{A}\mathcal{L}_\mathcal{A}(!\Sigma)|_{bad}) \iso
\P(?\Sigma)$ is the free commutative semigroup generated by the
elements of $\Sigma$.  Hence the result.  \epf

\section{Labelled coskeleton over a category of cubes}
\label{cosk}

We give in this section the generalization of the notion of labelled
coskeleton to any category of labelled precubical sets. The particular
case $\mathcal{A} = \square$ will give back the situation of
\cite{ccsprecub}. The unlabelled version, i.e. when $\Sigma =
\{\tau\}$ is the classical coskeleton functor, right adjoint to the
truncation functor \cite{Brown_cube}.

\subsection*{The unlabelled case}

\bp \label{i1} Let $n\geq 0$.
\begin{enumerate}
\item The functor $K\mapsto K_{\leq n}$ from
  $\mathcal{A}_{n+1}^{op}\set$ to $\mathcal{A}_n^{op}\set$ has a right
  adjoint denoted by $\cosk^\mathcal{A}_{n,n+1}:\mathcal{A}_n^{op}\set \rightarrow
  \mathcal{A}_{n+1}^{op}\set$. There is an inclusion of presheaves
  \[K\subset \cosk^\mathcal{A}_{n,n+1}(K)\] natural with respect to
  the $n$-dimensional $\mathcal{A}$-set $K$. This inclusion induces
  the isomorphism $K \iso \cosk^\mathcal{A}_{n,n+1}(K)_{\leq n}$.
\item The functor $K\mapsto K_{\leq n}$ from $\mathcal{A}^{op}\set$ to
  $\mathcal{A}_n^{op}\set$ has a right adjoint denoted by
  $\cosk^\mathcal{A}_n:\mathcal{A}_n^{op}\set \rightarrow
  \mathcal{A}^{op}\set$. There is an inclusion of presheaves $K\subset
  \cosk^\mathcal{A}_n(K)$ natural with respect to the $n$-dimensional
  $\mathcal{A}$-set $K$. This inclusion induces the isomorphism $K
  \iso \cosk^\mathcal{A}_n(K)_{\leq n}$.
\item Let $\cosk^\mathcal{A}_{n,n+p} = \cosk^\mathcal{A}_{n+p-1,n+p}
  \circ \dots \circ \cosk^\mathcal{A}_{n,n+1} \circ
  \cosk^\mathcal{A}_{n,n}$ where the functor \[\cosk^\mathcal{A}_{n,n}
  : \mathcal{A}_n^{op}\set \rightarrow \mathcal{A}_{n}^{op}\set\]
  denotes the identity functor. Then there is an isomorphism of
  functors \[\cosk^\mathcal{A}_n\iso \liminj
  \cosk^\mathcal{A}_{n,n+p}.\]
\end{enumerate}
\ep

\bpf Let us prove the first assertion. The functor $K\mapsto K_{\leq
  n}$ from $\mathcal{A}_{n+1}^{op}\set$ to $\mathcal{A}_n^{op}\set$ is
induced by the inclusion of categories $\mathcal{A}_n^{op} \subset
\mathcal{A}_{n+1}^{op}$. Thus, the right adjoint is obtained by taking
the right Kan extension along $\mathcal{A}_n^{op} \subset
\mathcal{A}_{n+1}^{op}$. The isomorphism of presheaves $K_{\leq n}
\iso K$ for an $n$-dimensional $\mathcal{A}$-set $K$ induces by
adjunction a natural map $K \rightarrow
\cosk^\mathcal{A}_{n,n+1}(K)$. Let $p\leq n$. There is a bijection
$\mathcal{A}^{op}\set(\mathcal{A}[p],K) \iso
\mathcal{A}^{op}\set(\mathcal{A}[p],\cosk^\mathcal{A}_{n,n+1}(K))$
because of the isomorphism $\mathcal{A}[p]_{\leq n} \iso
\mathcal{A}[p]$. Hence we obtain the desired inclusion. The proof of
the second assertion is similar to the above proof. The third
assertion is obvious. \epf

\bd \label{def_un_shell} Let $K$ be an $\mathcal{A}$-set. An
$(n+1)$-cube of $\cosk^\mathcal{A}_{n,n+1}(K_{\leq n})$, i.e. a map
$\de\mathcal{A}[n+1] \rightarrow K$ is called a {\rm $n$-dimensional
  shell or $n$-shell} of $K$. \ed

\subsection*{The labelled case}

Before giving the labelled version of Proposition~\ref{i1}, let us
prove the following general categorical fact.

\bp \label{relatif} Let $L:\C \leftrightarrows \D:R$ be a categorical
adjunction where $L$ is the left adjoint and $R$ the right one. Let us
suppose that $\C$ has all pullbacks. Let $A$ be an object of
$\C$. Then the functor $L_A:\C\ddownarrow A \rightarrow \D\ddownarrow
L(A)$ defined by $L_A(X\rightarrow A) := L(X)\rightarrow L(A)$ has a
right adjoint $R_A$ defined by the following pullback diagram of $\C$:
\[
\xymatrix{
R_A(Y) \cartesien \fd{} \fr{} && R(Y) \fd{}\\
&&\\
A \fr{} && R(L(A))}
\]
where the map $A\rightarrow R(L(A))$ is the unit of the
adjunction. \ep

Note that we are going to use Proposition~\ref{relatif} with $\C$ and
$\D$ locally presentable. In this situation, the categories
$\C\ddownarrow A$ and $\D\ddownarrow L(A)$ are both locally
presentable as well by \cite[Proposition~1.57]{MR95j:18001}. In
particular, the category $\C\ddownarrow A$ has a generator and is
co-wellpowered. The functor $L_A:\C\ddownarrow A \rightarrow
\D\ddownarrow L(A)$ is colimit-preserving since $L$ is
colimit-preserving. So by the opposite of the Special Adjoint Functor
Theorem, the functor $L_A$ has a right adjoint.

\bpf Let $X\rightarrow A$ be an object of $\C\ddownarrow A$. Let
$Y\rightarrow L(A)$ be an object of $\D\ddownarrow L(A)$. There is a
bijection between the commutative diagrams of the form 
\[
\xymatrix{
X\fr{} \fd{} && R_A(Y) \fd{} \\
&& \\
A \ar@{=}[rr] && A}
\]
and the  commutative diagrams of the form 
\[
\xymatrix{
X\fr{} \fd{} && R(Y) \fd{} \\
&& \\
A \fr{} && R(L(A))}
\]
because of the universal property of pullback. And there is a
bijection between the latter diagrams and the commutative diagrams of
the form
\[
\xymatrix{
L(X)\fr{} \fd{} && Y \fd{} \\
&& \\
L(A) \ar@{=}[rr] && L(A)}
\]
by universality of adjunction. Hence the result. \epf

Here is now the labelled analogue of Proposition~\ref{i1}.

\bp \label{i1_relatif} Let $n\geq 0$.
\begin{enumerate}
\item The functor $K\mapsto K_{\leq n}$ from
  $\mathcal{A}_{n+1}^{op}\set \ddownarrow
  \sh_\mathcal{A}\mathcal{L}_\mathcal{A}(!\Sigma)$ to
  $\mathcal{A}_n^{op}\set \ddownarrow
  \sh_\mathcal{A}\mathcal{L}_\mathcal{A}(!\Sigma)$ has a right adjoint
  denoted by $\cosk_{n,n+1}^{\mathcal{A},\Sigma}:\mathcal{A}_n^{op}\set \ddownarrow
  \sh_\mathcal{A}\mathcal{L}_\mathcal{A}(!\Sigma) \rightarrow
  \mathcal{A}_{n+1}^{op}\set \ddownarrow
  \sh_\mathcal{A}\mathcal{L}_\mathcal{A}(!\Sigma)$. There is an
  inclusion of presheaves
  \[K\subset \cosk_{n,n+1}^{\mathcal{A},\Sigma}(K)\] natural with
  respect to the $n$-dimensional labelled $\mathcal{A}$-set $K$. This
  inclusion induces the isomorphism $K \iso
  \cosk_{n,n+1}^{\mathcal{A},\Sigma}(K)_{\leq n}$.
\item The functor $K\mapsto K_{\leq n}$ from $\mathcal{A}^{op}\set
  \ddownarrow \sh_\mathcal{A}\mathcal{L}_\mathcal{A}(!\Sigma)$ to
  $\mathcal{A}_n^{op}\set \ddownarrow
  \sh_\mathcal{A}\mathcal{L}_\mathcal{A}(!\Sigma)$ has a right adjoint
  denoted by $\cosk_n^{\mathcal{A},\Sigma}:\mathcal{A}_n^{op}\set
  \ddownarrow \sh_\mathcal{A}\mathcal{L}_\mathcal{A}(!\Sigma)
  \rightarrow \mathcal{A}^{op}\set \ddownarrow
  \sh_\mathcal{A}\mathcal{L}_\mathcal{A}(!\Sigma)$. There is an
  inclusion of presheaves $K\subset \cosk_n^{\mathcal{A},\Sigma}(K)$
  natural with respect to the $n$-dimensional labelled
  $\mathcal{A}$-set $K$. This inclusion induces the isomorphism $K
  \iso \cosk_n^{\mathcal{A},\Sigma}(K)_{\leq n}$.
\item Let $\cosk_{n,n+p}^{\mathcal{A},\Sigma} =
  \cosk_{n+p-1,n+p}^{\mathcal{A},\Sigma} \circ \dots \circ
  \cosk_{n,n+1}^{\mathcal{A},\Sigma} \circ
  \cosk_{n,n}^{\mathcal{A},\Sigma}$ where the functor
  $\cosk_{n,n}^{\mathcal{A},\Sigma} : \mathcal{A}_n^{op}\set
  \ddownarrow \sh_\mathcal{A}\mathcal{L}_\mathcal{A}(!\Sigma)
  \rightarrow \mathcal{A}_{n}^{op}\set \ddownarrow
  \sh_\mathcal{A}\mathcal{L}_\mathcal{A}(!\Sigma)$ denotes the
  identity functor. Then there is an isomorphism of functors
  $\cosk_n^{\mathcal{A},\Sigma}\iso \liminj
  \cosk_{n,n+p}^{\mathcal{A},\Sigma}$.
\end{enumerate}
\ep

\bpf We note that the categories $\mathcal{A}_n^{op}\set \ddownarrow
\sh_\mathcal{A}\mathcal{L}_\mathcal{A}(!\Sigma)$ and
$\mathcal{A}_n^{op}\set \ddownarrow
(\sh_\mathcal{A}\mathcal{L}_\mathcal{A}(!\Sigma))_{\leq n}$ are
isomorphic. So the theorem is a consequence of Proposition~\ref{i1}
and Proposition~\ref{relatif}. \epf

Note that for every $n\geq 1$ and for every $n$-dimensional labelled
$\mathcal{A}$-set $K$, one has the pullback diagram of
$\mathcal{A}$-sets
\[
\xymatrix{
  \cosk_n^{\mathcal{A},\Sigma}(K) \fr{}\fd{}\cartesien && \cosk^\mathcal{A}_n(K) \fd{}\\
  &&\\
  \sh_\mathcal{A}\mathcal{L}_\mathcal{A}(!\Sigma) \fr{} && \cosk^\mathcal{A}_n((\sh_\mathcal{A}\mathcal{L}_\mathcal{A}(!\Sigma))_{\leq n}).}
\] 
Intuitively, this means that the labelled coskeleton functor keeps
from the unlabelled one only the shells which are compatibly
labelled. For example, the boundary of a square is compatibly labelled
if and only if opposite sides are labelled in the same way.

\bd Let $K$ be a labelled $\mathcal{A}$-set. An $(n+1)$-cube of
$\cosk_{n,n+1}^{\mathcal{A},\Sigma}(K_{\leq n})$ is called a {\rm
  labelled $n$-dimensional shell or $n$-shell} of $K$. \ed

The following proposition generalizes \cite[Definition~3.12 and
Proposition~3.13]{ccsprecub}.

\bp \label{explicitation} Let $K$ be a labelled $\mathcal{A}$-set. The
set of labelled $n$-dimensional shells of $K$ is in bijection with the
set of commutative diagrams of the form
\[
\xymatrix{
  \de\mathcal{A}[n+1] \fr{}\fd{} && K \fd{}\\
  &&\\
  \mathcal{A}[n+1]\fr{}&&
  \sh_\mathcal{A}\mathcal{L}_\mathcal{A}(!\Sigma).}
\] 
\ep 

Since $\sh_\mathcal{A}\mathcal{L}_\mathcal{A}(!\{\tau\})$ is the
terminal $\mathcal{A}$-set by Proposition~\ref{1objet}, the case
$\Sigma = \{\tau\}$ coincides with the unlabelled notion of
Definition~\ref{def_un_shell}.

\bpf Let $\mathcal{A}[n+1] \rightarrow
\cosk_{n,n+1}^{\mathcal{A},\Sigma}(K_{\leq n})$ be a labelled
$n$-shell of $K$. By adjunction, one obtains the commutative diagram
of labelled $\mathcal{A}$-sets
\[
\xymatrix{ \de\mathcal{A}[n+1] \fr{}\fd{} && K_{\leq
    n} \fd{}\\
  &&\\
  \mathcal{A}[n+1]\fr{}&& \cosk_{n,n+1}^{\mathcal{A},\Sigma}(K_{\leq n}).}
\] 
By composing with the labelling map
$\cosk_{n,n+1}^{\mathcal{A},\Sigma}(K_{\leq n}) \rightarrow
\sh_\mathcal{A}\mathcal{L}_\mathcal{A}(!\Sigma)$, one obtains the
commutative diagram of $\mathcal{A}$-sets
\[
\xymatrix{
  \de\mathcal{A}[n+1] \fr{}\fd{} && K \fd{}\\
  &&\\
  \mathcal{A}[n+1]\fr{}&& \sh_\mathcal{A}\mathcal{L}_\mathcal{A}(!\Sigma).}
\] 
Conversely, from such a diagram, one obtains the commutative diagram
of $\mathcal{A}$-sets
\[
\xymatrix{
  \mathcal{A}[n+1]_{\leq n} \fr{}\fd{} && K_{\leq n} \fd{}\\
  &&\\
  \sh_\mathcal{A}\mathcal{L}_\mathcal{A}(!\Sigma)\ar@{=}[rr]&&
  \sh_\mathcal{A}\mathcal{L}_\mathcal{A}(!\Sigma),}
\]
hence the result by adjunction.  \epf

\part{Mathematical treatment}
\label{math}

\section{Shell-complete category of cubes}
\label{forma}

The purpose of this combinatorial section is to address the following
question. Is it possible to find a category of cubes $\mathcal{A}$
such that $\cosk_1^{\mathcal{A},\Sigma}(\mathcal{A}[a_1, \dots,
a_n]_{\leq 1})$ is exactly the labelled $n$-cube $\mathcal{A}[a_1,
\dots, a_n]$ for every $n \geq 0$ and every $a_1, \dots, a_n \in
\Sigma$ ? Let us repeat one more time that there is always a strict
inclusion $\square[a_1, \dots, a_n] \subset
\cosk_1^{\square,\Sigma}(\square[a_1, \dots, a_n]_{\leq 1})$ for every
$n\geq 2$ by \cite[Proposition~3.15]{ccsprecub} and that this is the
reason for introducing in \cite{ccsprecub} the non-functorial
subobject of $\cosk_1^{\square,\Sigma}(\square[a_1, \dots, a_n]_{\leq
  1})$ called the labelled directed coskeleton of $\square[a_1, \dots,
a_n]_{\leq 1}$ (see Definition~\ref{def_directed}). For $\Sigma =
\{\tau\}$, i.e. for the unlabelled case, the previous equality reduces
to finding a category of cubes $\mathcal{A}$ such that
$\cosk_1^\mathcal{A}(\mathcal{A}[n]_{\leq 1}) \iso \mathcal{A}[n]$ for
every $n \geq 0$.  Such a category $\mathcal{A}$ will be called a
shell-complete category of cubes. We will see in
Theorem~\ref{but_math} that such a category of cubes answers the
question above. We will see in Theorem~\ref{identification_smallest}
that there exists one and only one such a category of cubes.

\subsection*{Definition and elementary properties}

\bp \label{yoneda_cube} Let $p,q\geq 0$. The natural
bijection \[\mathcal{A}^{op}\set(\mathcal{A}[p],\mathcal{A}[q]) \iso
\mathcal{A}([p],[q])\] induced by the mapping $f\mapsto
f_p(\id_{[p]})$ given by the Yoneda lemma takes $f:\mathcal{A}[p]
\rightarrow \mathcal{A}[q]$ to $f_0:[p]\iso \mathcal{A}([0],[p])
\rightarrow \mathcal{A}([0],[q]) \iso [q]$. \ep

\bpf Let $f:\mathcal{A}[p] \rightarrow \mathcal{A}[q]$ be a map of
$\mathcal{A}^{op}\set$. The map $\delta_p^{\epsilon_p}\dots
\delta_1^{\epsilon_1}:[0] \rightarrow [p]$ induces a commutative
square of sets
\[
\xymatrix{
\mathcal{A}([p],[p]) \fr{f_p}\fd{(\delta_p^{\epsilon_p}\dots
\delta_1^{\epsilon_1})^{*}} && \mathcal{A}([p],[q])\fd{(\delta_p^{\epsilon_p}\dots
\delta_1^{\epsilon_1})^{*}} \\
&&\\
\mathcal{A}([0],[p])\fr{f_0} && \mathcal{A}([0],[q])}
\]
for any $\epsilon_1,\dots,\epsilon_p\in \{0,1\}$ since $\mathcal{A}$
is a category of cubes.  So \[f_0(\delta_p^{\epsilon_p}\ldots
\delta_1^{\epsilon_1})^{*}(\id_{[p]}) = (\delta_p^{\epsilon_p}\ldots
\delta_1^{\epsilon_1})^{*} (f_p(\id_{[p]})).\] Therefore $f_0 =
f_p(\id_{[p]})$.  \epf

The following proposition motivates the notion of shell-complete
category of cubes.

\bp \label{tjrs_inc} For any $q\geq 0$, the canonical map
$\mathcal{A}[q] \rightarrow \cosk^\mathcal{A}_1(\mathcal{A}[q]_{\leq
  1})$ induced by the isomorphism $\mathcal{A}[q]_{\leq 1} \iso
\mathcal{A}[q]_{\leq 1}$ is an inclusion of presheaves. For $q = 0$ or
$q = 1$, this inclusion is always an equality for any category of
cubes $\mathcal{A}$. \ep

\bpf Let $x$ and $y$ be two $k$-cubes of $\mathcal{A}[q]$ having the
same image by the map \[\mathcal{A}[q] \longrightarrow
\cosk^\mathcal{A}_1(\mathcal{A}[q]_{\leq 1}).\] So one has the
commutative diagram of $\mathcal{A}$-sets
\[
\xymatrix{ \mathcal{A}[k] \ar@<.8ex>[rr]^{x}\ar@<-.8ex>[rr]_{y}&&
  \mathcal{A}[q] \fr{} && \cosk^\mathcal{A}_1(\mathcal{A}[q]_{\leq
    1}).}\] By adjunction, one obtains the commutative diagram of
$\mathcal{A}$-sets\[ \xymatrix{ \mathcal{A}[k]_{\leq 1}
  \ar@<.8ex>[rr]^{x_{\leq 1}}\ar@<-.8ex>[rr]_{y_{\leq 1}}&&
  \mathcal{A}[q]_{\leq 1} \fr{\iso} && \mathcal{A}[q]_{\leq 1}.}\] In
particular, the two set maps $x_0,y_0:[k] \rightrightarrows [q]$ are
equal.  Thus, by Proposition~\ref{yoneda_cube}, one obtains $x=y$. The
last assertion is a consequence of Proposition~\ref{i1}. \epf

Hence the definition:

\bd \label{shell-complete} A category of cubes $\mathcal{A}$ is {\rm
  shell-complete} if for every $p\geq 2$, the canonical inclusion of
presheaves $\mathcal{A}[p] \subset
\cosk^\mathcal{A}_1(\mathcal{A}[p]_{\leq 1})$ is an isomorphism.  \ed

The category of cubes $\square$ is of course not shell-complete by
\cite[Proposition~3.15]{ccsprecub}. For example, the precubical set
$\cosk^\square_1(\square[2]_{\leq 1})$ contains the $2$-cube
$x:\square[2] \rightarrow \cosk^\square_1(\square[2]_{\leq 1})$
corresponding by adjunction to the map $\de x:\de\square[2] =
\square[2]_{\leq 1} \rightarrow \square[2]_{\leq 1}$ characterized by
$x_0(\epsilon_1,\epsilon_2) = (\epsilon_2,\epsilon_1)$.  It is not a
$2$-cube of $\square[2]$ since the only $2$-cube of the precubical set
$\square[2]$ is the identity of $[2]$.

In general, for any $p,q\geq 2$, there exists at most one lift
$\overline{x}$ in the commutative diagram of solid arrows
\[
\xymatrix{
\de\mathcal{A}[p] \fr{x}\fd{} && \mathcal{A}[q] \fd{}\\
&&\\
\mathcal{A}[p] \fr{}\ar@{-->}[rruu]^-{\overline{x}} && \mathbf{1},
}
\]
where $\mathbf{1}$ is the terminal object. Indeed, by
Proposition~\ref{yoneda_cube}, the bijection of sets
\[\mathcal{A}^{op}\set(\mathcal{A}[p],\mathcal{A}[q]) \iso
\mathcal{A}([p],[q])\] takes $\overline{x}$ to
$x_0$. Shell-completeness means that this lift always exists.

\bth \label{condcompl} Let $\mathcal{A}$ be a category of cubes.  The
following conditions are equivalent:
\begin{enumerate}
\item The category $\mathcal{A}$ is shell-complete.
\item For any $p,q\geq 2$, for any map $x:\de\mathcal{A}[p]
  \rightarrow \mathcal{A}[q]$, the set map $x_0:[p] \rightarrow [q]$
  belongs to $\mathcal{A}$.
\item For any $p,q\geq 2$, any map $x:\de\mathcal{A}[p] \rightarrow
  \mathcal{A}[q]$ factors uniquely as a composite $x:\de\mathcal{A}[p]
  \rightarrow \mathcal{A}[p] \rightarrow \mathcal{A}[q]$.
\end{enumerate}
\eth 

\bpf Let us prove the implication $(1) \Longrightarrow (2)$.  Let
$x:\de\mathcal{A}[p] \rightarrow \mathcal{A}[q]$ be a map of
$\mathcal{A}$-sets with $p,q\geq 2$. One can suppose that $p\leq q$ by
Proposition~\ref{ex}. Then $x$ factors (uniquely) as a composite
\[x : \de\mathcal{A}[p] =\mathcal{A}[p]_{\leq p-1} \longrightarrow
\mathcal{A}[q]_{\leq p-1} \longrightarrow \mathcal{A}[q].\] One has
the isomorphisms 
\[\mathcal{A}[q]_{\leq p-1}\iso\cosk^\mathcal{A}_1(\mathcal{A}[q]_{\leq 1})_{\leq
  p-1} \iso \cosk^\mathcal{A}_{1,p-1}(\mathcal{A}[q]_{\leq 1})\] since
$\mathcal{A}$ is shell-complete and by Proposition~\ref{i1}. So $x$
factors as a composite \[x:\de\mathcal{A}[p] \longrightarrow
\mathcal{A}[p] \longrightarrow
\cosk^\mathcal{A}_{p-1,p}(\cosk^\mathcal{A}_{1,p-1}(\mathcal{A}[q]_{\leq
  1})) = \mathcal{A}[q]_{\leq p} \longrightarrow \mathcal{A}[q].\] So
$x_0:[p] \rightarrow [q]$ is a morphism of $\mathcal{A}$ by
Proposition~\ref{yoneda_cube}.

Let us prove now the implication $(2) \Longrightarrow
(1)$. Proposition~\ref{i1} and Proposition~\ref{tjrs_inc} imply that
there is an inclusion of presheaves \[\mathcal{A}[q]_{\leq p} \subset
(\cosk^\mathcal{A}_{1}(\mathcal{A}[q]_{\leq 1}))_{\leq p} =
\cosk^\mathcal{A}_{1,p}(\mathcal{A}[q]_{\leq 1})\] for any $p\geq
1$. This inclusion is trivially an equality for $p=1$. Let us prove by
induction on $p$ that this inclusion is an equality. This will
establish the shell-completeness of $\mathcal{A}$. Let us suppose the
equality proved for $p\geq 1$. Let $x:\mathcal{A}[p+1] \rightarrow
\cosk^\mathcal{A}_{1,p+1}(\mathcal{A}[q]_{\leq 1})$ be a
$(p+1)$-dimensional $\mathcal{A}$-cube of
$\cosk^\mathcal{A}_{1,p+1}(\mathcal{A}[q]_{\leq 1})$. By adjunction
and by induction hypothesis, one obtains a map \[\de
x:\mathcal{A}[p+1]_{\leq p} = \de\mathcal{A}[p+1] \longrightarrow
\cosk^\mathcal{A}_{1,p}(\mathcal{A}[q]_{\leq 1}) \iso
\mathcal{A}[q]_{\leq p} \subset \mathcal{A}[q].\] By hypothesis, the
map $x_0:[p+1] \rightarrow [q]$ belongs to $\mathcal{A}$. Thus, by
Proposition~\ref{yoneda_cube}, there is a commutative diagram of
$\mathcal{A}$-sets
\[\xymatrix{ \de\mathcal{A}[p+1] \fd{}\fr{} &&
  \cosk^\mathcal{A}_{1,p}(\mathcal{A}[q]_{\leq 1})
  \iso \mathcal{A}[q]_{\leq p} \fd{\subset} \\
  && \\
  \mathcal{A}[p+1] \fr{} && \mathcal{A}[q]_{\leq p+1}.}
\] 
Hence the equality for $p+1$.  

The equivalence $(2) \Longleftrightarrow (3)$ is a consequence of
Proposition~\ref{yoneda_cube}. \epf

\subsection*{Examples of shell-complete categories of cubes}

\bth \label{nonvide} The category of cubes $\widehat{\square}$
(i.e. the maximal category of cubes containing all
adjacency-preserving maps) is shell-complete.\eth

\bpf Let $x:\de\widehat{\square}[p] \rightarrow \widehat{\square}[q]$
be a morphism of $\widehat{\square}^{op}\set$ with $p,q\geq 2$. For
all $k$ such that $1\leq k\leq p-1$, one has the commutative diagram
of sets
\[
\xymatrix{
  \widehat{\square}([k],[p]) \fr{x_k } \fd{\de_1^{\epsilon_1}\dots \de_k^{\epsilon_k}} && \widehat{\square}([k],[q])\fd{\de_1^{\epsilon_1}\dots \de_k^{\epsilon_k}} \\
  && \\
  \widehat{\square}([0],[p])\fr{x_0} && \widehat{\square}([0],[q])}
\]
for all $\epsilon_1,\dots,\epsilon_k\in\{0,1\}$ since
$x:\de\widehat{\square}[p] \rightarrow \widehat{\square}[q]$ is a map
of $\widehat{\square}^{op}\set$ and where the set map
$\de_1^{\epsilon_1}\dots \de_k^{\epsilon_k}$ is induced by the
morphism $\delta_k^{\epsilon_k}\dots \delta_1^{\epsilon_1}:[0]
\rightarrow [k]$ of $\widehat{\square}$. With $\phi\in
\widehat{\square}([k],[p])$, that means that $x_k(\phi)
(\delta_k^{\epsilon_k}\dots \delta_1^{\epsilon_1}) = x_0(\phi
(\delta_k^{\epsilon_k}\dots \delta_1^{\epsilon_1}))$. Thus, one
obtains $x_k(\phi)(\epsilon_1,\dots,\epsilon_k) =
x_0(\phi(\epsilon_1,\dots,\epsilon_k))$. So $x_k(\phi)=x_0 \phi$ with
the identification $\widehat{\square}([0],[p])\iso [p]$.  Let
$(\epsilon_1,\dots,\epsilon_p)$ and $(\epsilon'_1,\dots,\epsilon'_p)$
be two elements of $[p]$ with $\epsilon_i=\epsilon'_i$ for all $i$ but
one denoted by $i_0$. Suppose moreover that $\epsilon_{i_0} = 0$ and
$\epsilon'_{i_0} = 1$. Since $p \geq 2$, there exists $i_1\in
\{1,\dots,p\} \backslash \{i_0\}$.  Consider
$\delta_{i_1}^{\epsilon_{i_1}}:[p-1] \rightarrow [p]$.  Then
$x_{p-1}(\delta_{i_1}^{\epsilon_{i_1}})=x_0
\delta_{i_1}^{\epsilon_{i_1}}$. The preceding equality applied to
$(\epsilon_1,\dots,\widehat{\epsilon_{i_1}},\dots,\epsilon_p)$ gives
\begin{align*}
  x_0(\epsilon_1,\dots,\epsilon_p) &=x_{p-1}(\delta_{i_1}^{\epsilon_{i_1}})(\epsilon_1,\dots,\widehat{\epsilon_{i_1}},\dots,\epsilon_p) \\
  &< x_{p-1}(\delta_{i_1}^{\epsilon_{i_1}})
  (\epsilon'_1,\dots,\widehat{\epsilon_{i_1}},\dots,\epsilon'_p) \\
  &= x_0(\epsilon'_1,\dots,\epsilon'_p) \end{align*} since the map
$x_{p-1}(\delta_{i_1}^{\epsilon_{i_1}}):[p-1] \rightarrow [q]$ is a
morphism of the small category $\widehat{\square}$. So the set map
$x_0:[p] \rightarrow [q]$ is adjacency-preserving, i.e. it belongs to
the small category $\widehat{\square}$.  Thus, the small category
$\widehat{\square}$ is shell-complete by Theorem~\ref{condcompl}. \epf

\bp \label{correspondance} Let $\mathcal{A}$ be a category of
cubes. Let $p,q\geq 2$.  The set
$\mathcal{A}^{op}\set(\de\mathcal{A}[p],\mathcal{A}[q])$ is equal to
the set of families $(f_i^\alpha:\mathcal{A}[p-1] \rightarrow
\mathcal{A}[q])$ of morphisms of $\mathcal{A}^{op}\set$ with $1\leq i
\leq p$ and $\alpha\in \{0,1\}$ with $(f_j^\beta)_0 \delta_i^\alpha =
(f_i^\alpha)_0 \delta_{j-1}^\beta$ for any $i<j$ and any
$\alpha,\beta\in \{0,1\}$.  \ep

\bpf Let $f:\de\mathcal{A}[p] \rightarrow \mathcal{A}[q]$ be a
morphism of $\mathcal{A}^{op}\set$. The $2p$ inclusions
$\mathcal{A}[p-1] \subset \de\mathcal{A}[p]$ with $1\leq i\leq p$ and
$\alpha\in\{0,1\}$ induce $2p$ maps $f_i^\alpha : \mathcal{A}[p-1]
\subset \de\mathcal{A}[p] \rightarrow \mathcal{A}[q]$ such that
$(f_i^\alpha)_0 = f_0 \delta_i^\alpha$ with $1\leq i \leq p$ and
$\alpha\in \{0,1\}$.  The equalities $(f_j^\beta)_0 \delta_i^\alpha =
(f_i^\alpha)_0 \delta_{j-1}^\beta$ for any $i<j$ and any
$\alpha,\beta\in \{0,1\}$ are then a consequence of the cocubical
relations.

Conversely, let $(f_i^\alpha:\mathcal{A}[p-1] \rightarrow
\mathcal{A}[q])$ be a family of morphisms of $\mathcal{A}^{op}\set$
with $1\leq i \leq p$ and $\alpha\in \{0,1\}$ such that $(f_j^\beta)_0
\delta_i^\alpha = (f_i^\alpha)_0 \delta_{j-1}^\beta$ for any $i<j$ and
any $\alpha,\beta\in \{0,1\}$.  Consider the set map $g:[p]
\rightarrow [q]$ defined by $g(\epsilon_1,\dots,\epsilon_p) :=
(f_p^{\epsilon_p})_0(\epsilon_1,\dots,\epsilon_{p-1})$. Then $g
\delta_p^\alpha = (f_p^\alpha)_0$ by definition of $g$ and for any
$0\leq i<p$, one has \[g
\delta_i^\alpha(\epsilon_1,\dots,\epsilon_{p-1}) =
(f_p^{\epsilon_{p-1}})_0\delta_i^\alpha(\epsilon_1,\dots,\epsilon_{p-2})
= (f_i^\alpha)_0
\delta_{p-1}^{\epsilon_{p-1}}(\epsilon_1,\dots,\epsilon_{p-2}) =
(f_i^\alpha)_0(\epsilon_1,\dots,\epsilon_{p-1})\] for any
$\alpha,\epsilon_{p-1}\in \{0,1\}$ thanks to the cocubical
relations. So one obtains $g \delta_i^\alpha = (f_i^\alpha)_0$ for
$0\leq i\leq p$ and $\alpha\in\{0,1\}$.  The mapping $\phi\mapsto g
\phi$ gives rise for each $0\leq k\leq p-1$ to a set map $g_k
:\de\mathcal{A}[p]_k:=\mathcal{A}([k],[p]) \rightarrow
\mathcal{A}[q]_k:=\mathcal{A}([k],[q])$. For any morphism $\psi:
[k']\rightarrow [k]$ of $\mathcal{A}$ with $0\leq k'\leq k\leq p-1$,
one obtains a diagram of sets
\[
\xymatrix{
\de\mathcal{A}[p]_k \fr{} \fd{} && \mathcal{A}[q]_k \fd{}\\
&& \\
\de\mathcal{A}[p]_{k'} \fr{} && \mathcal{A}[q]_{k'}}
\] 
which is commutative since the two boundaries of the square takes
$\phi\in \mathcal{A}[p]_k$ to $g \phi \psi$.  \epf

\bp \label{inclusion} Let $\mathcal{A}$ and $\mathcal{B}$ be two
categories of cubes such that $\mathcal{A} \subset \mathcal{B}$.  Let
$p,q\geq 2$. Then one has the inclusion
$\mathcal{A}^{op}\set(\de\mathcal{A}[p],\mathcal{A}[q]) \subset
\mathcal{B}^{op}\set(\de\mathcal{B}[p],\mathcal{B}[q])$ by identifying
the maps $f$ with the corresponding set maps $f_0$ from $[p]$ to
$[q]$.  \ep

\bpf This is a corollary of Proposition~\ref{correspondance} and of
the fact that $\mathcal{A}([p-1],[q]) \subset
\mathcal{B}([p-1],[q])$. \epf

\bth There exists a smallest shell-complete category, denoted by
$\widetilde{\square}$. \eth

\bpf Let $(\square^{(i)})_{i\in I}$ be the class of all shell-complete
small categories of cubes. This class is non-empty by
Theorem~\ref{nonvide}, and small since for any $i$, there is the
inclusion $\square^{(i)} \subset \poset$. Consider the small category
$\widetilde{\square} = \bigcap_{i\in I} \square^{(i)}$. Let
$f:\de\widetilde{\square}[p] \rightarrow \widetilde{\square}[q]$ be a
map of $\widetilde{\square}^{op}\set$ with $p,q\geq 2$.  By
Proposition~\ref{inclusion}, the morphism of presheaves $f$ gives rise
for each $i\in I$ to a morphism of presheaves
$f^{(i)}:\de\square^{(i)}[p] \rightarrow \square^{(i)}[q]$. By
Theorem~\ref{condcompl}, $f_0=(f^{(i)})_0$ is a morphism of
$\square^{(i)}$ for each $i\in I$. So by Theorem~\ref{condcompl}
again, the category $\widetilde{\square}$ is shell-complete.  \epf

\subsection*{Some combinatorial lemmas}

Let us recall that $\sigma_i:[n] \rightarrow [n]$ is the set map
defined for $1\leq i\leq n-1$ and $n\geq 2$ by $\sigma_i(\epsilon_1,
\dots, \epsilon_{n}) = (\epsilon_1, \dots,
\epsilon_{i-1},\epsilon_{i+1},\epsilon_{i},
\epsilon_{i+2},\dots,\epsilon_{n})$ (cf. Definition~\ref{def_sym}).

\bp (\cite{MR1988396} p195) \label{sigma_delta} Let $\sigma_i:[n] \rightarrow
[n]$ be the set map defined for $1\leq i\leq n-1$ and $n\geq 2$ by
$\sigma_i(\epsilon_1, \dots, \epsilon_{n}) = (\epsilon_1, \dots,
\epsilon_{i-1},\epsilon_{i+1},\epsilon_{i},
\epsilon_{i+2},\dots,\epsilon_{n})$.  One has the relations
$\sigma_i\delta_j^\alpha=\delta_j^\alpha\sigma_{i-1}$ for $j<i$,
$\sigma_i\delta_j^\alpha=\delta_{i+1}^\alpha$ for $j=i$,
$\sigma_i\delta_j^\alpha=\delta_i^\alpha$ for $j=i+1$ and
$\sigma_i\delta_j^\alpha=\delta_j^\alpha\sigma_{i}$ for $j>i+1$. \ep

\bp \label{m1} $\sigma_i\in \widetilde{\square}$. \ep

\bpf Let us prove by induction on $n\geq 2$ that the set maps
$\sigma_i:[n] \rightarrow [n]$ for $1\leq i\leq n-1$ belong to
$\widetilde{\square}$. The composite map $\de\sigma_1:\de\square[2]
\subset \square[2]\stackrel{\sigma_1}\longrightarrow \square[2]$
induces a map $\de\sigma_1:\de\widetilde{\square}[2] \rightarrow
\widetilde{\square}[2]$ by Proposition~\ref{inclusion} since $\square
\subset \widetilde{\square}$. So $\sigma_1:[2] \rightarrow [2]$ is a
map of $\widetilde{\square}$ by Theorem~\ref{condcompl} since
$\widetilde{\square}$ is shell-complete. Hence the proof is complete
for $n=2$.  Now assume that $n > 2$. By Proposition~\ref{sigma_delta}
and by induction hypothesis, the $2n$ set maps
$\sigma_i\delta_j^\alpha: [n-1] \subset [n] \rightarrow [n]$ belong to
$\widetilde{\square}$.  These $2n$ morphisms of $\widetilde{\square}$
induce a morphism $\de\widetilde{\square}[n] \rightarrow
\widetilde{\square}[n]$ by Proposition~\ref{correspondance}.  So
$\sigma_i: [n] \rightarrow [n]$ belongs to $\widetilde{\square}$ by
shell-completeness.  \epf

To our knowledge, the structure maps introduced below are new.  They
are related to the notion of \emph{connection} in the setting of
\emph{cubical} sets, see \cite{phd-Al-Agl} \cite{MR1929304}; indeed,
with their notation of $\varepsilon_i$ for degeneracies and
$\Gamma_i^\alpha$ for connections, one has
$\Gamma_i^+=\varepsilon_i\gamma_i$ and
$\Gamma_i^-=\varepsilon_{i+1}\gamma_i$. An example of use of these
connections in directed algebraic topology can be found in
\cite{homcat} and \cite{coin}.

\bd \label{def_trans} Let $\gamma_i:[n] \rightarrow [n]$ be the set map
defined for $1\leq i\leq n-1$ and $n\geq 2$ by
\[\gamma_i(\epsilon_1, \dots, \epsilon_{n}) = (\epsilon_1, \dots,
\epsilon_{i-1},\max(\epsilon_{i},\epsilon_{i+1}),
\min(\epsilon_{i},\epsilon_{i+1}),
\epsilon_{i+2},\dots,\epsilon_{n}).\] These maps are called the {\rm
  transverse degeneracy maps}.  \ed

\bp \label{gamma_delta} One has the relations $\gamma_j\delta_i^\alpha
= \delta_i^\alpha \gamma_j$ for $j<i-1$, $\gamma_j\delta_i^\alpha =
\delta_i^\alpha \gamma_{j-1}$ for $j\geq i+1$,
$\gamma_j\delta_i^\alpha = \delta_{i-\alpha}^\alpha$ for $j=i-1$ and
$\gamma_j\delta_i^\alpha = \delta_{i+1-\alpha}^\alpha$ for $j=i$.  \ep

\bpf The relation $\gamma_j\delta_i^\alpha = \delta_i^\alpha \gamma_j$
for $j<i-1$ is obvious.  One has
\[\gamma_j\delta_i^\alpha(\epsilon_1,\dots,\epsilon_{n-1}) =
\gamma_j(\epsilon_1,\dots,\epsilon_{i-1},\alpha,\epsilon_i,\dots,\epsilon_{n-1})
= \delta_i^\alpha\gamma_{j-1}(\epsilon_1,\dots,\epsilon_{n-1})\] for
$j\geq i+1$. For $j=i-1$, one has
\[\gamma_j\delta_i^1(\epsilon_1,\dots,\epsilon_{n-1}) =
\gamma_j(\epsilon_1,\dots,\epsilon_{i-1},1,\epsilon_i,\dots,\epsilon_{n-1})
= \delta_{i-1}^1(\epsilon_1,\dots,\epsilon_{n-1})\] and 
\[\gamma_j\delta_i^0(\epsilon_1,\dots,\epsilon_{n-1}) =
\gamma_j(\epsilon_1,\dots,\epsilon_{i-1},0,\epsilon_i,\dots,\epsilon_{n-1})
= \delta_{i}^0(\epsilon_1,\dots,\epsilon_{n-1}).\] Finally for $j=i$,
one has \[\gamma_j\delta_i^1(\epsilon_1,\dots,\epsilon_{n-1}) =
\gamma_j(\epsilon_1,\dots,\epsilon_{i-1},1,\epsilon_i,\dots,\epsilon_{n-1})
= \delta_i^1(\epsilon_1,\dots,\epsilon_{n-1})\] and
\[\gamma_j\delta_i^0(\epsilon_1,\dots,\epsilon_{n-1}) = 
\gamma_j(\epsilon_1,\dots,\epsilon_{i-1},0,\epsilon_i,\dots,\epsilon_{n-1})
= \delta_{i+1}^0(\epsilon_1,\dots,\epsilon_{n-1}).\] 
\epf

\bp \label{m2} $\gamma_i \in \widetilde{\square}$. \ep

\bpf The proof is mutatis mutandis the one of Proposition~\ref{m1}.
\epf

\bp \label{decomposition_distance} Let $0\leq m\leq n$. Every
adjacency-preserving (resp. adjacency-preser\-ving one-to-one) map
$f:[m] \rightarrow [n]$ factors uniquely as a composite $[m]
\stackrel{\psi}\longrightarrow [m] \stackrel{\phi}\longrightarrow [n]$
with $\phi\in \square$ and $\psi$ adjacency-preserving
(resp. adjacency-preserving one-to-one).  \ep

Note that by a cardinality argument, if $\psi:[m] \rightarrow [m]$ is
one-to-one, then it is bijective.

\bpf One has $d(f(0,\dots, 0),f(1,\dots, 1)) = m$. So by
Proposition~\ref{distance}, $f([m])$ is an $m$-subcube of $[n]$. So
the assertion is a consequence of
Proposition~\ref{example_preservation_distance}.  \epf

\subsection*{The uniqueness and a negative result}

\bth \label{identification_smallest} The category of cubes
$\widehat{\square}$ is the unique shell-complete category of
cubes. \eth

\bpf It suffices to prove that the inclusion $\widetilde{\square}
\subset \widehat{\square}$ is an equality since the category of cubes
$\widehat{\square}$ is the maximal category of cubes.  

For any $p,q\geq 0$, there is the inclusion $\square([p],[q]) \subset
\widetilde{\square}([p],[q]) \subset \widehat{\square}([p],[q])$ and
one wants to prove the equality $\widetilde{\square}([p],[q]) =
\widehat{\square}([p],[q])$. 

For $p>q$, one has $\square([p],[q]) = \widetilde{\square}([p],[q]) =
\widehat{\square}([p],[q]) = \varnothing$ by Proposition~\ref{ex}. One
has $\widetilde{\square}([0],[q]) = \widehat{\square}([0],[q]) = [q]$
and $\widetilde{\square}([1],[q]) = \widehat{\square}([1],[q])$ by
Proposition~\ref{ex} again. It remains to prove the equality
$\widetilde{\square}([p],[q])= \widehat{\square}([p],[q])$ for $2\leq
p\leq q$ for a fixed $q$ by induction on $p$. 

First of all, let us treat the case $p = 2$. Let $f \in
\widehat{\square}([2],[q])$. By
Proposition~\ref{decomposition_distance}, the set map $f$ factors
uniquely as a composite of set maps $f:[2]
\stackrel{\phi}\longrightarrow [2] \stackrel{\psi}\longrightarrow [q]$
with $\phi \in \widehat{\square}$ and $\psi\in \square$. It is easy to
see that the set $\widehat{\square}([2],[2])$ consists of the four set
maps $\id_{[2]}:(\epsilon_1,\epsilon_2) \mapsto
(\epsilon_1,\epsilon_2)$, $\sigma_1:(\epsilon_1,\epsilon_2) \mapsto
(\epsilon_2,\epsilon_1)$, $\gamma_1:(\epsilon_1,\epsilon_2) \mapsto
(\max(\epsilon_1,\epsilon_2),\min(\epsilon_1,\epsilon_2))$ and
$\sigma_1 \gamma_1:(\epsilon_1,\epsilon_2) \mapsto
(\min(\epsilon_1,\epsilon_2),\max(\epsilon_1,\epsilon_2))$. So
$\widetilde{\square}([2],[2]) = \widehat{\square}([2],[2])$ by
Proposition~\ref{m1} and Proposition~\ref{m2}. Therefore one obtains
$\widetilde{\square}([2],[q]) = \widehat{\square}([2],[q])$ for any
$q\geq 0$.

Let us now treat the case $p \geq 3$. Every set map $f\in
\widehat{\square}([p],[q])$ for $p\geq 3$ gives rise to a map
$x:\widehat{\square}[p] \rightarrow \widehat{\square}[q]$ such that
$x_0 = f$ by Proposition~\ref{yoneda_cube}. By composition, one
obtains a map $\de x:\de\widehat{\square}[p] \subset
\widehat{\square}[p] \rightarrow \widehat{\square}[q]$. By
Proposition~\ref{correspondance}, one obtains $2p$ maps
$x_i^\alpha:\widehat{\square}[p-1] \rightarrow \widehat{\square}[q]$
with $1\leq i \leq p$ and $\alpha\in \{0,1\}$ such that $(x_j^\beta)_0
\delta_i^\alpha = (x_i^\alpha)_0 \delta_{j-1}^\beta$ for any $i<j$ and
any $\alpha,\beta\in \{0,1\}$.  By Proposition~\ref{yoneda_cube}, the
$2p$ set maps $(x_i^\alpha)_0:[p-1] \rightarrow [q]$ for $1\leq i \leq
p$ and $\alpha\in \{0,1\}$ belong to $\widehat{\square}$. So by
induction hypothesis, the latter set maps belong to
$\widetilde{\square}$ as well. By Proposition~\ref{yoneda_cube} again,
one obtains $2p$ maps $y_i^\alpha:\widetilde{\square}[p-1] \rightarrow
\widetilde{\square}[q]$ with $1\leq i \leq p$ and $\alpha\in \{0,1\}$
such that $(y_j^\beta)_0 \delta_i^\alpha = (y_i^\alpha)_0
\delta_{j-1}^\beta$ for any $i<j$ and any $\alpha,\beta\in \{0,1\}$
and such that $(y_i^\alpha)_0=(x_i^\alpha)_0$ for all $1\leq i\leq p$
and $\alpha\in\{0,1\}$. So by Proposition~\ref{correspondance}, one
obtains a map $\de y:\de\widetilde{\square}[p] \rightarrow
\widetilde{\square}[q]$ such that $(\de
y)_0\delta_i^\alpha=y_i^\alpha$ for all $1\leq i\leq p$ and
$\alpha\in\{0,1\}$. By Theorem~\ref{condcompl} and since
$\widetilde{\square}$ is shell-complete, the set map $(\de
y)_0=x_0=f:[p] \rightarrow [q]$ then belongs to
$\widetilde{\square}$. The induction on $p$ is complete.  \epf

\bth \label{pas_assez_gd} The category of cubes $\overline{\square}$
generated by the $\delta_i^\alpha$, $\sigma_i$ and $\gamma_i$
operators is not shell-complete. In other terms, the inclusion of
small categories $\overline{\square} \subset \widehat{\square}$ is
strict. \eth

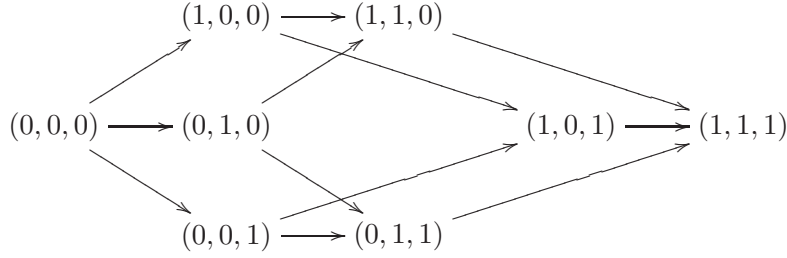
\begin{figure}
\[
\xymatrix{
  & (1,0,0) \ar@{->}[r] \ar@{->}[rrd] & (1,1,0) \ar@{->}[drr] && \\
  (0,0,0) \ar@{->}[r] \ar@{->}[ur] \ar@{->}[dr] & (0,1,0) \ar@{->}[ur]
  \ar@{->}[dr] && (1,0,1) \ar@{->}[r] &
  (1,1,1) \\
  & (0,0,1) \ar@{->}[r] \ar@{->}[rru] & (0,1,1) \ar@{->}[urr] && }
\] 
\caption{Poset of vertices of the $3$-cube}
\label{source}
\end{figure}

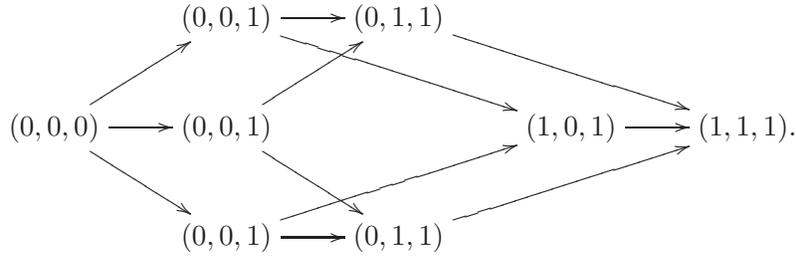
\begin{figure}
\[
\xymatrix{
  & (0,0,1) \ar@{->}[r] \ar@{->}[rrd] & (0,1,1) \ar@{->}[drr] && \\
  (0,0,0) \ar@{->}[r] \ar@{->}[ur] \ar@{->}[dr] & (0,0,1) \ar@{->}[ur]
  \ar@{->}[dr] && (1,0,1) \ar@{->}[r] &
  (1,1,1). \\
  & (0,0,1) \ar@{->}[r] \ar@{->}[rru] & (0,1,1) \ar@{->}[urr] && }
\] 
\caption{Image of the vertices of the $3$-cube by $f$}
\label{but}
\end{figure}

\bpf It suffices to find a morphism of $\widehat{\square}$ which does
not belong to $\overline{\square}$. Consider the set map $f:[3]
\rightarrow [3]$ sending the poset of vertices of the $3$-cube
(Figure~\ref{source}) to the poset depicted in Figure~\ref{but}.

It is clear that $f$ is adjacency-preserving, i.e.  $f\in
\widehat{\square}$. One has
\begin{itemize}
\item $f(0,0,0) = (0,0,0)$, $f(0,1,0) = (0,0,1)$, $f(0,0,1) =
  (0,0,1)$, $f(0,1,1) = (0,1,1)$, so $f\delta_1^0 =
  \delta_1^0\sigma_1\gamma_1$.
\item $f(1,0,0) = (0,0,1)$, $f(1,0,1) = (1,0,1)$, $f(1,1,0) =
  (0,1,1)$, $f(1,1,1) = (1,1,1)$, so $f\delta_1^1 =
  \delta_3^1\sigma_1$.
\item $f(0,0,0) = (0,0,0)$, $f(1,0,0) = (0,0,1)$, $f(0,0,1) =
  (0,0,1)$, $f(1,0,1) = (1,0,1)$, so $f\delta_2^0 =
  \delta_2^0\sigma_1\gamma_1$.
\item $f(0,1,0) = (0,0,1)$, $f(1,1,0) = (0,1,1)$, $f(0,1,1) =
  (0,1,1)$, $f(1,1,1) = (1,1,1)$, so $f\delta_2^1 =
  \delta_3^1\sigma_1\gamma_1$.
\item $f(0,0,0) = (0,0,0)$, $f(1,0,0) = (0,0,1)$, $f(0,1,0) =
  (0,0,1)$, $f(1,1,0) = (0,1,1)$, so $f\delta_3^0 =
  \delta_1^0\sigma_1\gamma_1$.
\item $f(0,0,1) = (0,0,1)$, $f(1,0,1) = (1,0,1)$, $f(0,1,1) =
  (0,1,1)$, $f(1,1,1) = (1,1,1)$, so $f\delta_3^1 = \delta_3^1$.
\end{itemize}
The six set maps $f\delta_1^0 $, $f\delta_1^1$, $f\delta_2^0$,
$f\delta_2^1$, $f\delta_3^0$ and $f\delta_3^1$ belong to
$\overline{\square}$, giving rise to a morphism of presheaves
$\de\overline{\square}[3] \rightarrow \overline{\square}[3]$ by
Proposition~\ref{correspondance}. Any set map $g:[2] \rightarrow [3]$
of $\widehat{\square}$ factors uniquely as a composite $[2]
\stackrel{g_1}\longrightarrow [2] \stackrel{g_2}\longrightarrow [3]$
with $g_1\in \widehat{\square}$ and $g_2\in \square$ by
Proposition~\ref{decomposition_distance}. So the set map
$f_0=\sigma_1\gamma_1 : [2] \rightarrow [2]$ is the unique set map of
$\widehat{\square}$ such that $f\delta_1^0 = \delta_1^0 f_0$. And the
set map $f_1=\sigma_1 : [2] \rightarrow [2]$ is the unique set map of
$\widehat{\square}$ such that $f\delta_1^1 = \delta_3^1 f_1$. Since
$f_0\neq f_1$, the set map $f:[3] \rightarrow [3]$ cannot be a
composite of $\sigma_i : [3] \rightarrow [3]$ and $\gamma_i : [3]
\rightarrow [3]$ with $i=1,2$ by Proposition~\ref{sigma_delta} and
Proposition~\ref{gamma_delta}. Therefore $f\notin \overline{\square}$.
\epf

In fact, we do not know any ``small'' presentation by generators and
relations of the small category $\widehat{\square}$. This is an
interesting and open question. It seems to be related to similar
questions arising in combinatorics: 

\bd An {\rm extremal path} of $[n]$ is a $n$-tuple $(A_1,\dots,A_n)$
of $[n]$ such that $A_1 = 0_n < A_2 < \dots < A_{n-1} < A_n =
1_n$. The set of extremal paths of $[n]$ is denoted by $P_n$. \ed

\begin{nota} Let $e_I$ be the element $(\epsilon_1,\dots,\epsilon_n)$
  of $[n]$ such that $\epsilon_i=1$ if and only if $i\in
  I$. \end{nota}

There is a bijection $p : \Sigma_n \rightarrow P_n$ from the $n$-th
symmetric group $\Sigma_n$ to $P_n$ defined by $p(\sigma) =
(e_\varnothing, e_{\{\sigma(1)\}}, e_{\{\sigma(1),\sigma(2)\}}, \dots,
e_{\{1,\dots,n\}})$.

\bp Let $\sigma \in \Sigma_n$. Let $s_i$ be the transposition $(i\
i+1)$. Then one has the equalities
\[p(\sigma s_i) = (\sigma_i(e_\varnothing),
\sigma_i(e_{\{\sigma(1)\}}), \sigma_i(e_{\{\sigma(1),\sigma(2)\}}),
\dots, \sigma_i(e_{\{1,\dots,n\}}))\] and
\[p(\sigma.\overline{\pi}_i) = (\gamma_i(e_\varnothing),
\gamma_i(e_{\{\sigma(1)\}}), \gamma_i(e_{\{\sigma(1),\sigma(2)\}}),
\dots, \gamma_i(e_{\{1,\dots,n\}}))\] where $\overline{\pi}_i$ are the
elementary increasing bubble sort operators (cf. \cite{HT}) defined by
$\sigma.\overline{\pi}_i = \sigma$ if $\sigma(i) < \sigma(i+1)$ and
$\sigma.\overline{\pi}_i = \sigma s_i$ otherwise.  \ep

\bpf Trivial. \epf 

As a corollary, the monoid $\overline{\square}([n],[n])$ is isomorphic
to the monoid $\langle s_i,\overline{\pi}_i\rangle_{i=1,\dots,n}$ of
set maps from $\Sigma_n$ to itself generated by the operators $s_i$
and $\overline{\pi}_i$. In particular, it satisfies the relations:

\bp \label{relation} Let $n \geq 1$. The monoid of set maps from $[n]$
to itself generated by the $\sigma_i$ and $\gamma_i$ operators
satisfies the relations:
\begin{itemize}
\item $\sigma_i\sigma_i=\id$, $\sigma_i\sigma_j\sigma_i =
  \sigma_j\sigma_i \sigma_j$ for $i=j-1$ and
  $\sigma_i\sigma_j=\sigma_j\sigma_i$ for $i<j-1$ (the Moore relations
  for symmetry operators)
\item $\gamma_i\gamma_i=\gamma_i$, $\gamma_i\gamma_j\gamma_i =
  \gamma_j\gamma_i \gamma_j$ for $i=j-1$ and
  $\gamma_i\gamma_j=\gamma_j\gamma_i$ for $i<j-1$ (the Moore relations
  for transverse degeneracy)
\item $\gamma_j\sigma_i=\sigma_i\gamma_j$ for $j>i+1$ and $j<i-1$,
  $\gamma_i\sigma_i=\gamma_i$, $\sigma_{i+1}\gamma_i\sigma_{i+1} =
  \sigma_{i}\gamma_{i+1}\sigma_{i}$.
\end{itemize}
\ep 

\bpf The Moore relations for symmetry operators are explained for
example in \cite[Theorem~8.1]{MR1988396}.  Let us prove the Moore
relations for transverse degeneracy maps. The relations
$\gamma_i\gamma_i = \gamma_i$ and $\gamma_i\gamma_j =
\gamma_j\gamma_i$ for $i<j-1$ are obvious. One wants to prove that
$\gamma_i\gamma_{i+1}\gamma_i = \gamma_{i+1}\gamma_i\gamma_{i+1}$. It
suffices to prove the identity
$\gamma_1\gamma_{2}\gamma_1(\epsilon_1,\epsilon_2,\epsilon_3) =
\gamma_{2}\gamma_1\gamma_{2}(\epsilon_1,\epsilon_2,\epsilon_3)$. One
has
\begin{align*}
& \gamma_1\gamma_{2}\gamma_1(\epsilon_1,\epsilon_2,\epsilon_3) \\ 
& = \gamma_1\gamma_{2}( \max(\epsilon_1, \epsilon_2), \min(\epsilon_1,
\epsilon_2), \epsilon_3)\\
& = \gamma_1(\max(\epsilon_1, \epsilon_2), \max(\min(\epsilon_1, \epsilon_2), \epsilon_3),
  \min(\epsilon_1, \epsilon_2, \epsilon_3))\\
&=  (\max(\epsilon_1, \epsilon_2, \min(\epsilon_1, \epsilon_2), \epsilon_3), \min(\max(\epsilon_1, \epsilon_2), \max(\min(\epsilon_1, \epsilon_2), \epsilon_3)), \min(\epsilon_1, \epsilon_2,
\epsilon_3))\\
&= (\max(\epsilon_1, \epsilon_2, \epsilon_3), \min(\max(\epsilon_1, \epsilon_2), \max(\min(\epsilon_1, \epsilon_2), \epsilon_3)), \min(\epsilon_1, \epsilon_2,
\epsilon_3))
\end{align*}
and 
\begin{align*}
&\gamma_2\gamma_{1}\gamma_2(\epsilon_1,\epsilon_2,\epsilon_3) \\
&= \gamma_2\gamma_{1}(\epsilon_1, \max(\epsilon_2, \epsilon_3),
\min(\epsilon_2, \epsilon_3))\\
&= \gamma_2(\max(\epsilon_1, \epsilon_2, \epsilon_3), \min(\epsilon_1, \max(\epsilon_2, \epsilon_3)),
   \min(\epsilon_2, \epsilon_3))\\
&= (\max(\epsilon_1, \epsilon_2, \epsilon_3), \max(
   \min(\epsilon_1, \max(\epsilon_2, \epsilon_3)), \min(\epsilon_2, \epsilon_3)),
   \min(\epsilon_1, \max(\epsilon_2, \epsilon_3), \epsilon_2, \epsilon_3))\\
&= (\max(\epsilon_1, \epsilon_2, \epsilon_3), \max(
   \min(\epsilon_1, \max(\epsilon_2, \epsilon_3)), \min(\epsilon_2, \epsilon_3)),
   \min(\epsilon_1, \epsilon_2, \epsilon_3)).
\end{align*}
It remains to check the equality 
\[\min(\max(\epsilon_1, \epsilon_2), \max(\min(\epsilon_1,
\epsilon_2), \epsilon_3)) = \max( \min(\epsilon_1, \max(\epsilon_2,
\epsilon_3)), \min(\epsilon_2, \epsilon_3))\] for any
$(\epsilon_1,\epsilon_2,\epsilon_3)\in\{0,1\}^3$. By distributivity of
$\min$ and $\max$ over each other, one has: 
\begin{align*}
  & \min(\max(\epsilon_1, \epsilon_2), \max(\min(\epsilon_1,
  \epsilon_2), \epsilon_3)) &\\
& = \max(\min(\epsilon_1, \epsilon_2),\min(\epsilon_1,
\epsilon_3),\min(\epsilon_2, \epsilon_3))& \\
  &= \max( \min(\epsilon_1, \max(\epsilon_2, \epsilon_3)),
  \min(\epsilon_2, \epsilon_3)). 
\end{align*}

The proof will be complete by establishing the relations between
transverse degeneracy maps and symmetry operators. The equalities
$\gamma_j\sigma_i=\sigma_i\gamma_j$ for $j>i+1$ and $j<i-1$ and
$\gamma_i\sigma_i=\gamma_i$ are obvious. One wants to prove that
$\sigma_{i+1}\gamma_i\sigma_{i+1} =
\sigma_{i}\gamma_{i+1}\sigma_{i}$. It suffices to prove the identity
$\sigma_{2}\gamma_1\sigma_{2}(\epsilon_1,\epsilon_2,\epsilon_3) =
\sigma_1\gamma_{2}\sigma_1(\epsilon_1,\epsilon_2,\epsilon_3)$. One has
\begin{align*}
   \sigma_{2}\gamma_1\sigma_{2}(\epsilon_1,\epsilon_2,\epsilon_3) &= \sigma_{2}\gamma_1(\epsilon_1,\epsilon_3,\epsilon_2)& \\
  &=
  \sigma_{2}(\max(\epsilon_1,\epsilon_3),\min(\epsilon_1,\epsilon_3),\epsilon_2)
  &\\
  &=
  (\max(\epsilon_1,\epsilon_3),\epsilon_2,\min(\epsilon_1,\epsilon_3))&
\end{align*}
and 
\begin{align*}
  \sigma_1\gamma_{2}\sigma_1(\epsilon_1,\epsilon_2,\epsilon_3)&=\sigma_1\gamma_{2}(\epsilon_2,\epsilon_1,\epsilon_3)&
  \\
  &=
  \sigma_1(\epsilon_2,\max(\epsilon_1,\epsilon_3),\min(\epsilon_1,\epsilon_3))&\\
  &=
  (\max(\epsilon_1,\epsilon_3),\epsilon_2,\min(\epsilon_1,\epsilon_3)).
\end{align*}
\epf

\cite[Conjecture~3.5 and Paragraph~3.1.1]{HT} suggest the
following conjecture:

\begin{conj} \label{conj0} Proposition~\ref{relation} gives a
  presentation by generators and relations of the monoid
  $\overline{\square}([n],[n])$ for every $n\geq 2$. \end{conj}

\subsection*{Functorial interpretation of the labelled directed
  coskeleton}

For $n\geq 2$, and for every $a_1, \dots, a_n \in \Sigma$, the
inclusion $\square[a_1,\dots,a_n] \subset
\cosk^{\square,\Sigma}_1(\square[a_1,\dots,a_n]_{\leq 1})$ is strict
by \cite[Proposition~3.15]{ccsprecub}.  The strictness of the latter
inclusion means that the concurrent execution of $n$ actions always
assemble in $\cosk^{\square,\Sigma}_1(\square[a_1,\dots,a_n]_{\leq
  1})$ to \emph{several} labelled $n$-cubes. To remedy this problem,
the labelled directed coskeleton construction is introduced in
\cite{ccsprecub}. Its main feature is to select one $n$-cube (the
\emph{non-twisted} one) for each multiset of $n$ actions running
concurrently:

\bd \label{def_directed} Let $K$ be a $1$-dimensional labelled
precubical set with $K_0 = [p]$ for some $p\geq 0$. The {\rm labelled
  directed coskeleton} of $K$ is the labelled precubical set
$\COSK^\Sigma(K)$ defined as the subobject of
$\cosk^{\square,\Sigma}_1(K)$ such that:
\begin{itemize}
\item $\COSK^\Sigma(K)_{\leq 1} = \cosk^{\square,\Sigma}_1(K)_{\leq
    1}$
\item for every $n\geq 2$, $x\in \cosk^{\square,\Sigma}_1(K)_n$ is an
  $n$-cube of $\COSK^\Sigma(K)$ if and only if the set map
  $x_0:[n]\rightarrow [p]$ is {\rm non-twisted}, i.e. $x_0: [n]
  \rightarrow [p]$ is a composite\footnote{The factorization is
    necessarily unique.}
\[x_0: [n] \stackrel{\phi}\longrightarrow [q]
\stackrel{\psi}\longrightarrow [p],\] where $\psi$ is a morphism of
the small category $\square$ and where $\phi$ is of the form
\[(\epsilon_1,\dots,\epsilon_{n}) \mapsto
(\epsilon_{i_1},\dots,\epsilon_{i_q})\] such that
$\{1,\dots,n\}\subset \{i_1,\dots,i_q\}$ and such that the first
appearance of $\epsilon_i$ is before the first appearance of
$\epsilon_{i+1}$ in $(\epsilon_{i_1},\dots,\epsilon_{i_q})$ for any
$1\leq i\leq n$ by reading from the left to the right.
\end{itemize} \ed

The fundamental property of the labelled directed coskeleton is then:

\bth \cite[Proposition~3.21]{ccsprecub} \label{reco} Let $n\geq
1$. Let $(a_1,\dots,a_n)\in \Sigma^n$. Then one has the isomorphism of
labelled precubical sets
\[\COSK^\Sigma(\square[a_1,\dots,a_n]_{\leq 1}) \iso
\square[a_1,\dots,a_n].\] \eth

The following theorem gives the functorial interpretation of the
labelled directed coskeleton construction:

\bth \label{but_math} The category of cubes $\mathcal{A} =
\widehat{\square}$ (i.e. the maximal category of cubes containing all
adjacency-preserving maps) is the only category of cubes such that for
every $n\geq 1$ and every $(a_1,\dots,a_n)\in \Sigma^n$, there is the
isomorphism of labelled $\mathcal{A}$-sets
\[\mathcal{L}_\mathcal{A}(\COSK^\Sigma(\square[a_1,\dots,a_n]_{\leq
  1})) \iso \cosk_1^{\mathcal{A},\Sigma}(\mathcal{A}[a_1,\dots,a_n]_{\leq
  1}) (\iso \mathcal{A}[a_1,\dots,a_n]).\]
\eth

\bpf This is a consequence of Theorem~\ref{reco},
Proposition~\ref{carre} and Theorem~\ref{identification_smallest}.
\epf

The commutative diagram of Figure~\ref{conclusion} proves that the
labelled directed coskeleton construction and the labelled transverse
symmetric coskeleton functor are equivalent from a directed algebraic
topological point of view.

\part{Computer-scientific application}
\label{cs}

A short introduction to process algebra can be found in
\cite{MR1365754}. An introduction to CCS (Milner's calculus of
communicating systems \cite{0683.68008}) for mathematicians is
available in \cite{ccsprecub}.

\section{Parallel composition (local case)}
\label{OK}

We want to explain in this section how it is possible to use the
labelled transverse symmetric coskeleton functor to model the parallel
composition in CCS of two labelled cubes representing two higher
dimensional transitions.

\subsection*{The fibered product in CCS}

\begin{nota} $\mathcal{L} := \mathcal{L}_{\widehat{\square}}$, $\sh :=
  \sh_{\widehat{\square}}$ and $\omega =
  \omega_{\widehat{\square}}$. \end{nota}

The set $\Sigma\backslash\{\tau\}$, which may be empty, is now
supposed to be equipped with an involution $a\mapsto \overline{a}$. In
Milner's calculus of communicating systems (CCS) \cite{0683.68008},
which is the only case treated of this paper, one has $a\neq
\overline{a}$.  However, this mathematical hypothesis is not used in
this paper. The involution on $\Sigma\backslash\{\tau\}$ is used only
in Definition~\ref{fibered} of the fibered product (and in the new
definition given with the proof of Theorem~\ref{cavavraiment}) of two
$1$-dimensional labelled (transverse symmetric) precubical sets over
$\Sigma$. For other examples of fibered products over other
synchronization algebras than the one of CCS, cf. \cite{0683.68008}
and \cite{MR1365754}.

\bd \label{fibered} Let $K$ and $L$ be two $1$-dimensional labelled
(transverse symmetric) precubical sets.  The {\rm fibered product} of
$K$ and $L$ {\rm over $\Sigma$} is the $1$-dimensional labelled
precubical set $K\p_\Sigma L$ defined as follows:
\begin{itemize}
\item $(K\p_\Sigma L)_0 = K_0\p L_0$,
\item $(K\p_\Sigma L)_1 = (K_1\p L_0) \sqcup (K_0\p L_1) \sqcup
  \{(x,y)\in K_1\p L_1, \overline{\ell(x)} = \ell(y)\}$,
\item $\de_1^\alpha(x,y) = (\de_1^\alpha(x),y)$ for any $(x,y)\in
  K_1\p L_0$,
\item $\de_1^\alpha(x,y) = (x,\de_1^\alpha(y))$ for any $(x,y)\in
  K_0\p L_1$,
\item $\de_1^\alpha(x,y) = (\de_1^\alpha(x),\de_1^\alpha(y))$ for any $(x,y)\in
  K_1\p L_1$,
\item $\ell(x,y)=\ell(x)$ for any $(x,y)\in K_1\p L_0$,
\item $\ell(x,y)=\ell(y)$ for any $(x,y)\in K_0\p L_1$,
\item $\ell(x,y)=\tau$ for any $(x,y)\in K_1\p L_1$ with
  $\overline{\ell(x)} = \ell(y)$.
\end{itemize}
The $1$-cubes $(x,y)$ of $(K\p_\Sigma L)_1\cap (K_1\p L_1)$ are called
{\rm synchronizations} of $x$ and $y$.  \ed

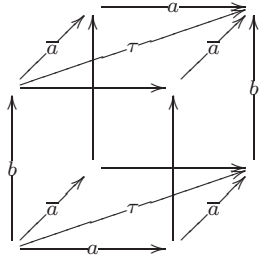
\begin{figure}
\[
\xymatrix{
&\ar@{->}[rr]|{a} &&\\
\ar@{->}[rrru]|-{\tau}\ar@{->}[ru]|-{\overline{a}}\ar@{->}[rr]&&\ar@{->}[ru]|-{\overline{a}}&\\
&\ar@{->}[uu] \ar@{->}[rr]&&\ar@{->}[uu]|-{b}\\
\ar@{->}[uu]|-{b}\ar@{->}[ru]|-{\overline{a}}\ar@{->}[rrru]|-{\tau}\ar@{->}[rr]|-{a}&&\ar@{->}[uu]\ar@{->}[ru]|-{\overline{a}}&
}\]
\caption{Representation of $\square[a,b]_{\leq 1}\p_\Sigma
  \square[\overline{a}]$, labelled over $\Sigma$}
\label{3sync00}
\end{figure}

The $1$-dimensional labelled precubical set $K\p_\Sigma L$ is the key
ingredient in the definition of the synchronized tensor product of
labelled precubical sets given in \cite{ccsprecub}, and recalled in
Section~\ref{compare}. Figure~\ref{3sync00} describes
$\square[a,b]_{\leq 1}\p_\Sigma \square[\overline{a}]$.

We want to prove in this section that for every $m\geq 0$ and $n\geq
0$, for every $a_1, \dots, a_{m+n} \in \Sigma$, the labelled
precubical set $\COSK^\Sigma(\square[a_1, \dots, a_m]_{\leq 1}
\p_\Sigma \square[a_{m+1}, \dots, a_{m+n}]_{\leq 1})$ can be
interpreted as a full labelled coskeleton in the category of labelled
transverse symmetric precubical sets.

\bp \label{pintuitive} Let $K$ be a precubical set. For any $p$-cube
$x$ of $\omega\mathcal{L}(K)$ with $p\geq 0$, there exists a $p$-cube
$y$ of $K \subset \omega\mathcal{L}(K)$ and a map $\mu \in
\widehat{\square}([p],[p])$ such that $x = \mu^*(y)$ where
$\mu^*:\mathcal{L}(K)_p \rightarrow \mathcal{L}(K)_p$ is the image of
$\mu$ by the presheaf $\mathcal{L}(K) \in \widehat{\square}^{op}\set$.
\ep

\bpf Let \[x\in \omega\mathcal{L}(K)_p \iso
\liminj_{\square[n]\rightarrow K} \widehat{\square}[n]_p \iso
\liminj_{\square[n]\rightarrow K} \widehat{\square}([p],[n]).\] Then
there exists an $n$-cube $z:\square[n]\rightarrow K$ and $\overline{x}
\in \widehat{\square}([p],[n])$ (the copy corresponding to $z$) such
that $z \circ \overline{x} = x$. By
Proposition~\ref{decomposition_distance}, $\overline{x} : [p]
\rightarrow [n]$ factors as a composite $[p] \stackrel{\mu}
\longrightarrow [p] \stackrel{\phi} \longrightarrow [n]$ with $\mu \in
\widehat{\square}$ and $\phi \in \square$. Then $\phi^*(z)$ is a
$p$-cube of $K$ and $\mu^*(\phi^*(z)) = x$. So $y = \phi^*(z)$ is a
solution. \epf

Note that the decomposition $x = \mu^*(y)$ is unique. But this fact
will not be used in the sequel. Indeed, let us consider another
decomposition $x = \mu'^*(y')$, $z':\square[n']\rightarrow K$, $z'
\circ \overline{x}' = x$ where $\overline{x}' : [p]
\stackrel{\mu'}\longrightarrow [p] \stackrel{\phi'}\longrightarrow
[n']$ belongs to the copy of $\widehat{\square}([p],[n'])$
corresponding to $z'$ and $y' = \phi'^*(z')$. Since $\overline{x} \in
\widehat{\square}([p],[n])$ and $\overline{x}' \in
\widehat{\square}([p],[n'])$ are equal in the colimit calculating
$\omega\mathcal{L}(K)_p$, the two sets $\widehat{\square}([p],[n])$
and $\widehat{\square}([p],[n'])$ are related in the colimit by a
zig-zag sequence of maps of $\square$ (this is the crucial point)
relating $\overline{x}$ and $\overline{x}'$. We can suppose that there
exists a map $h : [n] \rightarrow [n']$ such that $h\circ \overline{x}
= \overline{x}'$ and such that $z = z' \circ h$. Then the composite
$[p] \stackrel{\mu} \longrightarrow [p] \stackrel{\phi}
\longrightarrow [n] \stackrel{h} \longrightarrow [n']$ gives the
unique decomposition of $\overline{x}'$ as the composite of a map of
$\widehat{\square}([p], [p])$ followed by a map of $\square$ by
Proposition~\ref{decomposition_distance}. Thus, $\mu = \mu'$ and $h
\circ \phi = \phi'$. Therefore $y' = \phi'^*(z') = z' \circ h \circ
\phi = z \circ \phi = y$.

We will need the following combinatorial lemma twice in the sequel:

\bp \label{dec_non_twisted} Let $x : [p] \rightarrow [r]$ be a
strictly increasing set map.  Then there exists a unique decomposition
of $x$ as
\[[p] \stackrel{\mu} \longrightarrow [p'] \stackrel{\phi}
\longrightarrow [q] \stackrel{\psi} \longrightarrow [r]\] such that
$\phi$ is non-twisted, $\psi \in \square$ and $\mu =
(g_1,\dots,g_{p'})$ where the $g_i: [p] \rightarrow [1]$ are
non-constant and mutually distinct (i.e. $g_i = g_j$ implies $i =
j$). Moreover, $p\leq p'$, $\mu$ is strictly increasing, and it is
also adjacency-preserving if and only if $p = p'$.  \ep

\bpf Let $x = (x^{(1)}, \dots, x^{(r)})$ where the maps $x^{(i)} : [p]
\rightarrow [1]$ are the $r$ projection maps. The map $\psi$ is
necessarily the composite $\delta_{i_1}^{\alpha_1} \dots
\delta_{i_s}^{\alpha_s}$ where $\{i_1> \dots > i_s\} = \{i\in
\{1,\dots,r\} \mid x^{(i)} = 0 \hbox{ or } x^{(i)} = 1\}$ and where
$x^{(i_k)}(\epsilon_1,\dots,\epsilon_p) = \alpha_k$ for all
$(\epsilon_1,\dots,\epsilon_p) \in [p]$. Let $A\subset \{1,\dots,r\}$
be the subset of $i$ such that $x^{(i)}$ is a non-constant
map. Consider the equivalence relation on the set $A$ defined by
$i\sim j$ if and only if $x^{(i)}=x^{(j)}$. Let $p' = \card(A/\sim)$
where $\card(S)$ denotes the cardinality of the set $S$. The map $\mu
= (x^{(j_1)}, \dots, x^{(j_{p'})})$ is obtained by taking in each
equivalence class of $\sim$ the representative $x^{(j)}$ with the
smallest $j$ and by imposing $j_1 < \dots < j_{p'}$. The non-twisted
map $\phi$ is then defined so that the repetitions encode the
equivalence relation $\sim$.  Since $x$ is strictly increasing, the
set map $\mu$ is also strictly increasing. Therefore $p\leq p'$.
Since none of the set maps $x^{(j_k)}$ are constant, one has $\mu(0,
\dots, 0) = (0, \dots, 0)$ and $\mu(1, \dots, 1) = (1, \dots,
1)$. Thus, if $p < p'$, then $\mu$ cannot be adjacency-preserving. And
if $p = p'$, then $\mu$ is adjacency-preserving by
Proposition~\ref{example_preservation_distance}. This decomposition is
clearly unique.  \epf

\bth \label{inclusion_full_cosk} Let $\square[a_1,\dots,a_m]$ and
$\square[a_{m+1},\dots,a_{m+n}]$ be two labelled cubes with $m\geq 0$
and $n\geq 0$. Then there is an inclusion of presheaves
\begin{multline*} \mathcal{L} \lp\COSK^\Sigma\lp
  \square[a_1,\dots,a_m]_{\leq 1} \p_\Sigma
  \square[a_{m+1},\dots,a_{m+n}]_{\leq 1}\rp \rp \\\subset \cosk_1^{{\widehat{\square}},\Sigma}\lp
  {\widehat{\square}}[a_1,\dots,a_m]_{\leq 1} \p_\Sigma
  {\widehat{\square}}[a_{m+1},\dots,a_{m+n}]_{\leq 1}\rp .\end{multline*} Moreover,
when $\Sigma\backslash\{\tau\}$ is non-empty, there exist two labelled
cubes such that the above inclusion is strict.  \eth

\bpf Let $K$ be a labelled precubical set. Consider the composite set
map, natural with respect to $K$,
\begin{align*}
  &\lp \square^{op}\set\ddownarrow !\Sigma\rp \lp K,\COSK^\Sigma\lp
  \square[a_1,\dots,a_m]_{\leq 1} \p_\Sigma
  \square[a_{m+1},\dots,a_{m+n}]_{\leq 1}\rp \rp \\
  &\rightarrow \lp \square^{op}\set\ddownarrow!\Sigma\rp \lp
  K,\cosk_1^{\square,\Sigma}\lp \square[a_1,\dots,a_m]_{\leq 1}
  \p_\Sigma
  \square[a_{m+1},\dots,a_{m+n}]_{\leq 1}\rp \rp\\
  &\iso \lp \square_1^{op}\set\ddownarrow!\Sigma\rp \lp K_{\leq 1},
  \square[a_1,\dots,a_m]_{\leq 1} \p_\Sigma
  \square[a_{m+1},\dots,a_{m+n}]_{\leq 1}\rp \\
  &\iso \lp
  {\widehat{\square}}_1^{op}\set\ddownarrow\sh\mathcal{L}\lp
  !\Sigma\rp \rp \lp\mathcal{L}\lp K\rp _{\leq 1},
  {\widehat{\square}}[a_1,\dots,a_m]_{\leq 1} \p_\Sigma
  {\widehat{\square}}[a_{m+1},\dots,a_{m+n}]_{\leq 1}\rp \\
  &\iso \lp
  {\widehat{\square}}^{op}\set\ddownarrow\sh\mathcal{L}\lp
  !\Sigma\rp \rp \lp \mathcal{L}\lp K\rp
  ,\cosk_1^{{\widehat{\square}},\Sigma}\lp {\widehat{\square}}[a_1,\dots,a_m]_{\leq 1}
  \p_\Sigma {\widehat{\square}}[a_{m+1},\dots,a_{m+n}]_{\leq 1}\rp \rp,
\end{align*} the first and last isomorphisms by adjunction and the
second one by Proposition~\ref{2}. Take $K =
\COSK^\Sigma(\square[a_1,\dots,a_m]_{\leq 1} \p_\Sigma
\square[a_{m+1},\dots,a_{m+n}]_{\leq 1})$. The identity of $K$ yields
a map of labelled transverse symmetric precubical sets
\begin{multline*}f:\mathcal{L}\lp\COSK^\Sigma\lp
  \square[a_1,\dots,a_m]_{\leq 1} \p_\Sigma
  \square[a_{m+1},\dots,a_{m+n}]_{\leq 1}\rp \rp\\\longrightarrow
  \cosk_1^{{\widehat{\square}},\Sigma}\lp
  {\widehat{\square}}[a_1,\dots,a_m]_{\leq 1} \p_\Sigma
  {\widehat{\square}}[a_{m+1},\dots,a_{m+n}]_{\leq 1}\rp
  .\end{multline*} The case $K = \square[p]$ for $p\geq 0$ gives the
injection of sets
\begin{multline*}\lp\COSK^\Sigma\lp
  \square[a_1,\dots,a_m]_{\leq 1} \p_\Sigma
  \square[a_{m+1},\dots,a_{m+n}]_{\leq 1}\rp\rp_p \\\subset
  \lp\cosk_1^{{\widehat{\square}},\Sigma}\lp
  {\widehat{\square}}[a_1,\dots,a_m]_{\leq 1} \p_\Sigma
  {\widehat{\square}}[a_{m+1},\dots,a_{m+n}]_{\leq
    1}\rp\rp_p. \end{multline*} The set map $f_p$ is therefore
one-to-one for every $p\geq 0$ by Proposition~\ref{pintuitive} and
Proposition~\ref{dec_non_twisted}.  Suppose now that
$\Sigma\backslash\{\tau\}$ is non-empty. Let $a\in
\Sigma\backslash\{\tau\}$. The transverse symmetric precubical set
\[\cosk_1^{{\widehat{\square}},\Sigma}({\widehat{\square}}[a,a]_{\leq 1} \p_\Sigma
{\widehat{\square}}[\overline{a},\overline{a}]_{\leq 1})\] contains a
$2$-cube $x$ such that $x_0(0,0) = (0,0,0,0)$, $x_0(0,1) = (1,0,0,1)$,
$x_0(1,0) = (1,0,1,0)$ and $x_0(1,1) = (1,1,1,1)$ since all $1$-cubes
of $x$ are labelled by $\tau$. The set map $x_0:[2] \rightarrow [4]$
cannot be written as a composite $[2] \stackrel{\mu}\longrightarrow
[2] \stackrel{\phi} \longrightarrow [4]$ with $\mu\in
{\widehat{\square}}([2],[2])$ and $\phi:[2] \rightarrow [4]$
non-twisted (see Definition~\ref{def_directed}) since
$x_0=(x^{(1)},x^{(2)},x^{(3)},x^{(4)})$ where the set maps
$x^{(i)}:[2] \rightarrow [1]$ are four different set maps. So, by
Proposition~\ref{pintuitive}, one obtains
\[x\notin \mathcal{L}\lp\COSK^\Sigma(\square[a,a]_{\leq 1} \p_\Sigma
\square[\overline{a},\overline{a}]_{\leq 1})\rp.\] Therefore the
inclusion of presheaves
\[ \mathcal{L}\lp \COSK^\Sigma\lp \square[a,a]_{\leq 1}
\p_\Sigma \square[\overline{a},\overline{a}]_{\leq 1}\rp \rp \subset
\cosk_1^{{\widehat{\square}},\Sigma}\lp {\widehat{\square}}[a,a]_{\leq 1} \p_\Sigma
{\widehat{\square}}[\overline{a},\overline{a}]_{\leq 1}\rp \] is strict.  \epf

\subsection*{Functorial construction of the parallel composition}

\begin{figure}
\[
\xymatrix{
&\ar@{->}[rr]|{a} &&\\
\ar@{->}[rrru]|-{(1,3)}\ar@{->}[ru]|-{\overline{a}}\ar@{->}[rr]&&\ar@{->}[ru]|-{\overline{a}}&\\
&\ar@{->}[uu] \ar@{->}[rr]&&\ar@{->}[uu]|-{b}\\
\ar@{->}[uu]|-{b}\ar@{->}[ru]|-{\overline{a}}\ar@{->}[rrru]|-{(1,3)}\ar@{->}[rr]|-{a}&&\ar@{->}[uu]\ar@{->}[ru]|-{\overline{a}}&
}\]
\caption{Representation of $\square[a,b]_{\leq 1}\overline{\p}_\Sigma
  \square[\overline{a}]$, labelled over $\overline{\Sigma}= \Sigma
  \sqcup (\N^* \p \N^*)$}
\label{3sync01}
\end{figure}
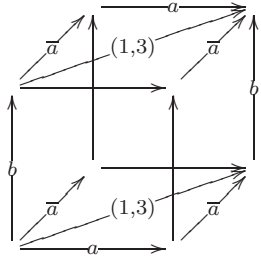

Theorem~\ref{inclusion_full_cosk} does not mean that the labelled
coskeleton functor of the category of labelled transverse symmetric
precubical sets is badly behaved. The coskeleton functor of
$\widehat{\square}^{op}\set$ does the job it is designed for: filling
all compatibly labelled shells. To avoid this problem, we have to keep
the memory of what is synchronized by $\tau$, as depicted in
Figure~\ref{3sync01}. By labelling the $1$-cube $x(0,*)$ by $(2,3)$
instead of $\tau$, the $1$-cube $x(1,*)$ by $(2,4)$ instead of $\tau$,
the $1$-cube $x(*,0)$ by $(1,3)$ instead of $\tau$ and the $1$-cube
$x(*,1)$ by $(1,4)$ instead of $\tau$, it becomes impossible to fill
the new shell since the opposite faces are not labelled anymore in the
same way. Hence the definition of the new labelling: 

\begin{nota} \label{newfibered} Let $a_1,\dots, a_{m+n} \in \Sigma$
  with $m\geq 0$ and $n\geq 0$. Let $\overline{\Sigma} := \Sigma
  \sqcup (\N^* \p \N^*)$ where $\N^*$ is the set of strictly positive
  integers. Let us define the $1$-dimensional labelled (transverse
  symmetric) precubical set $\widehat{\square}[a_1,\dots,a_m]_{\leq 1}
  \overline{\p}_\Sigma \widehat{\square}[a_{m+1},\dots,a_{m+n}]_{\leq
    1}$ as follows (the boxed part is the only new part):
\begin{itemize}
\item The underlying $1$-dimensional precubical set is the one of  
 \[\widehat{\square}[a_1,\dots,a_m]_{\leq 1} \p_\Sigma
  \widehat{\square}[a_{m+1},\dots,a_{m+n}]_{\leq 1}\]
\item The labelling map is defined by:  
\begin{itemize}
\item $\ell(x,y) = \ell(x)$ for any $(x,y)\in
  \widehat{\square}[a_1,\dots,a_m]_1\p
  \widehat{\square}[a_{m+1},\dots,a_{m+n}]_0$,
\item $\ell(x,y) = \ell(y)$ for any $(x,y)\in
  \widehat{\square}[a_1,\dots,a_m]_0\p
  \widehat{\square}[a_{m+1},\dots,a_{m+n}]_1$,
\item $\boxed{\ell(x,y) = (r,s)\in \N^*\p \N^*}$~\footnote{instead of
    $\ell(x,y) = \tau$} for any $(x,y)\in
  \widehat{\square}[a_1,\dots,a_m]_1\p
  \widehat{\square}[a_{m+1},\dots,a_{m+n}]_1$ with $\overline{\ell(x)}
  = \ell(y)$ where $1\leq r \leq m$ and $m+1 \leq s \leq m+n$ are the
  unique integers such that $(x_0(\alpha),y_0(\alpha)) =
  \delta_{s}^\alpha\delta_r^\alpha(X)$ for some $X\in [m+n-2]$ and for
  $\alpha = 0,1$.
\end{itemize}
\end{itemize}
\end{nota}

\begin{lem} \label{auplus2} Let $c : \widehat{\square}[1] \rightarrow
  \widehat{\square}[a_1,\dots,a_m]_{\leq 1} \overline{\p}_\Sigma
  \widehat{\square}[a_{m+1},\dots,a_{m+n}]_{\leq 1}$ be a $1$-cube of
  $\widehat{\square}[a_1,\dots,a_m]_{\leq 1} \overline{\p}_\Sigma
  \widehat{\square}[a_{m+1},\dots,a_{m+n}]_{\leq 1}$. Then the set map
  $c_0:[1] \rightarrow [m+n]$ satisfies $c_0(0) < c_0(1)$ and there
  are two mutually exclusive possibilities:
\begin{itemize}
\item $d(c_0(0),c_0(1)) = 1$ and $c_0(\alpha) = \delta_r^\alpha(X)$
  for some $X\in [m+n-1]$ with $1\leq r \leq m+n$ and for $\alpha =
  0,1$. In this case, $\ell(c) = a_r \in \Sigma$.
\item $d(c_0(0),c_0(1)) = 2$ and $c_0(\alpha) =
  \delta_{s}^\alpha\delta_r^\alpha(X)$ for some $X\in [m+n-2]$ with
  $1\leq r \leq m$ and $m+1 \leq s \leq m+n$ and for $\alpha =
  0,1$. In this case, $\ell(c) = (r,s) \in \N^* \p \N^*$.
\end{itemize}
\end{lem}

\bpf Obvious. \epf 

Note that Lemma~\ref{auplus2} holds for
$\widehat{\square}[a_1,\dots,a_m]_{\leq 1} \p_\Sigma
\widehat{\square}[a_{m+1},\dots,a_{m+n}]_{\leq 1}$ as well by
replacing in the last sentence $\ell(c) = (r,s) \in \N^* \p \N^*$ by
$\ell(c) = \tau$.

We are now ready to give the categorical interpretation of the
labelled directed coskeleton construction when applied to the fibered
product of two $1$-dimensional labelled precubical sets.

\bth \label{cavavraiment} Let $\square[a_1,\dots,a_m]$ and
$\square[a_{m+1},\dots,a_{m+n}]$ be two labelled cubes with $m\geq 0$
and $n\geq 0$. Then one has the isomorphism of labelled transverse
symmetric precubical sets
\begin{multline*}
  \mathcal{L}\lp\COSK^\Sigma\lp \square[a_1,\dots,a_m]_{\leq 1}
  \p_\Sigma \square[a_{m+1},\dots,a_{m+n}]_{\leq 1}\rp \rp\\\iso
  \cosk_1^{\widehat{\square},\overline{\Sigma}}\lp
  \widehat{\square}[a_1,\dots,a_m]_{\leq 1} \overline{\p}_\Sigma
  \widehat{\square}[a_{m+1},\dots,a_{m+n}]_{\leq 1}\rp
\end{multline*}
where the right-hand labelled transverse symmetric precubical set over
$\overline{\Sigma}$ is viewed as labelled over $\Sigma$ by composing
its labelling map with the morphism of transverse symmetric precubical
sets $\sh\mathcal{L}(!\overline{\Sigma}) \rightarrow
\sh\mathcal{L}(!\Sigma)$, where the set map $\overline{\Sigma}
\rightarrow \Sigma$ is defined as the identity on $\Sigma$ and by the
mapping $(p,q) \mapsto \tau$ on the complement. \eth

Note that with $m=0$ or $n=0$, we have the isomorphism of
Theorem~\ref{but_math}.

\bpf \boxed{\hbox{Injectivity}}.  There is an inclusion of presheaves
\begin{multline*} \cosk_1^{\widehat{\square},\overline{\Sigma}}\lp
  \widehat{\square}[a_1,\dots,a_m]_{\leq 1} \overline{\p}_\Sigma
  \widehat{\square}[a_{m+1},\dots,a_{m+n}]_{\leq 1}\rp \\\subset
  \cosk_1^{\widehat{\square},\Sigma}\lp
  \widehat{\square}[a_1,\dots,a_m]_{\leq 1} \p_\Sigma
  \widehat{\square}[a_{m+1},\dots,a_{m+n}]_{\leq 1}\rp \end{multline*}
since less shells are filled in the left-hand term than in the
right-hand term because of the labelling over
$\overline{\Sigma}$. Moreover, one has the equalities
\begin{multline*} \lp\mathcal{L} \lp\COSK^\Sigma\lp
  \square[a_1,\dots,a_m]_{\leq 1} \p_\Sigma
  \square[a_{m+1},\dots,a_{m+n}]_{\leq 1}\rp \rp\rp_{\leq 1} \\=
  \lp\cosk_1^{\widehat{\square},\overline{\Sigma}}\lp
  \widehat{\square}[a_1,\dots,a_m]_{\leq 1} \overline{\p}_\Sigma
  \widehat{\square}[a_{m+1},\dots,a_{m+n}]_{\leq
    1}\rp \rp_{\leq 1}\\
  =\lp\cosk_1^{\widehat{\square},\Sigma}\lp
  \widehat{\square}[a_1,\dots,a_m]_{\leq 1} \p_\Sigma
  \widehat{\square}[a_{m+1},\dots,a_{m+n}]_{\leq 1}\rp \rp_{\leq
    1} \end{multline*} by Proposition~\ref{restrict1_1} and
Proposition~\ref{i1_relatif}. Let \[x\in \lp\mathcal{L}\lp
\COSK^\Sigma(\square[a_1,\dots,a_m]_{\leq 1} \p_\Sigma
\square[a_{m+1},\dots,a_{m+n}]_{\leq 1})\rp\rp_p\] with $p\geq
2$. Then $x=\mu^*(y)$ where $\mu:[p] \rightarrow [p]$ is an
adjacency-preserving map and where $y$ is a $p$-cube of the labelled
precubical set 
\[\COSK^\Sigma(\square[a_1,\dots,a_m]_{\leq 1} \p_\Sigma
\square[a_{m+1},\dots,a_{m+n}]_{\leq 1})\] by
Proposition~\ref{pintuitive}. The map of $1$-dimensional precubical
sets
\[y_{\leq 1} : \square[p]_{\leq 1} \longrightarrow
(\COSK^\Sigma(\square[a_1,\dots,a_m]_{\leq 1} \p_\Sigma
\square[a_{m+1},\dots,a_{m+n}]_{\leq 1}))_{\leq 1}\] induces by
Proposition~\ref{2} a map of transverse symmetric $1$-dimensional
precubical sets \[y_{\leq 1}:\widehat{\square}[p]_{\leq 1} \rightarrow
\widehat{\square}[a_1,\dots,a_m]_{\leq 1} \p_\Sigma
\widehat{\square}[a_{m+1},\dots,a_{m+n}]_{\leq 1}.\] The latter
induces a unique map $\overline{y}: \widehat{\square}[p]_{\leq 1}
\rightarrow \widehat{\square}[a_1,\dots,a_m]_{\leq 1}
\overline{\p}_\Sigma \widehat{\square}[a_{m+1},\dots,a_{m+n}]_{\leq
  1}$ of $1$-dimensional precubical sets which is this time labelled
over $\overline{\Sigma}$ since the underlying precubical sets of
$\widehat{\square}[a_1,\dots,a_m]_{\leq 1} \p_\Sigma
\widehat{\square}[a_{m+1},\dots,a_{m+n}]_{\leq 1}$ and
$\widehat{\square}[a_1,\dots,a_m]_{\leq 1} \overline{\p}_\Sigma
\widehat{\square}[a_{m+1},\dots,a_{m+n}]_{\leq 1}$ are equal. The map
$\overline{y}$ induces by adjunction a unique $p$-dimensional
transverse symmetric cube of
\[\cosk_1^{\widehat{\square},\overline{\Sigma}}(\widehat{\square}[a_1,\dots,a_m]_{\leq
  1} \overline{\p}_\Sigma
\widehat{\square}[a_{m+1},\dots,a_{m+n}]_{\leq 1}).\] Thus the
inclusion \begin{multline*}\mathcal{L} \lp\COSK^\Sigma\lp
  \square[a_1,\dots,a_m]_{\leq 1} \p_\Sigma
  \square[a_{m+1},\dots,a_{m+n}]_{\leq 1}\rp \rp \\\subset
  \cosk_1^{\widehat{\square},\Sigma}\lp
  \widehat{\square}[a_1,\dots,a_m]_{\leq 1} \p_\Sigma
  \widehat{\square}[a_{m+1},\dots,a_{m+n}]_{\leq 1}\rp \end{multline*}
factors uniquely as a composite of inclusions
\begin{multline*}\mathcal{L} \lp\COSK^\Sigma\lp \square[a_1,\dots,a_m]_{\leq
    1} \p_\Sigma \square[a_{m+1},\dots,a_{m+n}]_{\leq 1}\rp \rp
  \\\subset \cosk_1^{\widehat{\square},\overline{\Sigma}}\lp
  \widehat{\square}[a_1,\dots,a_m]_{\leq 1} \overline{\p}_\Sigma
  \widehat{\square}[a_{m+1},\dots,a_{m+n}]_{\leq 1}\rp \\\subset
  \cosk_1^{\widehat{\square},\Sigma}\lp
  \widehat{\square}[a_1,\dots,a_m]_{\leq 1} \p_\Sigma
  \widehat{\square}[a_{m+1},\dots,a_{m+n}]_{\leq 1}\rp
  .\end{multline*} Let us call $f$ again the inclusion
\begin{multline*}\mathcal{L} \lp\COSK^\Sigma\lp \square[a_1,\dots,a_m]_{\leq 1} \p_\Sigma
  \square[a_{m+1},\dots,a_{m+n}]_{\leq 1}\rp \rp\\\subset
  \cosk_1^{\widehat{\square},\overline{\Sigma}}\lp
  \widehat{\square}[a_1,\dots,a_m]_{\leq 1} \overline{\p}_\Sigma
  \widehat{\square}[a_{m+1},\dots,a_{m+n}]_{\leq 1}\rp
  .\end{multline*} It then remains to prove that for every $p\geq 2$,
the set map $f_p$ is onto. 

\boxed{\hbox{Surjectivity}}. Let $x:\widehat{\square}[p] \rightarrow
\cosk_1^{\widehat{\square},\overline{\Sigma}}(\widehat{\square}[a_1,\dots,a_m]_{\leq
  1} \overline{\p}_\Sigma
\widehat{\square}[a_{m+1},\dots,a_{m+n}]_{\leq 1})$ be a
$p$-dimensional transverse symmetric cube of
$\cosk_1^{\widehat{\square},\overline{\Sigma}}(\widehat{\square}[a_1,\dots,a_m]_{\leq
  1} \overline{\p}_\Sigma
\widehat{\square}[a_{m+1},\dots,a_{m+n}]_{\leq 1})$ with $p\geq 2$.
Let $x_0 = (x^{(1)},\dots,x^{(m+n)})$ where the $x^{(i)}:[p]
\rightarrow [1]$ are the $m+n$ projections.  Let us apply the
decomposition of Proposition~\ref{dec_non_twisted}. Let $\mu =
(x^{(j_1)}, \dots, x^{(j_{p'})})$. If one had $p < p'$, then there
would exist a $1$-cube $c:[1] \rightarrow [p]$ such that
$d(\mu(c(0)),\mu(c(1))) > 1$. By Lemma~\ref{auplus2}, one would
have \[1 < d(\mu(c(0)),\mu(c(1))) \leq
d(\psi\phi\mu(c(0)),\psi\phi\mu(c(1))) \leq 2,\] and therefore
$d(\mu(c(0)),\mu(c(1))) = 2$.  Thus, one would have $\mu c(\alpha) =
\delta_v^\alpha \delta_u^\alpha(X)$ for some $u < v$, $X \in [p'-2]$
and for $\alpha = 0,1$. By Lemma~\ref{auplus2}, one obtains
$\psi\phi\mu c(\alpha) = \delta_{j_v}^\alpha \delta_{j_u}^\alpha(Z)$
for some cube $Z \in [m+n-2]$ and for $\alpha = 0,1$~\footnote{So far,
  the particular labelling of $\overline{\p}_\Sigma$ has not been used
  in the surjectivity part of the proof. In the counterexample of
  Theorem~\ref{inclusion_full_cosk}, one has $p = 2$ and $p' = 4$. So
  we cannot yet conclude that $p = p'$.}, and finally
$\ell(\psi\phi\mu c) = (j_u,j_v)$.

\underline{Use of the particular labelling of
  $\overline{\p}_\Sigma$}. The crucial point is that the labelling of
$\overline{\p}_\Sigma$ implies $x^{(j_u)} = x^{(j_v)}$, which
contradicts the definition of $\mu$. By Proposition~\ref{maximal}, the
commutative word $W = \ell(xc_1) \dots \ell(xc_p)$ of the free
commutative monoid without unit generated by $\overline{\Sigma}$ does
not depend on the maximal path $(c_1,\dots,c_p)$ of
$\widehat{\square}[p]$. And one of the labels is necessarily
$(j_u,j_v)$. If $(\epsilon_1, \dots, \epsilon_p) = (0,\dots, 0)$, then
$x^{(j_u)}(\epsilon_1, \dots, \epsilon_p) = x^{(j_v)}(\epsilon_1,
\dots, \epsilon_p) =0$. Let us suppose now that $(\epsilon_1, \dots,
\epsilon_p) \neq (0,\dots, 0)$. By Proposition~\ref{maximal}, for
every maximal path $(c_1,\dots,c_r)$ of the $r$-subcube from
$(0,\dots, 0)$ to $(\epsilon_1, \dots, \epsilon_p)$, the commutative
word $W' = \ell(xc_1) \dots \ell(xc_r)$ is a subword of $W$ which does
not depend on $(c_1,\dots,c_r)$. If $(j_u,j_v)$ belongs to $W'$, then
$x^{(j_u)}(\epsilon_1, \dots, \epsilon_p) = x^{(j_v)}(\epsilon_1,
\dots, \epsilon_p) = 1$. If $(j_u,j_v)$ does not belong to $W'$, then
it belongs to the complement of $W'$ in $W$. So $x^{(j_u)}(\epsilon_1,
\dots, \epsilon_p) = x^{(j_v)}(\epsilon_1, \dots, \epsilon_p) =
0$. Hence, $x^{(j_u)} = x^{(j_v)}$, which is the desired
contradiction.

\underline{End of the proof}. Hence, one obtains the equality $p = p'$
thanks to the particular labelling of $\overline{\p}_\Sigma$. The map
$\mu$ is therefore adjacency-preserving by
Proposition~\ref{example_preservation_distance}. Note that $x_0$ has
no reason to be adjacency-preserving.  By definition of the labelled
directed coskeleton, there exists a $p$-cube \[y:\square[p]
\longrightarrow \COSK^\Sigma\lp \square[a_1,\dots,a_m]_{\leq 1}
\p_\Sigma \square[a_{m+1},\dots,a_{m+n}]_{\leq 1}\rp \] such that $y_0
= \psi\phi$. Then $\mathcal{L}(y)$ is a $p$-cube of the labelled
transverse symmetric precubical set $\mathcal{L}\lp
\COSK^\Sigma(\square[a_1,\dots,a_m]_{\leq 1} \p_\Sigma
\square[a_{m+1},\dots,a_{m+n}]_{\leq 1})\rp$ such that
$\mathcal{L}(y)_0 = \psi\phi$. Then $\mu^{*}(\mathcal{L}(y))$ is a
$p$-cube of the labelled transverse symmetric precubical
set \[\mathcal{L}\lp \COSK^\Sigma(\square[a_1,\dots,a_m]_{\leq 1}
\p_\Sigma \square[a_{m+1},\dots,a_{m+n}]_{\leq 1})\rp\] such that
$(\mu^{*}(\mathcal{L}(y)))_0 = \psi\phi\mu$. By construction of $f$,
the $p$-cube $f(\mu^{*}(\mathcal{L}(y)))$ of the labelled transverse
symmetric precubical set
$\cosk_1^{\widehat{\square},\overline{\Sigma}}(\widehat{\square}[a_1,\dots,a_m]_{\leq
  1} \overline{\p}_\Sigma
\widehat{\square}[a_{m+1},\dots,a_{m+n}]_{\leq 1})$ satisfies
$(f(\mu^{*}(\mathcal{L}(y))))_0 = \psi\phi\mu = x_0$. Since there is
at most one $1$-cube between two vertices of
$\widehat{\square}[a_1,\dots,a_m]_{\leq 1} \overline{\p}_\Sigma
\widehat{\square}[a_{m+1},\dots,a_{m+n}]_{\leq 1}$, this implies
$(f(\mu^{*}(\mathcal{L}(y))))_{\leq 1} = x_{\leq 1}$ and therefore
$f(\mu^{*}(\mathcal{L}(y)))=x$ by adjunction. So $f$ is an isomorphism
of labelled transverse symmetric precubical sets.  \epf

Theorem~\ref{cavavraiment} is of course false for any other category
of cubes than $\widehat{\square}$. Indeed, the particular case $n=0$
and $a_1=\dots=a_m=\tau$ gives back the inclusion of
presheaves \[\mathcal{A}[m] \iso \mathcal{L}_\mathcal{A}\lp
\COSK^\Sigma\lp \square[m]_{\leq 1}\rp\rp \subset
\cosk_1^{\mathcal{A}}\lp \mathcal{A}[m]_{\leq 1}\rp
\]
which is an equality if and only if the category of cubes
$\mathcal{A}$ is shell-complete, so if and only if $\mathcal{A} =
\widehat{\square}$ by Theorem~\ref{identification_smallest}. The
crucial point in the proof of Theorem~\ref{cavavraiment} is that the
map $\mu:[p] \rightarrow [p]$ must belong to $\mathcal{A}$. Therefore,
it is really needed to work with the whole category
$\widehat{\square}$ of all adjacency-preserving maps.

\section{Parallel composition (global case)}
\label{compare}

We can now relate the synchronized tensor product of labelled
precubical sets with the synchronized tensor product of labelled
transverse symmetric precubical sets. First of all, let us give the
definition of these two synchronized tensor products.

\subsection*{Definition}

\bd \cite{ccsprecub} Let $K$ and $L$ be two labelled precubical
sets. The {\rm tensor product with synchronization} (or {\rm
  synchronized tensor product}) of $K$ and $L$ is
\[K \ot_\Sigma L := \liminj_{\square[m]\rightarrow K}
\liminj_{\square[n]\rightarrow L} \COSK^\Sigma(\square[m]_{\leq 1}
\p_\Sigma \square[n]_{\leq 1}).\] \ed

\bd Let $K$ and $L$ be two labelled transverse symmetric precubical
sets. The {\rm tensor product with synchronization} (or {\rm
  synchronized tensor product}) of $K$ and $L$ is
\[K \ot_\Sigma L := \liminj_{\widehat{\square}[m]\rightarrow K}
\liminj_{\widehat{\square}[n]\rightarrow L}
\cosk_1^{\widehat{\square},\overline{\Sigma}}(\widehat{\square}[m]_{\leq 1}
\overline{\p}_\Sigma \widehat{\square}[n]_{\leq 1}).\] \ed 

\subsection*{The two constructions coincide}

For the sequel, the category of small categories is denoted by
$\cat$. Let $H: I \longrightarrow \cat$ be a functor from a small
category $I$ to $\cat$. The \textit{Grothendieck construction} $I
\intop H$ is the category defined as follows~\cite{MR510404}: the
objects are the pairs $(i,a)$ where $i$ is an object of $I$ and $a$ is
an object of $H(i)$; a morphism $(i,a) \rightarrow (j,b)$ consists in
a map $\phi: i \rightarrow j$ and in a map $h : H(\phi)(a) \rightarrow
b$.

\begin{lem} \label{colim_id} Let $\mathcal{A}$ be a category of
  cubes. Let $I$ be a small category, and $i\mapsto K^i$ be a functor
  from $I$ to the category of labelled $\mathcal{A}$-sets.  Let $K =
  \liminj_i K^i$. Let $H: I \rightarrow \cat$ be the functor defined
  by $H(i) = \mathcal{A} \ddownarrow K^i$.  Then the functor $\iota :
  I\intop H \rightarrow \mathcal{A} \ddownarrow K$ defined by
  $\iota(i,\mathcal{A}[m] \rightarrow K^i) = (\mathcal{A}[m]
  \rightarrow K)$ is final in the sense of \cite{MR1712872}; that is
  to say the comma category $k \ddownarrow \iota$ is nonempty and
  connected for all objects $k$ of $\mathcal{A} \ddownarrow
  K$.  \end{lem}

\bpf The proof is similar to the proof of
\cite[Lemma~A.1]{ccsprecub}. \epf

\bp \label{preserve_colim} Let $\mathcal{A}$ be a category of
cubes. Let $F:\mathcal{A}\p \mathcal{A} \rightarrow \C$ be a functor
where $\C$ is a cocomplete category. Let
$\widehat{F}:(\mathcal{A}^{op}\set\ddownarrow
\sh_\mathcal{A}\mathcal{L}_\mathcal{A}(!\Sigma))\p
(\mathcal{A}^{op}\set\ddownarrow
\sh_\mathcal{A}\mathcal{L}_\mathcal{A}(!\Sigma)) \rightarrow \C$ be
the functor defined by
\[\widehat{F}(K,L):=\liminj_{\mathcal{A}[m]\rightarrow K}
\liminj_{\mathcal{A}[n]\rightarrow L} F([m],[n]).\] Then for any
labelled $\mathcal{A}$-set $L$, the two functors
\[\widehat{F}(L,-):\mathcal{A}^{op}\set \ddownarrow
\sh_\mathcal{A}\mathcal{L}_\mathcal{A}(!\Sigma)\longrightarrow \C\] and
\[\widehat{F}(-,L):\mathcal{A}^{op}\set \ddownarrow
\sh_\mathcal{A}\mathcal{L}_\mathcal{A}(!\Sigma) \longrightarrow \C\]
are colimit-preserving. \ep

\bpf The proof is similar to the proof of
\cite[Proposition~A.2]{ccsprecub}. Let $K = \liminj_i K^i$ be a colimit
of labelled $\mathcal{A}$-sets. By definition, one has the isomorphism
\[ \liminj_i \widehat{F}(K^i,L) \iso \liminj_i
\liminj_{\mathcal{A}[m]\rightarrow K^i}
\liminj_{\mathcal{A}[n]\rightarrow L} F([m],[n]).\] Consider the
functor $H: I \longrightarrow \cat$ defined by $H(i) = \mathcal{A}
\ddownarrow K^i$.  Consider the functors $F_i: H(i) \longrightarrow
\C$ defined by
\[
F_i(\mathcal{A}[m] \rightarrow K^i) =
\liminj_{\mathcal{A}[n]\rightarrow L} F([m],[n]).
\]
Consider the functor $\overline{F} : I\intop H \longrightarrow \C$
defined by
\[\overline{F}(i,\mathcal{A}[m] \rightarrow K^i) =
\liminj_{\mathcal{A}[n]\rightarrow L} F([m],[n]).\] Then the composite
$H(i) \subset I\intop H \rightarrow \C$ is exactly $F_i$. Therefore
one has the isomorphism
\[\liminj_i \liminj_{\mathcal{A}[m]\rightarrow
  K^i} \liminj_{\mathcal{A}[n]\rightarrow L} F([m],[n]) \iso
\liminj_{(i,\mathcal{A}[m]\rightarrow K^i)}
\liminj_{\mathcal{A}[n]\rightarrow L} F([m],[n])\] by
\cite[Proposition 40.2]{monographie_hocolim}. The functor $\iota:
I\intop H \rightarrow \mathcal{A} \ddownarrow K$ defined by
$\iota(i,\mathcal{A}[m] \rightarrow K^i) = (\mathcal{A}[m] \rightarrow
K)$ is final in the sense of~\cite{MR1712872} by Lemma~\ref{colim_id}.
Therefore by~\cite[p.\ 213, Theorem 1]{MR1712872} or~\cite[Theorem
14.2.5]{ref_model2}, one has the isomorphism
\[\liminj_{(i,\mathcal{A}[m]\rightarrow K^i)}
\liminj_{\mathcal{A}[n]\rightarrow L} F([m],[n]) \iso
\liminj_{\mathcal{A}[m]\rightarrow K}
\liminj_{\mathcal{A}[n]\rightarrow L} F([m],[n]) =:
\widehat{F}(K,L).\] \epf

\begin{cor} \label{la} Let $I$ be a small category. Let $i\mapsto K^i$
  be a functor from $I$ to the category of labelled transverse
  symmetric precubical sets, and let $L$ be a labelled transverse
  symmetric precubical set.  Then one has the natural isomorphism
  \[(\liminj_i K^i) \ot_\Sigma L \iso \liminj_i (K^i \ot_\Sigma
  L).\] \end{cor}

\bth \label{cavavraiment2} Let $K$ and $L$ be two labelled precubical
sets. Then there is the natural isomorphism of labelled transverse
symmetric precubical sets \[\mathcal{L}(K\ot_\Sigma L) \iso
\mathcal{L}(K) \ot_\Sigma \mathcal{L}(L).\] \eth

\bpf One has 

\begin{align*}
  \mathcal{L}(K\ot_\Sigma L) & \iso \mathcal{L}\lp
  \liminj_{\square[m]\rightarrow K} \liminj_{\square[n]\rightarrow L}
  \COSK^\Sigma(\square[m]_{\leq 1}
  \p_\Sigma \square[n]_{\leq 1})\rp & \\
  &\iso \liminj_{\square[m]\rightarrow K}
  \liminj_{\square[n]\rightarrow L}
  \mathcal{L}\lp\COSK^\Sigma(\square[m]_{\leq 1} \p_\Sigma
  \square[n]_{\leq 1})\rp & \hbox{ since $\mathcal{L}$ is a
    left adjoint}\\
  & \iso \liminj_{\square[m]\rightarrow K}
  \liminj_{\square[n]\rightarrow L}
  \cosk_1^{\widehat{\square},\overline{\Sigma}}(\widehat{\square}[m]_{\leq
    1} \overline{\p}_\Sigma \widehat{\square}[n]_{\leq 1}) & \hbox{ by
    Theorem~\ref{cavavraiment}}\\
  & \iso \liminj_{\square[m]\rightarrow K}
  \liminj_{\square[n]\rightarrow L} \widehat{\square}[m] \ot_\Sigma
  \widehat{\square}[n] & \hbox{ by definition of $\ot_\Sigma$}\\
  & \iso \liminj_{\square[m]\rightarrow K}
  \liminj_{\square[n]\rightarrow L} \mathcal{L}(\square[m]) \ot_\Sigma
  \mathcal{L}(\square[n])& \hbox{ by Proposition~\ref{carre}}\\
  & \iso \lp \liminj_{\square[m]\rightarrow K} \mathcal{L}(\square[m])
  \rp \ot_\Sigma \lp \liminj_{\square[n]\rightarrow L}
  \mathcal{L}(\square[n])\rp &\hbox{ by Corollary~\ref{la}}\\
  &\iso \mathcal{L}(K) \ot_\Sigma \mathcal{L}(L) &\hbox{ since
    $\mathcal{L}$ is a left adjoint.}
\end{align*}

\epf 

\subsection*{Associativity}

As in \cite{ccsprecub}, it is also possible to prove that the
synchronized tensor product of labelled transverse symmetric
precubical sets is associative.

\bth Let $K$, $L$ and $M$ be three labelled transverse symmetric
precubical sets. Then there is a canonical isomorphism of labelled
transverse symmetric precubical sets \[(K\ot_\Sigma L)\ot_\Sigma M
\iso K\ot_\Sigma (L\ot_\Sigma M).\] \eth

\bpf One has 

\begin{align*} 
  & K \ot_\Sigma (L \ot_\Sigma M) &\\
  & \iso \lp \liminj_{\widehat{\square}[p]\rightarrow K}
  \widehat{\square}[p] \rp \ot_\Sigma \lp
  \liminj_{\widehat{\square}[q]\rightarrow L}
  \liminj_{\widehat{\square}[r]\rightarrow
    M}  \widehat{\square}[q] \ot_\Sigma \widehat{\square}[r]\rp & \\
  & \iso \liminj_{\widehat{\square}[p]\rightarrow
    K}\liminj_{\widehat{\square}[q]\rightarrow L}
  \liminj_{\widehat{\square}[r]\rightarrow M} \widehat{\square}[p]
  \ot_\Sigma(\widehat{\square}[q] \ot_\Sigma \widehat{\square}[r]) &
  \hbox{
    by Corollary~\ref{la}}\\
  & \iso \liminj_{\widehat{\square}[p]\rightarrow
    K}\liminj_{\widehat{\square}[q]\rightarrow L}
  \liminj_{\widehat{\square}[r]\rightarrow M} \mathcal{L}\lp
  \square[p] \ot_\Sigma(\square[q] \ot_\Sigma \square[r])\rp & \hbox{
    by Theorem~\ref{cavavraiment2}}\\
  & \iso \liminj_{\widehat{\square}[p]\rightarrow
    K}\liminj_{\widehat{\square}[q]\rightarrow L}
  \liminj_{\widehat{\square}[r]\rightarrow M} \mathcal{L}\lp
  (\square[p] \ot_\Sigma\square[q]) \ot_\Sigma \square[r]\rp & \hbox{
    by \cite[Proposition~A.3]{ccsprecub}}\\
  & \iso \liminj_{\widehat{\square}[p]\rightarrow
    K}\liminj_{\widehat{\square}[q]\rightarrow L}
  \liminj_{\widehat{\square}[r]\rightarrow M} (\widehat{\square}[p]
  \ot_\Sigma\widehat{\square}[q]) \ot_\Sigma \widehat{\square}[r] &
  \hbox{
    by Theorem~\ref{cavavraiment2}}\\
  &\iso (K \ot_\Sigma L) \ot_\Sigma M& \hbox{ by Corollary~\ref{la}.}
\end{align*} 
\epf

As already pointed out in \cite{ccsprecub}, it is false in general
that the two labelled precubical sets $K\ot_\Sigma L$ and $L\ot_\Sigma
K$ are isomorphic as labelled precubical sets. Indeed, let us suppose
that $\Sigma\backslash\{\tau\}$ contains an element $a$. Then
$\square[a]\ot_\Sigma \square[\tau] \iso \square[a,\tau]$ and
$\square[\tau]\ot_\Sigma \square[a] \iso \square[\tau,a]$. Because of
the lack of symmetry operators, the two labelled $2$-cubes
$\square[a,\tau]$ and $\square[\tau,a]$ cannot be isomorphic as
labelled precubical sets. However, the two underlying precubical sets
are of course isomorphic, as already pointed out in
\cite{ccsprecub}. In the category of transverse symmetric precubical
sets, the situation is much better. Indeed, one has the isomorphisms
of labelled transverse symmetric precubical sets
\[\widehat{\square}[a,\tau] \iso \mathcal{L}(\square[a,\tau]) \iso 
\widehat{\square}[\tau,a] \iso \mathcal{L}(\square[\tau,a]).\]

\bp \label{commutativite} Let $K$ and $L$ be two labelled transverse
symmetric precubical sets. Then there is a natural isomorphism of
labelled transverse symmetric precubical sets $K\ot_\Sigma L \iso
L\ot_\Sigma K$. \ep

\bpf[Sketch of proof] It suffices to use Corollary~\ref{la} together
with the isomorphism \[\widehat{\square}[a_1,\dots,a_m]\ot_\Sigma
\widehat{\square}[b_1,\dots,b_n] \iso
\widehat{\square}[b_1,\dots,b_n]\ot_\Sigma
\widehat{\square}[a_1,\dots,a_m]\] (built using the symmetry
operators) for all labelled full cubes
$\widehat{\square}[a_1,\dots,a_m]$ and
$\widehat{\square}[b_1,\dots,b_n]$.  \epf

\section{Comparison of the two semantics of CCS}
\label{compare2}

\subsection*{Interpreting CCS as labelled precubical sets}

The \textit{CCS process names} are generated by the following syntax:
\[
P::=nil \ |\ a.P \ |\ (\nu a)P \ |\ P + P \ |\ P|| P \ |\
\rec(x)P(x)\] where $P(x)$ means a process name with one free variable
$x$. The variable $x$ must be \textit{guarded}, that is it must lie in
a prefix term $a.P'(x)$ for some $a\in\Sigma$.

\bd A labelled precubical set $\ell:K\rightarrow !\Sigma$ {\rm
  decorated by process names} is a labelled precubical set together
with a set map $d:K_0 \rightarrow \proc_\Sigma$ called the {\rm
  decoration}. \ed

Let us define by induction on the syntax of the CCS process name $P$
the decorated labelled precubical set $\square\llbracket P\rrbracket$
(see \cite{ccsprecub} for further explanations). The labelled
precubical set $\square\llbracket P\rrbracket$ has a unique initial
state canonically decorated by the process name $P$, and its other
states will be decorated in an inductive way. Therefore for every
process name $P$, $\square\llbracket P\rrbracket$ is an object of the
double comma category $\{i\}\ddownarrow \square^{op}\set \ddownarrow
!\Sigma$.  One has $\square\llbracket nil\rrbracket:=\square[0]$,
$\square\llbracket \mu.nil\rrbracket:=\mu.nil
\stackrel{(\mu)}\longrightarrow nil$, $\square\llbracket P+Q\rrbracket
:= \square\llbracket P\rrbracket \oplus \square\llbracket Q\rrbracket$
with the binary coproduct taken in $\{i\}\ddownarrow \square^{op}\set
\ddownarrow !\Sigma$, the pushout diagram of precubical sets
\[\xymatrix{
    \square[0]=\{0\} \ar@{->}[r]^-{0\mapsto nil} \ar@{->}[d]^-{0\mapsto P} & \square\llbracket \mu.nil\rrbracket \ar@{->}[d] \\
    \square\llbracket P\rrbracket \ar@{->}[r] & \cocartesien
    {\square\llbracket \mu.P\rrbracket},}\]
  the pullback diagram of precubical sets
\[\xymatrix{
    \square\llbracket (\nu a) P\rrbracket \ar@{->}[r] \ar@{->}[d] \cartesien & \square\llbracket P\rrbracket \ar@{->}[d] \\
    !(\Sigma\backslash \{a,\overline{a}\}) \ar@{->}[r] & !\Sigma,}
\]
the formula giving the interpretation of the parallel composition with
synchronization \[\square\llbracket P||Q\rrbracket :=
\square\llbracket P\rrbracket \ot_\Sigma \square\llbracket
Q\rrbracket\] and finally $\square\llbracket \rec(x)P(x)\rrbracket$
defined as the least fixed point of $P(-)$.

The prefix operator, the direct sum and the restriction operator are
$\omega$-continuous, that is to say they preserve the upper bounds of
ascending $\omega$-chains of labelled precubical sets $K^0 \subset K^1
\subset K^2 \subset \dots$, since they are finitely accessible and
since the upper bound is given by the colimit of the chain.  The
synchronized tensor product is also $\omega$-continuous since it is
colimit-preserving by \cite[Proposition~A.2]{ccsprecub}. Moreover, the
condition imposed on $P(x)$ implies that for all process names $Q_1$
and $Q_2$ with $\square \llbracket Q_1\rrbracket \subset \square
\llbracket Q_2\rrbracket$, one has $\square \llbracket
P(Q_1)\rrbracket \subset \square \llbracket P(Q_2)\rrbracket$.
Therefore the mapping $P(-)$ is $\omega$-continuous and
non-decreasing. Thus, the labelled precubical set \[\square\llbracket
\rec(x)P(x)\rrbracket:=\liminj\limits_n \square\llbracket
P^n(nil)\rrbracket \iso \bigcup_{n\geq 0} \square\llbracket
P^n(nil)\rrbracket\] will be equal to the least fixed point of
$P(-)$. This is a particular case of the Kleene fixed-point theorem on
a directed complete partial order.

\subsection*{Interpreting CCS as labelled transverse symmetric
  precubical sets}

Let us give now the new semantics of CCS in terms of labelled
transverse symmetric precubical sets.

\bd A labelled transverse symmetric precubical set $\ell:K\rightarrow
\sh\mathcal{L}(!\Sigma)$ {\rm decorated by process names} is a
labelled transverse symmetric precubical set together with a set map
$d:K_0 \rightarrow \proc_\Sigma$ called the {\rm decoration}. \ed

The interpretation of a CCS process name $P$ in terms of a decorated
labelled transverse symmetric precubical set
$\widehat{\square}\llbracket P\rrbracket$ is defined by induction on
the syntax of $P$, as for the case of labelled precubical sets.  The
only differences with the latter case are the pullback diagram
\[\xymatrix{
    \widehat{\square}\llbracket (\nu a) P\rrbracket \ar@{->}[r] \ar@{->}[d] \cartesien & \widehat{\square}\llbracket P\rrbracket \ar@{->}[d] \\
    \sh\mathcal{L}(!(\Sigma\backslash \{a,\overline{a}\})) \ar@{->}[r] &
    \sh\mathcal{L}(!\Sigma),}
\]
and the equation $\widehat{\square}\llbracket P||Q\rrbracket :=
\widehat{\square}\llbracket P\rrbracket \ot_\Sigma
\widehat{\square}\llbracket Q\rrbracket$ where $\ot_\Sigma$ is now the
synchronized tensor product of labelled transverse symmetric
precubical sets. Corollary~\ref{la} enables us to construct the least
fixed point of $P(-)$ in the same way as in the case of labelled
precubical sets.

\subsection*{The two semantics have same geometric realization}

\bth \label{LLL} For every CCS process name $P$, there is an
isomorphism of labelled transverse symmetric precubical sets
$\widehat{\square}\llbracket P\rrbracket \iso
\mathcal{L}(\square\llbracket P\rrbracket)$ and an isomorphism of
(labelled) flows $|\widehat{\square}\llbracket P\rrbracket| \iso
|\square\llbracket P\rrbracket|$. \eth

\bpf Let $K$ be a labelled precubical set. Let $a\in \Sigma\backslash
\{\tau\}$. Let $(\nu a)K$ be the labelled precubical set defined by
the pullback diagram
\[
\xymatrix{
(\nu a)K \cartesien \fr{\subset} \fd{} && K \fd{}\\
&&\\
!(\Sigma\backslash\{a,\overline{a}\}) \fr{} && !\Sigma.}
\]
One obtains the commutative diagram of labelled transverse symmetric
precubical sets
\[
\xymatrix{
\mathcal{L}((\nu a)K) \fr{\subset} \fd{} && \mathcal{L}(K) \fd{}\\
&&\\
\mathcal{L}(!(\Sigma\backslash\{a,\overline{a}\})) \fr{}\fd{} && \mathcal{L}(!\Sigma)\fd{}\\
&&\\
\sh\mathcal{L}(!(\Sigma\backslash\{a,\overline{a}\})) \fr{} && \sh\mathcal{L}(!\Sigma).}
\]
The map $(\nu a)K \rightarrow K$ is an inclusion of presheaves: the
labelled precubical set $(\nu a)K$ is the subobject of $K$ containing
the labelled cubes of $K$ not containing $a$ nor $\overline{a}$ as
label. By Proposition~\ref{carre}, the transverse symmetric precubical
set $\mathcal{L}((\nu a)K)$ is the subobject of $\mathcal{L}(K)$
containing the $p$-cubes of $\mathcal{L}(K)$ of the form $\mu^*(x)$
where $\mu : [p] \rightarrow [p]$ is a map of $\widehat{\square}$ and
$x$ is a $p$-cube of $(\nu a)K$. Therefore the map $\mathcal{L}((\nu
a)K) \rightarrow \mathcal{L}(K)$ is an inclusion of presheaves as
well. Consider now a commutative diagram of labelled transverse
symmetric precubical sets
\[
\xymatrix{
Z \ar@{-->}[rd]^-{k}\ar@/^10pt/[rrrd]^-{f}\ar@/_10pt/[rddd]&&&\\
&\mathcal{L}((\nu a)K)  \fr{\subset} \fd{} && \mathcal{L}(K)  \fd{}\\
&&&\\
&\sh\mathcal{L}(!(\Sigma\backslash\{a,\overline{a}\})) \fr{} && \sh\mathcal{L}(!\Sigma)}
\]
Every $p$-cube $x$ of $Z$ is taken to a $p$-cube $f(x)$ of
$\mathcal{L}(K)$. By Proposition~\ref{pintuitive}, $f(x) = \mu^*(y)$
for some $p$-cube $y\in K$ and for some map $\mu : [p] \rightarrow
[p]$ of $\widehat{\square}$. By construction, $y$ does not use the
labels $a$ or $\overline{a}$. Thus $y\in (\nu a)K$. Therefore $f(x)$
is a $p$-cube of $\mathcal{L}((\nu a)K)$. Hence $k$ exists and is
unique since the map $\mathcal{L}((\nu a)K) \rightarrow
\mathcal{L}(K)$ is an inclusion of presheaves.  Thus, the diagram of
labelled transverse symmetric precubical sets
\[
\xymatrix{
\mathcal{L}((\nu a)K)  \fr{\subset} \fd{} && \mathcal{L}(K) \fd{}\\
&&\\
\sh\mathcal{L}(!(\Sigma\backslash\{a,\overline{a}\})) \fr{} && \sh\mathcal{L}(!\Sigma)}
\]
is a pullback. So the isomorphism $\widehat{\square}\llbracket
P\rrbracket \iso \mathcal{L}(\square\llbracket P\rrbracket)$ implies
the isomorphism $\widehat{\square}\llbracket (\nu a)P\rrbracket \iso
\mathcal{L}(\square\llbracket (\nu a)P\rrbracket)$. Therefore, the
isomorphism of labelled transverse symmetric precubical sets
$\widehat{\square}\llbracket P\rrbracket \iso
\mathcal{L}(\square\llbracket P\rrbracket)$ is proved by induction on
the syntax of the process name $P$, using Theorem~\ref{cavavraiment2}
and the fact that the functor $\mathcal{L}$ preserves colimits since
it is a left adjoint. The isomorphism of labelled flows
$|\widehat{\square}\llbracket P\rrbracket| \iso |\square\llbracket
P\rrbracket|$ is a consequence of Proposition~\ref{factor_rea}. \epf

\part{Appendix}

\section{The case of labelled symmetric precubical sets}
\label{labsym}

By Theorem~\ref{identification_smallest}, the category of cubes
$\square_S$ is not shell-complete. It is interesting anyway for the
three following reasons. 1) It is possible to give an explicit
description of the symmetric precubical sets of labels with
Proposition~\ref{descr_explicit}. Such a description is still an open
problem for the transverse symmetric precubical set of labels (cf.
Conjecture~\ref{conj0}). 2) The category of cubes $\square_S$ is the
smallest category of cubes $\mathcal{A}$ such that the labelled cubes
$\mathcal{A}[a_{\sigma(1)}, \dots, a_{\sigma(n)}]$ with $a_1, \dots,
a_n \in \Sigma$ for $\sigma$ running over the set of permutations of
$\{1, \dots, n\}$ belong to the same isomorphism class. Let us recall
that the labelled precubical sets $\square[a_{\sigma(1)}, \dots,
a_{\sigma(n)}]$ and $\square[a_{\sigma'(1)}, \dots, a_{\sigma'(n)}]$
are not isomorphic as soon as $(a_{\sigma(1)}, \dots, a_{\sigma(n)})
\neq (a_{\sigma'(1)}, \dots, a_{\sigma'(n)})$.  3) There is a strong
link between labelled symmetric precubical sets and higher dimensional
transition systems in the sense of Cattani and Sassone
\cite{MR1461821}, see \cite{hdts}.  Indeed, it turns out that the
category of higher dimensional transition systems in the sense of
Cattani and Sassone is equivalent to a full reflective subcategory of
that of labelled symmetric precubical sets.

\subsection*{Description of the symmetric precubical set of labels}

The following combinatorial lemma is well-known (see \cite{poset_tool}
for a survey).

\begin{lem} \label{perm} Let $p\geq 1$. The group of automorphisms of
  the poset $[p]$ is isomorphic to the symmetric group on
  $\{1,...,p\}$.  In other terms, let $f$ be an automorphism of the
  poset $[p]$. Then there exists a permutation $\pi$ of the set
  $\{1,\dots,p\}$ such that $f(\epsilon_1, \dots, \epsilon_p) =
  (\epsilon_{\pi(1)}, \dots, \epsilon_{\pi(p)})$.
\end{lem}

\bpf Let $I\subset \{1,\dots,p\}$. Let $e_I$ be the element
$(\epsilon_1,\dots,\epsilon_p)$ of $[p]$ such that $\epsilon_i=1$ if
and only if $i\in I$. Since $f$ is bijective and strictly increasing,
it preserves the distance of Proposition~\ref{distance}. The distance
between $e_{\varnothing}$ and $f(e_{\{i\}})$ is $1$. So there exists a
permutation $\pi$ of $\{1,\dots,p\}$ such that $f(e_{\{i\}}) =
e_{\{\pi(i)\}}$.  Let $g(\epsilon_1,\dots,\epsilon_p) =
(\epsilon_{\pi^{-1}(1)},\dots,\epsilon_{\pi^{-1}(p)})$. Then
$g(f(e_{\{i\}})) = g(e_{\{\pi(i)\}}) = e_{\{i\}}$. It then suffices to
prove by induction on the cardinality $c$ of $I$ that $g(f(e_{I})) =
e_{I}$. Let $c\geq 2$ with $c\leq p$. Assume that $g(f(e_{J})) =
e_{J}$ for all subsets $J$ of $\{1,\dots,p\}$ of cardinality
$c-1$. Let $I$ be a subset of $\{1,\dots,p\}$ of cardinality $c$. Then
the distance between $g(f(e_{I\backslash\{i\}}))$ and $g(f(e_{I}))$ is
$1$ for all $i\in I$. By induction hypothesis, one has
$g(f(e_{I\backslash\{i\}})) = e_{I\backslash\{i\}}$. So the only
possibility is $g(f(e_{I})) = e_{I}$.  \epf

\bp \label{inj} An adjacency-preserving map $f:[m] \rightarrow [n]$
belongs to $\square_S$ if and only if $f$ is one-to-one. \ep

\bpf It is clear that any map of $\square_S$ is
one-to-one. Conversely, let $f:[m] \rightarrow [n]$ be a one-to-one
adjacency-preserving map. Then by
Proposition~\ref{decomposition_distance}, $f$ factors uniquely as a
composite $[m] \stackrel{\psi}\longrightarrow [m]
\stackrel{\phi}\longrightarrow [n]$ with $\phi\in \square$ and $\psi$
adjacency-preserving one-to-one. A cardinality argument implies that
$\psi$ is a bijection. Therefore $f\in \square_S$ by
Lemma~\ref{perm}. \epf

As for precubical sets, let $\de_i^\alpha = (\delta_i^\alpha)^*$. And
let $s_i = (\sigma_i)^*$.

\bp \label{descr_explicit} The symmetric precubical set of labels
$\sh_{\square_S}\mathcal{L}_{\square_S}(!\Sigma)$ is isomorphic to the
following symmetric precubical set, denoted by $!^S\Sigma$:
\begin{itemize}
\item $(!^S\Sigma)_0=\{()\}$ (the empty word)
\item for $n\geq 1$, $(!^S\Sigma)_n=\Sigma^n$
\item $\de_i^0(a_1,\dots,a_n) = \de_i^1(a_1,\dots,a_n) =
  (a_1,\dots,\widehat{a_i},\dots,a_n)$ where the notation
  $\widehat{a_i}$ means that $a_i$ is removed.
\item $s_i(a_1,\dots,a_n) =
  (a_1,\dots,a_{i-1},a_{i+1},a_i,a_{i+2},\dots,a_n)$ for $1\leq i\leq n$.
\end{itemize}
\ep

\bpf The category of cubes $\square_S$ is the small category freely
generated by the $\delta_i^\alpha$ and $\sigma_i$ operators and by the
cocubical relations, the algebraic relations of
Proposition~\ref{sigma_delta}, and the Moore relations for symmetry
operators $\sigma_i\sigma_i=\id$, $\sigma_i\sigma_j\sigma_i =
\sigma_j\sigma_i \sigma_j$ for $i=j-1$ and
$\sigma_i\sigma_j=\sigma_j\sigma_i$ for $i<j-1$ by
\cite[Theorem~8.1]{MR1988396}. 

It is easy to prove that the $s_i$ and $\de_i^\alpha$ operators of
$!^S\Sigma$ satisfy the dual of these algebraic relations. So
$!^S\Sigma$ together with the $\de_i^\alpha$ and $s_i$ operators is a
well-defined symmetric precubical set. 

The identity of $!\Sigma$ yields a map of precubical sets $!\Sigma
\rightarrow \omega_{\square_S}(!^S\Sigma)$. Hence by adjunction, one
obtains a map $\mathcal{L}_{\square_S}(!\Sigma) \rightarrow
!^S\Sigma$. The symmetric precubical set $!^S\Sigma$ is orthogonal to
the set of morphisms
$\{\square_S[p]\sqcup_{\de\square_S[p]}\square_S[p] \rightarrow
\square_S[p], p\geq 2\}$ for the same reason as $!\Sigma$ is
orthogonal to the set of morphisms
$\{\square[p]\sqcup_{\de\square[p]}\square[p] \rightarrow \square[p],
p\geq 2\}$. Hence by adjunction, one obtains a map of symmetric
precubical sets $f : \sh_{\square_S}\mathcal{L}_{\square_S}(!\Sigma)
\rightarrow !^S\Sigma$ which is clearly onto: an inverse image of
$(a_1,\dots,a_n)\in (!^S\Sigma)_n$ for $n\geq 1$ is given by the image
of $(a_1,\dots,a_n)\in \mathcal{L}_{\square_S}(!\Sigma)_n$ by the
canonical map $\mathcal{L}_{\square_S}(!\Sigma) \rightarrow
\sh_{\square_S}\mathcal{L}_{\square_S}(!\Sigma)$. 

Let us prove by induction on $p\geq 1$ that the map $f_{\leq p} :
(\sh_{\square_S}\mathcal{L}_{\square_S}(!\Sigma))_{\leq p} \rightarrow
(!^S\Sigma)_{\leq p}$ is one-to-one. The map induces the isomorphism
$f_{\leq 1} : (\sh_{\square_S}\mathcal{L}_{\square_S}(!\Sigma))_{\leq
  1} \rightarrow (!^S\Sigma)_{\leq 1}$ by
Proposition~\ref{restrict1_1} and Proposition~\ref{restrict1_2}. Hence
the proof is complete for $p=1$. Let us suppose that the map $f_{\leq
  p} : (\sh_{\square_S}\mathcal{L}_{\square_S}(!\Sigma))_{\leq p}
\rightarrow (!^S\Sigma)_{\leq p}$ is an isomorphism for $p \geq
1$. Let $x,y \in
(\sh_{\square_S}\mathcal{L}_{\square_S}(!\Sigma))_{p+1}$ be two
$(p+1)$-cubes having the same image in $!^S\Sigma$. Then they have the
same boundary in $(!^S\Sigma)_{\leq p}$, and therefore $x$ and $y$
have the same boundary $\de x = \de y$ by induction hypothesis. One
obtains a commutative square of solid arrows
\[
\xymatrix{ \square_S[p+1]\sqcup_{\de\square_S[p+1]}\square_S[p+1]
  \fr{x\sqcup_{\de x}y}\fd{} &&
  \sh_{\square_S}\mathcal{L}_{\square_S}(!\Sigma)\fd{} \\
  &&\\
  \square_S[p+1] \ar@{-->}[rruu]^-{k}\fr{} && \mathbf{1}.}
\] 
The lift $k$ exists and is unique. So $x = k = y$. The induction
hypothesis is therefore proved for $p+1$.  \epf

\subsection*{The labelled directed symmetric coskeleton construction}

The following proposition is similar to Proposition~\ref{pintuitive}.

\bp \label{pintuitive1} Let $K$ be a precubical set. For any $p$-cube
$x$ of $\omega_{\square_S}\mathcal{L}_{\square_S}(K)$ with $p\geq 0$,
there exists a $p$-cube $y$ of $K \subset
\omega_{\square_S}\mathcal{L}_{\square_S}(K)$ and a map $\mu \in
\square_S([p],[p])$ such that $x = \mu^*(y)$ where
$\mu^*:\mathcal{L}_{\square_S}(K)_p \rightarrow
\mathcal{L}_{\square_S}(K)_p$ is the image of $\mu$ by the presheaf
$\mathcal{L}_{\square_S}(K) \in \square_S^{op}\set$.  \ep

Note that as in Proposition~\ref{pintuitive}, the decomposition is
actually unique.

\bpf With the notations of the proof of
Proposition~\ref{pintuitive}. By
Proposition~\ref{decomposition_distance} and Proposition~\ref{inj},
the set map $\overline{x}:[p] \rightarrow [n]$ factors as a composite
$[p] \stackrel{\mu}\longrightarrow [p] \stackrel{\phi}\longrightarrow
[n]$ with $\mu \in \square_S$ and $\phi \in \square$. \epf

By Proposition~\ref{i1_relatif}, the truncation functor
\[\square_S^{op}\set \ddownarrow
\sh_{\square_S}\mathcal{L}_{\square_S}(!\Sigma)\rightarrow
(\square_S)_{n}^{op}\set\ddownarrow
\sh_{\square_S}\mathcal{L}_{\square_S}(!\Sigma)\] has a right
adjoint \[\cosk_n^{\square_S,\Sigma}:(\square_S)_{n}^{op}\set\ddownarrow
\sh_{\square_S}\mathcal{L}_{\square_S}(!\Sigma) \rightarrow
\square_S^{op}\set\ddownarrow
\sh_{\square_S}\mathcal{L}_{\square_S}(!\Sigma).\]

\bd (Compare with Definition~\ref{def_directed}) Let $K$ be a
$1$-dimensional labelled symmetric precubical set with $K_0 = [p]$ for
some $p\geq 0$. The {\rm labelled symmetric directed coskeleton} of
$K$ is the labelled precubical set $\COSK_S^\Sigma(K)$ defined as the
subobject of $\cosk^{\square_S,\Sigma}_1(K)$ such that:
\begin{itemize}
\item $\COSK_S^\Sigma(K)_{\leq 1} = \cosk^{\square_S,\Sigma}_1(K)_{\leq 1}$
\item for every $n\geq 2$, $x\in \cosk^{\square_S,\Sigma}_1(K)_n$ is an $n$-cube
  of $\COSK_S^\Sigma(K)$ if and only if the set map
  $x_0:[n]\rightarrow [p]$ is {\rm non-twisted}, i.e. $x_0: [n]
  \rightarrow [p]$ is a composite\footnote{The factorization is
    necessarily unique.}
\[x_0: [n] \stackrel{\phi}\longrightarrow [q]
\stackrel{\psi}\longrightarrow [p],\] where $\psi$ is a morphism of
the small category $\square$ and where $\phi$ is of the form
\[(\epsilon_1,\dots,\epsilon_{n}) \mapsto
(\epsilon_{i_1},\dots,\epsilon_{i_q})\] such that
$\{1,\dots,n\}\subset \{i_1,\dots,i_q\}$.
\end{itemize} \ed

The link with labelled precubical sets is:

\bp Let $K$ be a $1$-dimensional labelled (symmetric) precubical set
with $K_0 = [p]$ for some $p\geq 0$. Then there is the isomorphism of
labelled symmetric precubical sets 
\[\mathcal{L}_{\square_S}(\COSK^\Sigma(K)) \iso \COSK_S^\Sigma(K).\] 
\ep

\bpf By a proof similar to the one of
Theorem~\ref{inclusion_full_cosk}, one obtains the inclusion of
presheaves
\[\mathcal{L}_{\square_S}(\COSK^\Sigma(K)) \subset
\cosk_1^{\square_S,\Sigma}(K).\]
It is clear that the inclusion above factors as the composite of
inclusions
\[\mathcal{L}_{\square_S}(\COSK^\Sigma(K)) \subset \COSK_S^\Sigma(K)
\subset \cosk_1^{\square_S,\Sigma}(K).\] The left-hand inclusion is an
equality by Proposition~\ref{pintuitive1}.  \epf

Thanks to Proposition~\ref{factor_rea}, one obtains the isomorphism of
flows \[|\COSK^\Sigma(K)| \iso |\COSK_S^\Sigma(K)|.\]

\subsection*{Interpreting CCS as labelled symmetric precubical sets}

\bd Let $K$ and $L$ be two labelled symmetric precubical
sets. The {\rm tensor product with synchronization} (or {\rm
  synchronized tensor product}) of $K$ and $L$ is
\[K \ot_\Sigma L := \liminj_{\square_S[m]\rightarrow K}
\liminj_{\square_S[n]\rightarrow L}
\COSK_S^{\Sigma}(\square_S[m]_{\leq 1} \p_\Sigma
\square_S[n]_{\leq 1}).\] \ed

One can then easily adapt the semantics of CCS to the case of labelled
symmetric precubical sets. The interest of this setting is that it is
simpler than the one of transverse symmetric precubical sets, and
that, as in Proposition~\ref{commutativite}, there is an isomorphism
of labelled symmetric precubical sets $K\ot_\Sigma L \iso L\ot_\Sigma
K$ for all labelled symmetric precubical sets $K$ and $L$. The
synchronized tensor product of symmetric precubical sets is also
colimit-preserving by Proposition~\ref{preserve_colim} and therefore
associative.


\begin{thebibliography}{AABS02}

\bibitem[AA89]{phd-Al-Agl}
Fahd A.~A. Al-Agl.
\newblock {\em Aspects of multiple categories}.
\newblock PhD thesis, University of Wales, Department of Pure Mathematics,
  University College of North Wales, Bangor, Gwynedd LL57 1UT, U.K., September
  1989.

\bibitem[AABS02]{MR1929304}
F.~A. Al-Agl, R.~Brown, and R.~Steiner.
\newblock Multiple categories: the equivalence of a globular and a cubical
  approach.
\newblock {\em Adv. Math.}, 170(1):71--118, 2002.

\bibitem[AR94]{MR95j:18001}
J.~Ad{\'a}mek and J.~Rosick{\'y}.
\newblock {\em Locally presentable and accessible categories}.
\newblock Cambridge University Press, Cambridge, 1994.

\bibitem[Bek00]{MR1780498}
T.~Beke.
\newblock Sheafifiable homotopy model categories.
\newblock {\em Math. Proc. Cambridge Philos. Soc.}, 129(3):447--475, 2000.

\bibitem[BH81]{Brown_cube}
R.~Brown and P.~J. Higgins.
\newblock On the algebra of cubes.
\newblock {\em J. Pure Appl. Algebra}, 21(3):233--260, 1981.

\bibitem[BHR84]{0628.68025}
S.~D. Brookes, C.~A.~R. Hoare, and A.~W. Roscoe.
\newblock {A theory of communicating sequential processes.}
\newblock {\em J. Assoc. Comput. Mach.}, 31:560--599, 1984.

\bibitem[Bro06]{MR2273730}
R.~Brown.
\newblock {\em Topology and groupoids}.
\newblock BookSurge, LLC, Charleston, SC, 2006.
\newblock Third edition of {\it Elements of modern topology} [McGraw-Hill, New
  York, 1968].

\bibitem[CS96]{MR1461821}
G.~L. Cattani and V.~Sassone.
\newblock Higher-dimensional transition systems.
\newblock In {\em 11th Annual IEEE Symposium on Logic in Computer Science (New
  Brunswick, NJ, 1996)}, pages 55--62. IEEE Comput. Soc. Press, Los Alamitos,
  CA, 1996.

\bibitem[CS02]{monographie_hocolim}
W.~Chach{\'o}lski and J.~Scherer.
\newblock Homotopy theory of diagrams.
\newblock {\em Mem. Amer. Math. Soc.}, 155(736):x+90, 2002.

\bibitem[DS95]{MR1361887}
W.~G. Dwyer and J.~Spali{\'n}ski.
\newblock Homotopy theories and model categories.
\newblock In {\em Handbook of algebraic topology}, pages 73--126.
  North-Holland, Amsterdam, 1995.

\bibitem[Fah05a]{fahrenberg05-hda-long}
U.~Fahrenberg.
\newblock Bisimulation for higher-dimensional automata. {A} geometric
  interpretation.
\newblock Research report R-2005-01, Department of Mathematical Sciences,
  Aalborg University, 2005.
\newblock {http://www.math.aau.dk/research/reports/R-2005-01.ps}. Extended
  version of \cite{fahrenberg05-hda}.

\bibitem[Fah05b]{fahrenberg05-hda}
U.~Fahrenberg.
\newblock A category of higher-dimensional automata.
\newblock In {\em Proc.\ FOSSACS'05}, volume 3441 of {\em Lecture Notes in
  Computer Science}, pages 187--201. Springer-Verlag, 2005.

\bibitem[FGR98]{MR1683333}
L.~Fajstrup, E.~Goubault, and M.~Rau{\ss}en.
\newblock Detecting deadlocks in concurrent systems.
\newblock In {\em CONCUR'98: concurrency theory (Nice)}, volume 1466 of {\em
  Lecture Notes in Comput. Sci.}, pages 332--347. Springer, Berlin, 1998.

\bibitem[FR08]{FR}
L.~Fajstrup and J.~Rosick\'y.
\newblock A convenient category for directed homotopy.
\newblock {\em Theory and Applications of Categories}, 21(1):pp 7--20, 2008.

\bibitem[Gau00]{homcat}
P.~Gaucher.
\newblock Homotopy invariants of higher dimensional categories and concurrency
  in computer science.
\newblock {\em Math. Structures Comput. Sci.}, 10(4):481--524, 2000.

\bibitem[Gau01]{coin}
P.~Gaucher.
\newblock Combinatorics of branchings in higher dimensional automata.
\newblock {\em Theory Appl. Categ.}, 8(12):324--376 (electronic), 2001.

\bibitem[Gau03]{model3}
P.~Gaucher.
\newblock A model category for the homotopy theory of concurrency.
\newblock {\em Homology, Homotopy and Applications}, 5(1):p.549--599, 2003.

\bibitem[Gau08]{ccsprecub}
P.~Gaucher.
\newblock Towards a homotopy theory of process algebra.
\newblock {\em Homology Homotopy Appl.}, 10(1):353--388 (electronic), 2008.

\bibitem[Gau09]{hdts}
P.~Gaucher.
\newblock Directed algebraic topology and higher dimensional transition
  systems.
\newblock Preprint, 2009.

\bibitem[GG03]{diCW}
P.~Gaucher and E.~Goubault.
\newblock Topological deformation of higher dimensional automata.
\newblock {\em Homology, Homotopy and Applications}, 5(2):39--82, 2003.

\bibitem[GM03]{MR1988396}
M.~Grandis and L.~Mauri.
\newblock Cubical sets and their site.
\newblock {\em Theory Appl. Categ.}, 11(8):185--211 (electronic), 2003.

\bibitem[Gou02]{labelled}
E.~Goubault.
\newblock Labelled cubical sets and asynchronous transistion systems: an
  adjunction.
\newblock Presented at CMCIM'02, 2002.

\bibitem[Gou03]{survol}
E.~Goubault.
\newblock Some geometric perspectives in concurrency theory.
\newblock {\em Homology, Homotopy and Applications}, 5(2):95--136, 2003.

\bibitem[Gra03]{mg}
M.~Grandis.
\newblock Directed homotopy theory. {I}.
\newblock {\em Cah. Topol. G\'eom. Diff\'er. Cat\'eg.}, 44(4):281--316, 2003.

\bibitem[Hir03]{ref_model2}
P.~S. Hirschhorn.
\newblock {\em Model categories and their localizations}, volume~99 of {\em
  Mathematical Surveys and Monographs}.
\newblock American Mathematical Society, Providence, RI, 2003.

\bibitem[Hov99]{MR99h:55031}
M.~Hovey.
\newblock {\em Model categories}.
\newblock American Mathematical Society, Providence, RI, 1999.

\bibitem[HT09]{HT}
F.~Hivert and N.~M. Thi{\'e}ry.
\newblock The {H}ecke group algebra of a {C}oxeter group and its representation
  theory.
\newblock {\em J. Algebra}, 321(8):2230--2258, 2009.

\bibitem[Kri08]{SK}
S.~Krishnan.
\newblock A convenient category of locally preordered spaces.
\newblock {\em Applied Categorical Structures}, 17(5):1--22, 2008.
\newblock doi:10.1007/s10485-008-9140-9.

\bibitem[Lew78]{Ref_wH}
L.~G. Lewis.
\newblock {\em The stable category and generalized Thom spectra}.
\newblock PhD thesis, University of Chicago, 1978.

\bibitem[May99]{MR2000h:55002}
J.~P. May.
\newblock {\em A concise course in algebraic topology}.
\newblock University of Chicago Press, Chicago, IL, 1999.

\bibitem[Mil89]{0683.68008}
R.~Milner.
\newblock {\em {Communication and concurrency.}}
\newblock {Prentice Hall International Series in Computer Science. New York
  etc.: Prentice Hall. XI, 260 p. }, 1989.

\bibitem[ML98]{MR1712872}
S.~Mac~Lane.
\newblock {\em Categories for the working mathematician}.
\newblock Springer-Verlag, New York, second edition, 1998.

\bibitem[MLM94]{MR1300636}
S.~Mac~Lane and I.~Moerdijk.
\newblock {\em Sheaves in geometry and logic}.
\newblock Universitext. Springer-Verlag, New York, 1994.
\newblock A first introduction to topos theory, Corrected reprint of the 1992
  edition.

\bibitem[Pra91]{Pratt}
V.~Pratt.
\newblock Modeling concurrency with geometry.
\newblock In ACM Press, editor, {\em Proc. of the 18th ACM Symposium on
  Principles of Programming Languages}, 1991.

\bibitem[Tho79]{MR510404}
R.~W. Thomason.
\newblock Homotopy colimits in the category of small categories.
\newblock {\em Math. Proc. Cambridge Philos. Soc.}, 85(1):91--109, 1979.

\bibitem[vG06]{rvg}
R.J. van Glabbeek.
\newblock {On the expressiveness of higher dimensional automata.}
\newblock {\em Theor. Comput. Sci.}, 356(3):265--290, 2006.

\bibitem[Wac07]{poset_tool}
M.~L. Wachs.
\newblock Poset topology: Tools and applications.
\newblock In {\em Miller, Ezra (ed.) et al., Geometric combinatorics.
  Providence, RI: American Mathematical Society (AMS); Princeton, NJ: Institute
  for Advanced Studies. IAS/Park City Mathematics Series 13}, pages 497--615.
  2007.

\bibitem[WN95]{MR1365754}
G.~Winskel and M.~Nielsen.
\newblock Models for concurrency.
\newblock volume~4 of {\em Handb. Log. Comput. Sci.}, pages 1--148. Oxford
  Univ. Press, New York, 1995.

\bibitem[Wor04]{exHDA}
K.~Worytkiewicz.
\newblock Synchronization from a categorical perspective.
\newblock ArXiv cs.PL/0411001, 2004.

\end{thebibliography}
\end{document}